\documentclass[10pt]{amsart}

\usepackage{amssymb,amsmath,amscd,amsfonts,a4,psfrag,graphicx,
latexsym,epsfig,color}
\usepackage[active]{srcltx}
\usepackage[all]{xy}

\everymath{\displaystyle}

\newtheorem{teo}{Theorem}[section]
\newtheorem{lem}[teo]{Lemma}
\newtheorem{cor}[teo]{Corollary}
\newtheorem{exa}[teo]{Example}
\newtheorem{prop}[teo]{Proposition}
\newtheorem{defi}[teo]{Definition}

\newtheorem{remark}[teo]{Remark}
\newtheorem{remarks}[teo]{Remarks}

\newtheorem{prob}[teo]{Problem}

\newcommand{\fonc}[5]{                     % fonction
            \begin{array}{lcll}#1 :& #2 & \longrightarrow & #3 \\   %
                         &#4 &\longmapsto & #5          %
            \end{array}}

\newcommand{\cvd}{\hfill$\Box$}

\newcommand{\mr}{\mathbb{R}}
\newcommand{\mc}{\mathbb{C}}
\newcommand{\mz}{\mathbb{Z}}
\newcommand{\mh}{\mathbb{H}}

\newcommand{\mn}{\mathbb{N}}

\newcommand{\C}{\mathbb{C}}
\newcommand{\Z}{\mathbb{Z}}
\newcommand{\N}{\mathbb{N}}

\newcommand{\bw}{{\bf w}}
\newcommand{\bu}{{\bf u}}
\newcommand{\bv}{{\bf v}}

\newcommand{\Aa}{{\mathcal A}}
\newcommand{\Bb}{{\mathcal B}}

\newcommand{\Cc}{\mathbb{C}_*}

\newcommand{\Dd}{{\mathcal D}}

\newcommand{\Ff}{{\mathcal F}}

\newcommand{\Hh}{{\mathcal H}}
\newcommand{\Ii}{{\mathcal I}}

\newcommand{\Ll}{{\mathcal L}}

\newcommand{\Nn}{{\mathcal N}}

\newcommand{\Pp}{{\mathcal P}}
\newcommand{\Rr}{{\mathcal R}}
\newcommand{\Ss}{{\mathcal S}}
\newcommand{\Tt}{{\mathcal T}}

\newcommand{\Ww}{{\mathcal W}}

\newcommand{\Qq}{{\mathcal Q}}

\newcommand{\CG}{{\mathfrak G}}

\newcommand{\hG}{{\mathfrak h}}
\newcommand{\lra}{\longrightarrow}

\newcommand{\ra}{\rightarrow}
\newcommand{\Dim}{{\it Proof.\ }}

\def\cvd{\hfill$\Box$}

\title[Analytic families of quantum hyperbolic invariants]{Analytic
  families of quantum hyperbolic invariants}

\author{St\'ephane Baseilhac$^1$, Riccardo Benedetti$^2$}%, Charles Frohman$^3$}

\begin{document}
%\today

\maketitle

%\vspace{0.6cm}

\noindent $^1$ Institut de Math\'ematiques et de Mod\'elisation, 
Universit\'e Montpellier 2, Case Courrier 51,
34095 Montpellier Cedex 5, France (sbaseilh@math.univ-montp2.fr)
\smallskip

\noindent $^2$ Dipartimento di Matematica, Universit\`a di Pisa, Largo
Bruno Pontecorvo 5, 56127 Pisa, Italy (benedett@dm.unipi.it)

\begin{abstract}  
  We organize the quantum hyperbolic invariants (QHI) of $3$-manifolds
  into sequences of rational functions indexed by the odd integers
  $N\geq 3$ and defined on moduli spaces of geometric structures
  refining the character varieties. In the case of one-cusped hyperbolic
  $3$-manifolds $M$ we generalize the QHI and get rational
  functions $\Hh_N^{h_f,h_c,k_c}$ depending on a finite set of
  cohomological data $(h_f,h_c,k_c)$ called {\it weights}. These functions are regular on a determined Abelian covering of degree $N^2$ of a Zariski open subset, canonically associated to $M$, of the geometric
  component of the variety of augmented $PSL(2,\mc)$-characters of
  $M$. New combinatorial ingredients are a weak version of branchings which exists on every triangulation, and state sums over weakly branched triangulations, including a sign correction which eventually fixes the sign ambiguity of the QHI. We describe in detail the invariants of three cusped manifolds, and present the results of numerical computations showing that the functions $\Hh_N^{h_f,h_c,k_c}$ depend on the weights as $N\rightarrow + \infty$, and recover the volume for some specific choices of the weights.  %We provide also a cohomological formula of the state sum symmetrization factor, and this eventually leads to a factorization of the QHI into reduced invariants. 
\end{abstract}

\tableofcontents

 \section{Introduction}\label{I-INTRO} 
 
 In the series of papers \cite{Top, GT,AGT} we defined a family of complex valued quantum invariants $\Hh_N(\Pp)$ of ``patterns'' $\Pp$ of geometric nature, called {\it quantum hyperbolic invariants} (QHI) and indexed by the odd integers $N\geq 3$. Roughly, a pattern consists of a compact oriented $3$-manifold equipped with a representation of the fundamental group into $PSL(2,\mc)$, plus some compatible cohomological datas (see the next subsections for a complete definition).  The QHI generalize the Kashaev invariants of links in $S^3$ \cite{LINK}, subject to the celebrated Volume Conjecture. They define a $(2+1)$-dimensional quantum field theory, and, when applied to mapping cylinders of surface diffeomorphisms, they coincide with the invariants derived from the local version of finite dimensional quantum Teichm\"uller theory (see \cite{Bai}, \cite{B+} and \cite{AGT}).   
 \smallskip

In \cite{Top, GT,AGT}, the invariance of the QHI was proved only up to sign and multiplication by $N$-th roots of unity. The eventual existence and the meaning of such a phase anomaly, as well as the determination of the asymptotical behaviour of $\Hh_N(\Pp)$ as $N\to +\infty$, are two main open issues of the theory. In order to tackle them, it seemed necessary to develop a functional approach that would clarify the intrinsic nature of the various combinatorial and geometric ingredients involved in the definition of the QHI. To this aim, we achieve in this paper the following goals: 

\smallskip

{\bf QHI of cusped manifolds:} We extend the QHI of any one-cusped hyperbolic $3$-manifold $M$ to invariants defined on the geometric component $X_0(M)$ of the variety of augmented
$PSL(2,\mc)$-characters of $M$. In \cite{GT, AGT} the QHI of cusped hyperbolic manifolds were defined only at the hyperbolic holonomy.

\smallskip

{\bf State sums over weakly branched triangulations:} In order to achieve the previous goal it
is necessary to introduce state sums over triangulations that do not support any branching, replaced by relaxed structures called {\it weak branchings}. These states sums give rise to QHI generalizing all the previously defined ones.

\smallskip

{\bf Analytic families of QHI:} We recast all the QHI into sequences $\{\Aa_N(Y)\}$ of families of complex
analytic spaces and maps, indexed by the odd integers $N\geq 3$. Each family $\Aa_N(Y)$ is associated to a {\it topological support} $Y$,
and provides concrete models of geometric structures over $Y$ called
{\it patterns}. Patterns over cusped manifolds have an
intrinsic meaning in terms of the $PSL(2,\mc)$-version of the
$A$-polynomial, and natural relationships with Chern-Simons theory. 
\smallskip
 
{\bf Fixing the sign ambiguity:} We include a sign correction in the state sum formulas of the QHI, which removes their sign ambiguity under a mild assumption on the ``bulk $c$-weight'' (see below). The sign correction depends on the combinatorics of the weak branching and is a by-product of \cite{BP2}. It becomes trivial when dealing (when possible) with branched triangulations, so that the QHI of \cite{Top,GT,AGT} are eventually defined up to multiplication by $N$-th roots of unity.

\smallskip

In a sequel to this paper (in collaboration with C. Frohman), we will develop our approach concerning the asymptotical behaviour of the QHI. Basically, given a topological support $Y$ and a sequence $\{\Pp_N\}$ of patterns over $Y$, we study the limit of the family $\{\Aa_N(Y)\}$ as $N\rightarrow +\infty$, instead of $\{\Hh_N(\Pp_N)\}$ for a single sequence $\{\Pp_N\}$. Then, assuming that $\textstyle \Hh_\infty(\{\Pp_N\}):=  \limsup_{N\rightarrow \infty} (\log|\Hh_N(\Pp_N)|/N)$ is finite (this is the case for all natural sequences $\{\Pp_N\}$), we consider the following problems: 
\begin{enumerate}
\item {\it Determine the nature (regularity) of $\Hh_\infty(\{\Pp_N\})$ as a
    function of the patterns over $Y$. }

\item {\it Describe the asymptotical behaviour of $\{\Hh_N(\Pp_N)\}$ as
    $N\rightarrow \infty$ in terms of classical geometric invariants
    of $Y$: Chern-Simons invariants, torsions, twisted cohomology,
    etc.}
\end{enumerate}

The Kashaev-Murakami-Murakami Volume Conjecture is a particular case of (2), for constant sequences of patterns associated to links in $S^3$ (see Theorem \ref{QH=J} below).
\smallskip

In the rest of this Introduction we describe with more details the content of the paper.

\subsection{QHI of cusped manifold patterns}
\label{I-QHI-cusped} 
In this paper we call {\it cusped manifold} an oriented, connected,
non--compact complete hyperbolic $3$-manifold of finite volume with
{\it one} cusp. Hence a cusped manifold $M$ is diffeomophic to the
interior of a compact $3$-manifold denoted by $V$, with one torus
boundary component.

A {\it pattern $\Pp=(Y_\Pp,\rho,(h,k))$ over $M$} consists of a
{\it topological support} $Y_\Pp$ together with additional geometric
structures determined by the couple $(\rho,(h,k))$.  The topological
support takes the form $Y_\Pp=(V,(h_c,k_c))$, where $(h_c,k_c)\in H^1(V;\Z/2\Z)\times H^1(\partial V; \Z)$ is a so-called {\it c-weight}, defined by a ``bulk $c$-weight'' $h_c$ and a ``boundary $c$-weight'' $k_c$ satisfying
\begin{equation}\label{cconstraint}
 r(k_c)=i^*(h_c)
 \end{equation}
 where $r: H^1(\partial V; \Z)\rightarrow H^1(\partial V;\Z/2\Z)$ is
 the reduction mod$(2)$, and the map $i^*: H^1(V; \Z/2\mz)
 \rightarrow H^1(\partial V; \Z/2\mz)$ is induced by the inclusion map
 $i:\partial V \to V$. The pattern $\Pp$ is obtained by completing
 $Y_\Pp$ with a couple $(\rho,(h_f,k_f))$ where $\rho$ is a $PSL(2,\C)$-{\it character of $V$}, i.e. a conjugacy
class of representations of $\pi_1(V)$ in $PSL(2,\C)$, and
$(h_f,k_f)\in H^1(V;\Z/2\Z) \times H^1(\partial V ; \C)$ is a so-called
{\it $f$-weight} (relative to $\rho$), defined by a ``bulk
$f$-weight'' $h_f$ and a ``boundary $f$-weight'' $k_f$ satisfying the
following constraint.  Up to conjugacy the restriction of $\rho$ to
the torus $\partial V$ is valued in the group of complex affine
transformations of the plane; the linear part of this restriction
defines a class in $H^1(\partial V; \mc^*)$.  Let $d\in H^1(\partial
V; \C/2i\pi\mz)$ be the log of this class, with imaginary part in
$]-\pi,\pi]$. One requires that for all $a\in H_1(\partial V; \Z)$,
\begin{equation}\label{fconstraint}
k_f(a)= d(a) \ {\rm mod}(i\pi); \ \ 
(k_f(a)-d(a))/i\pi = i^*(h_f)(a) \ \ {\rm mod} (2) . 
\end{equation}
Collecting the bulk and boundary weights we will often write $\Pp$ as
$(V,\rho,(h,k))$, where
$$(h,k)=((h_f,h_c),(k_f,k_c)).$$

{\bf Notation.} For every $n\in \N$, we write ``$a\equiv_n b$'' to mean
that $a$ and $b$ are equal up to multiplication by a power of
$\exp(2i\pi/n)$. If $n$ is odd, then $a\equiv_{2n} b$ if and only if
$a\equiv_n \pm b$.  We denote by $\mu_{n}$ the group of $n$-th roots of
unity, acting on $\mc$ by multiplication.
\smallskip

In \cite{GT, AGT}, for every cusped manifold $M$ and odd
$N\geq 3$ we defined quantum hyperbolic invariants
\begin{equation}\label{QHIcuspedpattern}
\Hh_N(M) := \Hh_N(V,\rho_{hyp},(0,0)) \ \in\mc/\mu_{2N}
\end{equation}
that is, for the pattern $\Pp=(V,\rho_{hyp},(0,0))$ where $\rho_{hyp}$
is the hyperbolic holonomy of $M$ and all weights vanish. The following theorem summarizes our new results for patterns over $M$ (all terms are defined in Section \ref{I-cusped}).

\begin{teo}\label{more_cusped} 
  Let $M$ be an arbitrary cusped manifold, $X(M)$ the variety of {\rm
    augmented} $PSL(2,\C)$-characters of $M$, and $X_0(M)\subset X(M)$
  the irreducible component of $\rho_{\rm hyp}$. There is a canonical
  non empty Zariski open subset $\Omega(M)$ of $X_0(M)$ containing
  $\rho_{\rm hyp}$ such that:
\smallskip
  
  (1) For every odd integer $N\geq 3$ and every pattern
  $(V,\rho, (h, k))$ such that $\rho\in \Omega(M)$, there is a quantum hyperbolic invariant $\Hh_N(V,\rho,(h,k))\in \mc/\mu_{2N}$ satisfying $\Hh_N(V,\rho_{hyp}, (0, 0)) =   \Hh_N(M)$.
\smallskip

  (2) {\rm (Analytic Families)} Fix a topological support $(V,(h_c,k_c))$ and $h_f\in H^1(V;\mz/2\mz)$. For every odd integer $N\geq 3$, the invariants $\Hh_N(V,\rho, (h, k))$ define a regular rational function $\Hh_N^{h_f,h_c,k_c}: \tilde{\Omega}(M)_N\ra \mc/\mu_{2N}$ on a
  determined $(\mz/N\mz)^2$-covering space $\tilde{\Omega}(M)_N$ of
  $\Omega(M)$.
  
 (3) {\rm (Resolution of the sign ambiguity)} If $N\equiv 1$ mod$(4)$, or
 $N\equiv 3$ mod$(4)$ and $h_c=0$, then the statements (1) and (2) above hold true by replacing $\mu_{2N}$ with $\mu_N$.
\end{teo}
 {\bf Comments:}

{\bf a)} $X(M)$ is a complex algebraic variety and $\rho_{\rm hyp}$ a
regular point of $X(M)$. Hence there is a unique irreducible component
$X_0(M)$ of $X(M)$ containing $\rho_{\rm hyp}$. As $M$ has only one
cusp, $X_0(M)$ is an algebraic curve and $\Omega(M)$ the complement of
a finite set of points.

{\bf b)} Theorem \ref{more_cusped} (2) is a rough qualitative
formulation of the $N$-th analytic family $\Aa_N(Y)$ associated to the
topological support $Y=(V,(h_c,k_c))$ of $M$.

%{\bf c)} We introduced $\alpha_N(\Tt)$ in \cite{GT} to solve a ``symmetrization problem'' implying a full invariance of state sums under versions of triangulation moves called {\it transits} (see Section \ref{I-enhanced}). The invariance of $\alpha_N(V,\rho, (h, k)):\equiv_2 \alpha_N(\Tt)$ follows from the closed formula \eqref{globalformula}. On the other hand, the invariance of $\Hh_{N,red}(V,\rho, (h, k)):\equiv_{2N} \Hh_{N,red}(\Tt)$ is obtained indirectly, from that of $\Hh_N(V,\rho,(h,k)):\equiv_{2N} \Hh_N(\Tt)$ (proved via transit arguments) and $\alpha_N(V,\rho, (h,k))$. We found no direct proof because the state sums $\Hh_{N,red}(\Tt)$ have a non local behaviour under weak branching changes whereas the transit arguments are local by nature (see Remark \ref{no-local-basic} for details).

{\bf c)} The arguments of Theorem \ref{more_cusped} (1)-(2) apply verbatim to prove the invariance of simplicial formulas defining the $PSL(2,\mc)$-Chern-Simons section, realized as an analytic equivariant function $\mathfrak{h}^*\Hh_1^{\Ff_{\rm EP}, h_f}: \tilde{\Omega}(M)_\infty \ra \mc$ on a $\mz^2$-covering $\tilde{\Omega}(M)_\infty$ of $\Omega(M)$ (see Section \ref{qCS}).

%{\bf e)} It is not clear to the authors whether or not it is possible to remove the sign ambiguity of the invariants $\alpha_N(V,\rho, (h,k))$ and $\Hh_{N,red}(V,\rho, (h, k))$.

\subsection{QHI of other patterns}\label{I-QHI-recall}
In order to explain the nature of the QHI in Theorem
\ref{more_cusped} (1), and why this result is not straightforward, it is useful
to recall a few general facts from \cite{Top,GT,AGT}.  

\smallskip

{\bf QHFT partition functions.} The QHI of QHFT patterns have been defined in \cite{AGT}.  
The topological supports of QHFT patterns have the form $(V,L,(h_c,k_c))$ where:
\begin{itemize}
\item $V$ is a compact oriented connected $3$-manifold with (possibly
  empty) boundary $\partial V$ made by torus components; if $\partial
  V = \emptyset$ we will use the notation $W$ instead of $V$.
\item $L$ is a {\it non--empty} link in the interior of $V$;
\item the bulk $\&$ boundary $c$-weights $(h_c,k_c)\in H^1(V;\Z/2\Z)\times H^1(\partial V; \Z)$ satisfy \eqref{cconstraint}.
\end{itemize}
The QHFT patterns are obtained by completing $(V,L,(h_c,k_c))$ with a couple
$(\rho,(h_f,k_f))$, where $\rho$ is any $PSL(2,\C)$-character of $V$ and
$(h_f,k_f)$ are bulk and boundary $f$-weights satisfying \eqref{fconstraint} with respect to $\rho$. 

When $V=W$ is a closed $3$-manifold, $(k_c,k_f)$ disappears so that $h=(h_c,h_f)\in H^1(W;\Z/2\Z)^2$. In \cite{Top,GT} we
defined the QHI $\Hh_N(W,L,\rho,h)$ in that situation, for every character $\rho$
and weight $h$. A specialization is $H_N(S^3,L):= \Hh_N(S^3,L,\rho_{\rm triv},0)$, where $\rho_{\rm triv}$ is the trivial character of $S^3$. In \cite{LINK} we obtained the following result, which establishes a connection with Jones invariants. 
\begin{teo}\label{QH=J} For every link $L$ in $S^3$ and every odd integer $N\geq 3$ we have
$$H_N(S^3,L) \equiv_{N}\ <L>_N \ = J_N(L)(e^{2i\pi/N})$$ 
where $<L>_N$ is the link invariant defined by the
enhanced Yang-Baxter operator extending the Kashaev R-matrix, and
$J_N(L)(q)\in \Z[q^{\pm 1}]$ is the colored Jones polynomial, normalized so that
$J_N(K_U)(q)=1$ on the unknot $K_U$.
\end{teo}    

\begin{remark}{\rm The second equality is due to \cite{MM}. In
    \cite{Top, GT, AGT} we quoted occasionally the first one as a
    motivating fact. Later we realized that we were unable to derive a
    complete proof from the existing literature (in particular
    \cite{K0,K1}), so we provided an independent one in \cite{LINK}, under the ambiguity $\equiv_{2N}$.  The above statement with $\equiv_{N}$ follows from the state sum sign correction introduced in the present paper (see also Remarks \ref{o-graph} and \ref{norm}).}
\end{remark}
When $\partial V\ne \emptyset$ the QHFT partition
functions are more sophisticated (see \cite{AGT}). In particular the
link $L$ contains an essential simple closed curve on
each boundary component which actually encodes a Dehn filling
instruction. Anyway, also in this case the invariants
$\Hh_N(V,L,\rho,(h,k))$ are defined for arbitrary characters and
weights.

\medskip

{\bf Relation with the QHI of cusped manifolds.}  The
patterns over cusped manifolds and the QHFT patterns are complementary in the sense that the link $L$ is {\it empty} in the former.  In \cite{GT,
  AGT} the proof of invariance of $\Hh_N(M)=\Hh_N(V,\rho_{hyp},(0,0))$ differs to many extents
from the one for the QHFT partition functions. It uses the ``volume
rigidity'' for cusped manifolds (see eg. \cite{F}), Thurston's
hyperbolic Dehn filling theorem, a construction of certain auxiliary
invariants $\Hh_N(V,\rho_{\rm hyp},(0,0),a)$ that depend a priori on
an additional datum ``$a$'', and finally a surgery formula. Set $H_N(W,L,\rho):= \Hh_N(W,L,\rho,0)$, where $W$ is closed (see our notations above). By
combining all these results we proved:
\begin{teo}\label{fillings}  {\rm (\cite{AGT}, Section 6.2)} Let $W_n$ be a sequence of closed hyperbolic Dehn fillings of $M$ whose holonomies $\rho_n$ (considered as $PSL(2,\mc)$-characters on $M$) converge to $\rho_{hyp}$ in $X(M)$. Denote by $L_n$ the geodesic core of the solid torus that fills $V$ to produce
  $W_n$. For every odd $N\geq 3$ and every additional datum ``$a$'' we have
  \begin{equation}\label{limitcusped}
  \lim_{n\rightarrow \infty} H_N(W_n,L_n,\rho_n)\equiv_{2N} \Hh_N(V, \rho_{\rm
    hyp},(0,0),a) . 
    \end{equation} 
    Hence ``$a$'' is eventually immaterial, and the limit
    defines $\Hh_N(M)\in \mc/\mu_{2N}$.
  \end{teo}
  The normalization $(h,k)=(0,0)$ on
  the right side of \eqref{limitcusped} is a by-product of the
  proof. Under some additional assumptions on $M$ (for
  instance if $M$ is ``very gentle'' according to \cite{GT, AGT}; then
  ``$a$"$=\emptyset$), we could avoid the delicate surgery argument
  and define the invariants $\Hh_N(V,\rho_{\rm hyp},h,k)$
  for arbitrary weights relative to the hyperbolic holonomy.

\subsection{State sums over weakly branched triangulations}\label{I-En-Ss} One complication with the construction of the QHI of cusped manifolds in \cite{GT,AGT} depended on a technical difficulty that we overcome in the present paper, and that we are going to illustrate. 

For every topological support $Y_\Pp$,
denote by $\hat V$ the compact space obtained by filling each boundary
component of $\partial V$ with the cone over it. If $\partial V =
\emptyset$ then $\hat V = W$; $\hat V$ has a finite set of
non-manifold points, the vertices of the filling cones. For every
pattern $\Pp$ supported by $Y_\Pp$, $\Hh_N(\Pp)$ is computed by state
sums over certain ``decorated'' triangulations $T$ of $\hat V$,
depending on the choice of a point $z_\rho$ in the associated gluing
variety $G(T)$ such that $z_\rho$ represents the character $\rho$ (see
Section \ref{I-triang}), and on a suitable encoding of the weights
$(h,k)$. Moreover, $T$ is equipped with a {\it
  branching}; equivalently, $T$ carries a structure of $\Delta$-{\it
  complex} in the sense of \cite{HATCHER}.

In the case of QHFT partition functions this is not so demanding: $T$
can be a ``quasi-regular'' triangulation (every edge has distinct
endpoints), a branching $b$ can be induced for example by a total
ordering of the vertices, and $z_\rho$ can be realized by means of a
so-called ``idealization'' of $PSL(2,\C)$-valued 1-cocycles on
$(T,b)$.

On another hand, for a cusped manifold $M$ we use {\it ideal} triangulations $T$ of
$\hat V$ such that the gluing variety $G(T)$ contains a point $z_h$
representing the hyperbolic holonomy $\rho_{\rm hyp}$, and having
coordinates with non-negative imaginary part (we say that $z_h$ is {\it non-negative}). Such triangulations exist for
every $M$, for instance any maximal subdivision of the canonical
Epstein-Penner cell decomposition has this property. However we do not know if
every cusped manifold has such a triangulation $T$ admitting a branching $b$. For
instance, the canonical Epstein-Penner decomposition of the ``figure-8-knot sister'' (m003 in Snappea's census) is made of two regular hyperbolic ideal tetrahedra and does not carry any
branching. The ``very gentle'' manifolds $M$ mentioned after Theorem \ref{fillings} admit {\it by
definition} a branched triangulation $(T,b)$ with a non-negative point $z_h$ in $G(T)$.

In order to get Theorem \ref{more_cusped} (1) we relax branchings to
{\it weak branchings} which exist on every triangulation, and this
leads us to include {\it $2$-face tensors} in the state sum
formulas. In this setup, as well as to cover arbitrary characters of
$M$ in $\Omega(M)$, the proof of the state sum invariance requires
additional arguments with respect to \cite{GT,AGT}.
 
\subsection{Plan of the paper}\label{plan} Let $V$ be as in Section \ref{I-QHI-recall} and $\hat V$ as in
Section \ref{I-En-Ss}.
\smallskip

In Section \ref{I-triang} we recall a few general facts about
triangulations endowed with {\it pre-branchings}, {\it weak branchings}, or {\it branchings}, and the
associated {\it gluing varieties} of $M$.

In Section \ref{I-CAC} we construct the {\it analytic configuration} $\Aa_N(T,\tilde b, c)$ for every odd integer $N\geq 3$,  weakly branched triangulation $(T,\tilde b)$ of $\hat V$, and {\it rough charge} $c$ on $T$ (suitably specialized
``global charges'' will eventually encode the $c$-weights below). In particular $\Aa_N(T,\tilde b, c)$ contains an infinite Abelian covering of the gluing
variety, $p_\infty: G(T,\tilde b)_\infty\to G(T)$, and an analytic function $\Hh_{N}(T,\tilde b,c): G(T,\tilde b)_\infty \to \C$. In the case of QHFT patterns it is described qualitatively by:
\begin{prop}\label{informalPF} 
  For every topological support $(V,L,(h_c,k_c))$ there is a
  weakly branched triangulation $(T,\tilde b)$ of $\hat V$ and a
  global charge $c$ on $T$ such that:
\begin{enumerate}
\item $c$ encodes the $c$-weight $(h_c,k_c)$.
\item For every QHFT pattern $(V,L,\rho,(h,k))$ with topological
  support $(V,L,(h_c,k_c))$, there is a point $u\in G(T,\tilde b)_\infty$ such that
  $p_\infty(u)$ represents the character $\rho$ and for every odd
  $N\geq 3$ the scalar $\Hh_N(V,L,\rho,(h,k))\equiv_{2N} \Hh_{N}(T,\tilde b,c)(u)$ does not depend on the choice of $(T,\tilde b,c)$ and $u$.
\end{enumerate}
\end{prop}
In the case of patterns over a cusped manifold $M$, we have:
\begin{prop}\label{informal-cusp} For every topological support 
  $Y=(V,(h_c,k_c))$ over $M$, there is a
  determined Zariski open subset $\Omega(M)$ of $X_0(M)$ containing
  the hyperbolic holonomy $\rho_{hyp}$, and there is a weakly branched {\rm ideal}
  triangulation $(T,\tilde b)$ of $\hat V$, and a global charge $c$ on
  $T$ such that:
\begin{enumerate}
\item $c$ encodes the $c$-weight $(h_c,k_c)$.
\item The gluing variety $G(T)$ contains a non-negative point $z_h$ representing $\rho_{hyp}$.
\item There is an irreducible component $Z$ of $G(T)$, a Zariski open
  subset $\Omega_Z$ of $Z$ containing $z_h$, and a homeomorphism $\rho:\Omega_Z\to \Omega(M)$ extending a regular rational isomorphism between Zariski open subsets, such that $\rho(z)$ is the holonomy of $V$ represented by $z$.
\item For every pattern $(V,\rho,(h,k))$ over $M$ with
  topological support $Y$ and for every $\rho\in \Omega(M)$, there is
  a point $u\in Z_\infty := p_\infty^{-1}(Z)$ such that $\rho =
  \rho(p_\infty(u))$, and for every odd $N\geq 3$ the scalar $\Hh_N(V,\rho,(h,k))\equiv_{2N} \Hh_{N}(T,\tilde b,c)(u)$ does not depend on the choice of $(T,\tilde b,c)$ and $u$.
\end{enumerate}
\end{prop} 

We will use concrete models of the finite coverings and regular
rational maps in Theorem \ref{more_cusped} (2) by considering a
suitable factorization of $\Hh_N(T,\tilde b,c):Z_\infty \to \C$.
\smallskip

In Section \ref{I-cusped} we develop the content of Proposition
\ref{informal-cusp} and prove Theorem \ref{more_cusped} (1)-(2) and the analogous result when $N=1$ (see the comment (c) and Corollary \ref{ratinv}). 

In Section \ref{I-PF} we indicate briefly how to deal with QHFT
partition functions. 

In Section \ref{N-calculus} we collect a few facts about a diagrammatic
calculus for weakly branched triangulations. This calculus is used in Section \ref{I-enhanced}, which contains the invariance proof of the state sums defined over such triangulations.

In Section \ref{sign} we prove Theorem \ref{more_cusped} (3) and the analogous result for QHFT partition functions. As a by product we show that global compensations of the local sign ambiguities imply that the QHI $\Hh_N(\Pp)$ defined by means of branched triangulations have no sign ambiguity. This applies to QHFT partition functions and to very gentle cusped manifolds.

In Section \ref{EXAMPLES} we describe the quantum hyperbolic invariants of three cusped manifolds: the figure-eight knot complement, its ``sister", and the complement of the knot $5_2$. We present the results of numerical computations showing that the functions $\Hh_N^{h_f,h_c,k_c}$ depend on the weights as $N\rightarrow + \infty$, and recover the volume for some specific choices of the weights.

\smallskip

{\bf Notations.}  In general we denote by $N$ the odd integers $\geq 3$, including the case $N=1$ when specified, and we put $m:=(N-1)/2$ and $\zeta := \exp(2\pi i/N)$. The set $\Ii_N = \{0,\dots,
N-1\}$ is identified to the group $\mz/N\mz$ and $[n]_N\in \Ii_N$
denotes the residue modulo $N$ of $n\in \mn$. We let $\delta_N(n) :=
1$ (resp. $0$) if $[n]_N =0$ (resp. $[n]_N \neq 0$).  

\medskip

{\bf Acknowledgments.}  The first author's work was supported by the
ANR projects {\it G\'eom\'etrie et Topologie Quantiques} (ANR-08-JCJC-0114-01) and "Extensions des th\'eories de Teichm\"{u}ller-Thurston" (ANR-09-BLAN-0116-01).  Many thanks are due
to Columbia University, Fukuoka University, the MSC at Tsinghua
University and the University of Utah at Moab where parts of this
work was presented. The second author's work is
supported by the italian FIRB project {\it Geometry and topology of
  low-dimensional manifolds}.  The numerical computations presented in Section \ref{EXAMPLES} were obtained by using the Maple software.

\medskip

\section{Structured triangulations and gluing varieties}
\label{I-triang}

{\bf Triangulations.} Let $V$ be as in Section \ref{I-QHI-recall},
that is, a compact oriented connected $3$-manifold with (possibly
empty) boundary $\partial V$ made by toric components. Denote by $\hat
V$ the space obtained by taking the cone over each boundary component
of $\partial V$.  A triangulation $T$ of $\hat V$ is a
collection of oriented tetrahedra $\Delta_1,\dots,
\Delta_s$ together with a complete system $\sim$ of pairings of their
$2$-faces via orientation reversing affine isomorphisms, such that the
oriented quotient space $$T := \coprod_{i=1}^s \Delta_i/\sim $$ is
homeomorphic to $\hat V$, preserving the orientations.  We will
distinguish between the $2$-faces, edges and vertices of
the disjoint union $\textstyle \coprod_{i=1}^s \Delta_i$, and the ones of $T$ after the
$2$-face pairings. In particular we denote by $E(\{\Delta_i\})$ and
$E(T)$ the set of edges of $\textstyle \coprod_{i=1}^s \Delta_i$ and
$T$ respectively, and we write $E\to e$ to mean that an edge $E\in E(\{\Delta_i\})$ is identified to $e \in E(T)$ under the $2$-face pairings. The $2$-faces of each tetrahedron $\Delta_i$ have the boundary orientation defined by the rule: {\it ``first the outgoing
  normal''}. When $\partial V \neq \emptyset$, the non-manifold points of $\hat V$ are necessarily vertices of $T$.  A triangulation $T$
of $\hat V$ is called {\it ideal} if the set of vertices of $T$
coincides with the set of non-manifold points of $\hat V$.

\medskip

{\bf Gluing varieties.}  Let $T$ be as above.  For every tetrahedron $\Delta_j$ choose a vertex $v^j$. Order the edges of the opposite $2$-face $F^j$ so that the induced cyclic ordering is the {\it opposite} of the boundary orientation. Denote these edges by $E^j_0,E^j_1, E^j_2$. Give a label $u^j_r\in \Cc$ to $E^j_r$ and the opposite edge, where $\Cc := \C \setminus \{0,1\}$ and $r\in \{0,1,2\}$. The gluing variety of $T$ is the algebraic subset of $\Cc^{3s}$ with coordinates $(u^1_0,u^1_1,u^1_2,\ldots,u^s_0,u^s_1,u^s_2)$ and defining equations ($r$ is mod$(3)$):
$$\forall\ j\in \{1,\ldots,s\}, r\in\{0,1,2\}, u^j_{r+1}(1-u^j_r)=1$$
$$\forall \ e\in E(T), \prod_{E\to e} u(E)=1.$$
By the first set of $s$ ``tetrahedral relations'', we see that the gluing variety is the graph of an explicit regular rational map defined on an algebraic subset of $\Cc^{s}$ defined by the second set of ``edge relations''. The auxiliary choices of ordered edges $E^j_r$ being immaterial, these algebraic varieties are canonically isomorphic. We denote them by $G(T)$. 

Every point $u\in G(T)$ represents a
$PSL(2,\C)$-character $\rho(u)$ of $V$. Its components $u^j_0$, $u^j_1$, $u^j_2$ can be interpreted as
the cross-ratio parameters of an isometry class of oriented
hyperbolic ideal tetrahedra associated to $\Delta_j$, that we denote by $(\Delta_j,u^j)$. Their imaginary parts have a same sign $\epsilon^j\in \{-1,0,1\}$ (by convention $\epsilon^j= 0$ if the imaginary parts are zero). By a classical result of Schl\"affli, the algebraic volume of $(\Delta_j,u^j)$ is given by 
\begin{equation}\label{algvol}
{\rm Vol}_{\rm alg}(\Delta_j,u^j) := \epsilon (u^j) {\rm Vol}(\Delta_j,u^j) = {\rm D}_2(u^j_0)
\end{equation}
where Vol is the geometric positive volume and D$_2$ the Bloch-Wigner dilogarithm. When the components $u^j_r$ are real, $(\Delta_j,u^j)$ is
degenerate and both sides of \eqref{algvol} vanish. By summing the
algebraic volumes of the $(\Delta_j,u^j)$s we get a {\it
  volume function} 
  \begin{equation}\label{volf}{\rm Vol}: G(T)\to  \mr . 
\end{equation}
If $G(T)$ is non empty, for every point $u\in G(T)$, ${\rm Vol}(u)$ coincides with the (intrinsically defined) volume ${\rm Vol}(\rho(u))$ of the character $\rho(u)$.  In general $G(T)$ might be empty, but $\hat V$ has always triangulations $T$ such that $G(T)$ is non trivial. A first general result concerns its dimension.
\begin{teo}\label{cusped_dimension} {\rm (\cite{NZ}, \cite{N0}; see
    also \cite{BP})} Assume that $V$ has one torus boundary
  component. Let $T$ be an ideal triangulation of $\hat V$. If the
  gluing variety $G(T)$ is non empty, then it is a complex algebraic
  set of dimension $\geq 1$.
\end{teo}
This result depends on the combinatorial properties of $T$. If the interior of $V$ is a cusped
manifold $M$, then it has the canonical Epstein-Penner (EP) cell
decomposition by embedded convex hyperbolic ideal polyhedra (see for
instance \cite{BP}). Then we have:

\begin{prop} \label{EP-triang} The maximal subdivisions of the EP cell
  decomposition of $M$ define a {\rm finite} set $\Tt_{EP}(M)$ of
  ideal triangulations of $\hat V$, such that for every $T\in
  \Tt_{EP}(M)$ the gluing variety $G(T)$ contains a non-negative
  point $u_h\in G(T)$ such that $\rho(u_h)=\rho_{\rm hyp}$, and ${\rm
    Vol}(u_h) = {\rm Vol}(M)$.
\end{prop}

Every non-degenerate hyperbolic ideal tetrahedron of $(T,u_h)$ has strictly positive (geometric) volume, but in general one cannot avoid some degenerate tetrahedra. 

For other manifolds the nature of $G(T)$ is not so well-known. For instance consider the case of a closed $3$-manifold $W$. A triangulation $T$ of
$W$ is called {\it quasi-regular} if every edge of $T$ has distinct
vertices in $T$. It is clear that every triangulation has a quasi-regular subdivision. Take one, and fix a total ordering of the vertices. For every edge $e$ with endpoints $v$ and $v'$, orient $e$ from $v$ to $v'$ if
$v<v'$. Every simplicial $PSL(2,\C)$-valued $1$-cocycle $z$ on $T$,
defined by using this edge orientation, represents a character
$\rho(z)$ of $W$. Then, to any sufficiently generic cocycle $z$ one can associate a point $u\in G(T)$ such that $\rho(z)=\rho(u)$, and we have (\cite{Top}):
\begin{prop}\label{totalGT} Let $T$ be a quasi-regular triangulation   
of $W$. For every character $\rho$ of $W$ there is a point
$u\in G(T)$ such that $\rho=\rho(u)$.
\end{prop}

Clearly, the point $u$ is far to be unique. For instance, if $W=S^3$ every point of $G(T)$ represents the trivial character. A similar, slightly more elaborated result holds for all other topological supports of QHFT patterns; it uses triangulations of $\hat V$ obtained from quasi-regular relative triangulations $(T,\partial T)$ of $(V,\partial
V)$ by adding a cone over each component of $\partial T$.  

The method used to prove Proposition \ref{totalGT} is reminiscent of Thurston's spinning construction, and is strictly related to it when $W$ is hyperbolic.  In that case, the following result, which is proved by using the spinning construction, agrees with Proposition \ref{totalGT} in the case of quasi regular triangulations.

\begin{prop}\label{spinning}(\cite{LUO2})
  Let $T$ be a triangulation of a closed oriented hyperbolic
  $3$-manifold $W$ such that no edge is a null-homotopic loop in $W$. Then there exists $u\in G(T)$ such that ${\rm Vol}(u)={\rm Vol}(W)$, and moreover its holonomy $\rho(u)$ is the hyperbolic holonomy. 
\end{prop}

{\bf Variations on branched triangulations.}  Define a {\it
  pre-branched tetrahedron} $(\Delta,\sigma)$ as an oriented
tetrahedron $\Delta$ with a choice $\sigma$ of co-orientations of the $2$-faces, such that two co-orientations are ingoing and two are outgoing. As every $2$-face has the boundary orientation, by
duality $\sigma$ can be interpreted as a system of $2$-face
orientations. 
\begin{figure}[ht]
\begin{center}
 \includegraphics[width=3.5cm]{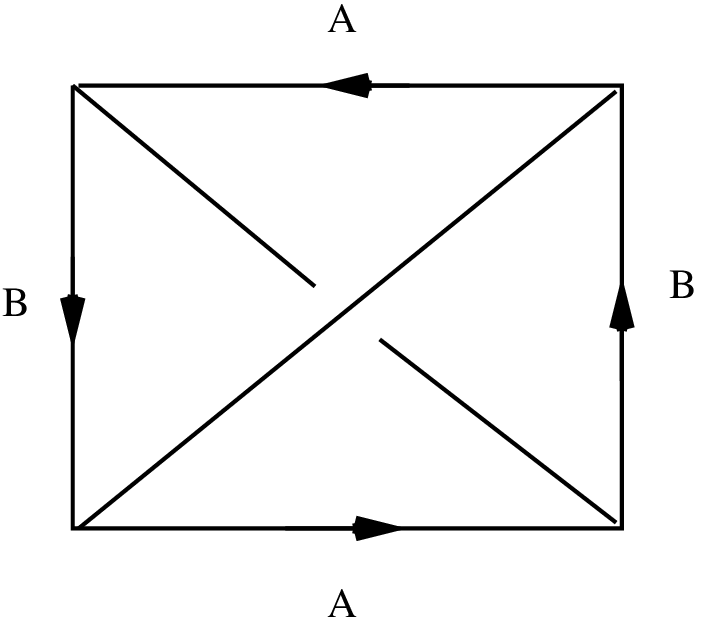}
\caption{\label{Taut_Delta} Pre-branched tetrahedron.}
\end{center}
\end{figure}

Figure \ref{Taut_Delta} shows a pre-branched tetrahedron
$(\Delta,\sigma)$ embedded in $\mr^3$, with coordinates $(x_1,x_2,x_3)$
such that the plane of the picture is $\{x_3=0\}$. We put on $\Delta$ the orientation induced from $\mr^3$, and assume that the two $2$-faces above (resp. below) the plane are those with outgoing (resp. ingoing) co-orientations. This specifies two {\it diagonal edges} and four {\it square edges}. Every square edge is oriented as the common boundary edge of
two $2$-faces with opposite co-orientations. So the square edges form an oriented quadrilateral. Using the orientation of $\Delta$, one can also distinguish among the square edges two pairs of opposite edges, called {\it $A$-edges} and {\it $B$-edges}. The orientation of the diagonal edges is not determined.

\begin{figure}[ht]
\begin{center}
 \includegraphics[width=8cm]{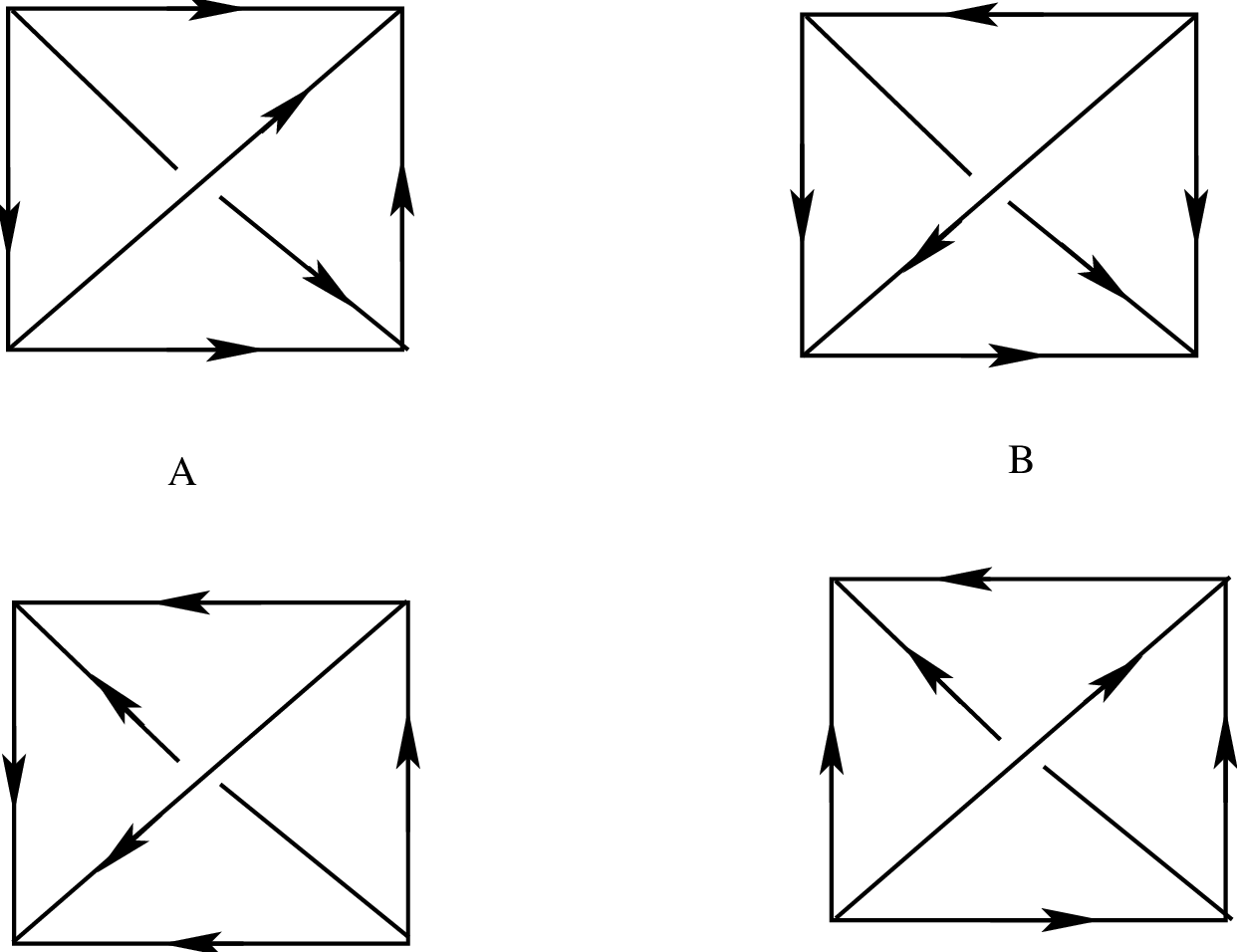}
\caption{\label{Branched_Delta} Branched tetrahedra inducing the same
pre-branched tetrahedron.}
\end{center}
\end{figure}

An oriented tetrahedron $\Delta$ becomes a $3$-simplex by ordering its vertices. This is equivalent to a system $b$ of orientations of the edges, called
(local) {\it branching}, such that the vertex $v_j$ has $j$ incoming edges ($j\in \{0,\ldots,3\}$). The
$2$-faces of $(\Delta, b)$ are ordered as the opposite
vertices, and $b$ induces a branching $b_F$ on each $2$-face $F$. The branchings $b$
and $b_F$ define orientations on $\Delta$ and $F$ respectively, the
{\it $b$-} and {\it $b_F$-orientations}. The $b$-orientation may coincide or not with the orientation of
$\Delta$. We encode this by a sign, $*_b\in \{-1,+1\}$. The boundary orientation and the $b_F$-orientation agree on two $2$-faces. Hence $b$ defines a pre-branching $\sigma_b$. On another hand, given a pre-branching $\sigma$ on $\Delta$ there are exactly four branchings such that $\sigma_b =\sigma$. They can be obtained by choosing an $A$ (resp. $B$) edge, reversing its orientation, and completing the resulting square edge orientations to define a branching $b$ (this can be done in a single way; see Figure \ref{Branched_Delta}). Note that $*_b=1$ (resp. $*_b=-1$) if and only if we have chosen an $A$ (resp. $B$) square edge, and the square edge is $[v_0,v_3]$. The diagonal edges are $[v_0,v_2]$ and $[v_1,v_3]$.

One can extend these notions to triangulations $T$ of $\hat V$:
\begin{itemize}
\item A {\it pre-branched triangulation} $(T,\sigma)$ is formed by pre-branched tetrahedra $(\Delta_j,\sigma_j)$ such that the $2$-face co-orientations match under the $2$-face pairings.
\item A {\it weakly-branched triangulation} $(T,\tilde b)$ is formed by branched tetrahedra $(\Delta_j, b_j)$
such that the induced pre-branched tetrahedra $(\Delta_j,\sigma_{b_j})$ form a pre-branched
triangulation $(T,\sigma)$. 
\item A {\it branched triangulation} $(T,b)$ is formed by branched tetrahedra $(\Delta_j, b_j)$ such that the branchings (ie. the edge orientations) match under the 2-face pairings. 
\end{itemize}

\begin{remark}\label{DELTA-COMP} {\rm Branched
    triangulations of $\hat V$ and $\Delta$-{\it complexes} over $\hat V$ \cite{HATCHER} are equivalent notions. In particular the simplicial $3$-chain $\textstyle \sum_j *_{b_j}(\Delta_j,b_j)$ represents the
    fundamental class in $H_3(\hat V;\Z)$.}
\end{remark}
\smallskip

Let $T$ be a triangulation of $\hat V$. Denote by $\bar V$ the compact
$3$-manifold with boundary obtained by removing a small open
$3$-ball around every vertex of $T$ which is a manifold point.
Clearly, $V=\bar V$ if and only if $T$ is an ideal triangulation. A pre-branched
triangulation $(T,\sigma)$ of $\hat V$ can be described in a very
concrete way in terms of the standard spine $P$ of $\bar V$ dual to
$T$:

\begin{lem} \label{dualwbspine} There is a $1$-to-$1$ duality correspondence between
  the sets of pre-branchings of $T$ and those of $P$, where a prebranching of $P$ is defined as an orientation of the singular locus $\ {\rm Sing}(P)$ such that every vertex has two outgoing and two ingoing edges.
\end{lem}

The proof is evident, as every edge of ${\rm Sing}(P)$ is dual to a $2$-face of $T$. The notion of (weakly) branched
triangulation $(T,\tilde b)$ has a natural dual counterpart as
{\it (weakly) branched spine} $(P,\bar b)$.

\smallskip

\begin{figure}[ht]
\begin{center}
 \includegraphics[width=9cm]{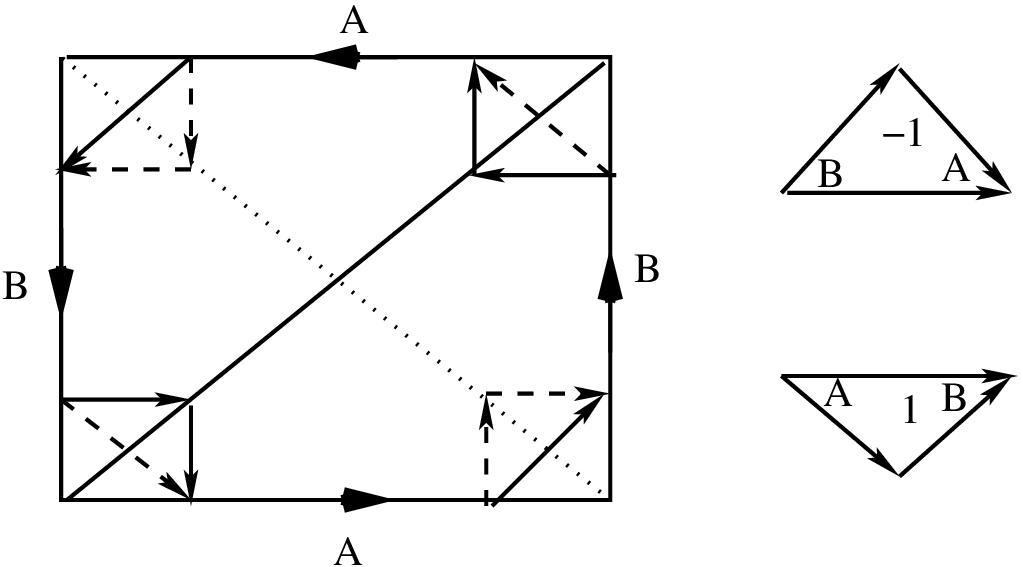}
\caption{\label{pb-boundary} Branched boundary of pre-branched triangulations.}
\end{center}
\end{figure}

{\bf The branched boundary of a pre-branched triangulation.}  The triangulation $T$ of $\hat V$ yields a decomposition of  $\bar V$ into truncated tetrahedra. Their triangular $2$-faces form a triangulation $\partial T$ of $\partial \bar V$. If $(T,\sigma)$ is a pre-branched triangulation, then
$\partial T$ has a branching $b_{\partial \sigma}$ defined on the $2$-faces of $\partial T$ as in Figure \ref{pb-boundary}.
 
\noindent The corners of every $2$-simplex of $(\partial T, b_{\partial \sigma})$
have a label in $\{A,B,\emptyset\}$, and a sign $\pm 1$
obtained by comparing the boundary and the $b_{\partial \sigma}$-orientations. Such a $\{A,B,\emptyset\}$-labelling may be defined more generally for any branched triangulation of a closed oriented surface. 

\begin{lem}\label{even} Each vertex of a branched triangulation of a closed
  oriented surface has an even number of adjacent corners with the empty label $\emptyset$. Hence every edge $e$ of a pre-branched triangulation $(T,\sigma)$ of $\hat V$ has an even number of diagonal edges $E\in E(\{\Delta_i\})$ such that $E\to e$.
\end{lem} 
\Dim For each $2$-face the edges adjacent to the corner with label $\emptyset$
have opposite $b$-orientations, while they agree elsewhere. The existence of an orientation on any transverse loop implies the first claim. The second follows by considering the triangulation $\partial T$. \cvd

\begin{remark}\label{taut/2}{\rm Give the labels $1\in \Z/2\Z$ to the diagonal
edges and $0\in \Z/2\Z$ to the square edges of every pre-branched tetrahedron $(\Delta,\sigma)$ of $(T,\sigma)$. Then Lemma \ref{even} and Lemma \ref{pre-b-exist} below imply that $\Z_2$-valued {\it taut} structures on $T$ always exist (a notion borrowed from F. Luo's work).}
\end{remark}

{\bf Networks and $\Nn$-graphs.} Consider an oriented graph $({\rm Sing}(P),\bar\sigma)$ as in Lemma \ref{dualwbspine}. Put around every vertex $v$ the dual branched tetrahedron $(\Delta_v,b_v)$ so that its $2$-faces intersect transversely the edges of ${\rm Sing}(P)$. Thus each edge connects a $2$-face
 $F^i$ of the ``initial'' branched tetrahedron
 $(\Delta_{v_i},b_{v_i})$ with a $2$-face $F^f$ of the ``final'' one
 $(\Delta_{v_f},b_{v_f})$, identified in $T$. The gluing map $\phi_e :F^i\rightarrow F^f$ is determined by a color $r(e) \in \Z/3\Z$ defined as follows. Set $J_3=\{0,1,2\}$. Denote by $S(J_3)$ the symmetric group on $J_3$, by $A(J_3)$ the subgroup of even permutations, and by $u^i_j$ and $u^f_j$ the vertices of $F^i$ and $F^f$ ($j\in J_3$). The map $\phi_e$ is determined by the permutation $\tau_e \in A(J_3)$ such that $\phi_e(u^i_j)=u^f_{\tau_e(j)}$. Then we put $$r(e) := \alpha^{-1}(\tau_e)\in \Z/3\Z$$
 where $\alpha: \Z/3\Z \to A(J_3)$ is the isomorphism given by $\alpha(j)=(012)^j$. We define $\Nn(T,\tilde b)$ as the oriented graph $({\rm Sing}(P),\bar\sigma)$ endowed with the $\Z/3\Z$-edge colors $r(e)$, and with the correspondence, for every vertex $v$, between the $2$-faces of $(\Delta_v,b_v)$ and the germs of edges adjacent to $v$. Clearly we have:
\begin{lem}\label{r-b} 
  A weakly branched triangulation $(T,\tilde b)$ is branched
  if and only if all the $\Z/3\Z$-edge colors $r(e)$ are equal to $0$.
\end{lem}

\begin{lem}\label{pre-b-exist} Any triangulation $T$ of $\hat V$
admits pre-branchings and hence compatible weak-branchings.
\end{lem}  

\begin{figure}[ht]
\begin{center}
 \includegraphics[width=8cm]{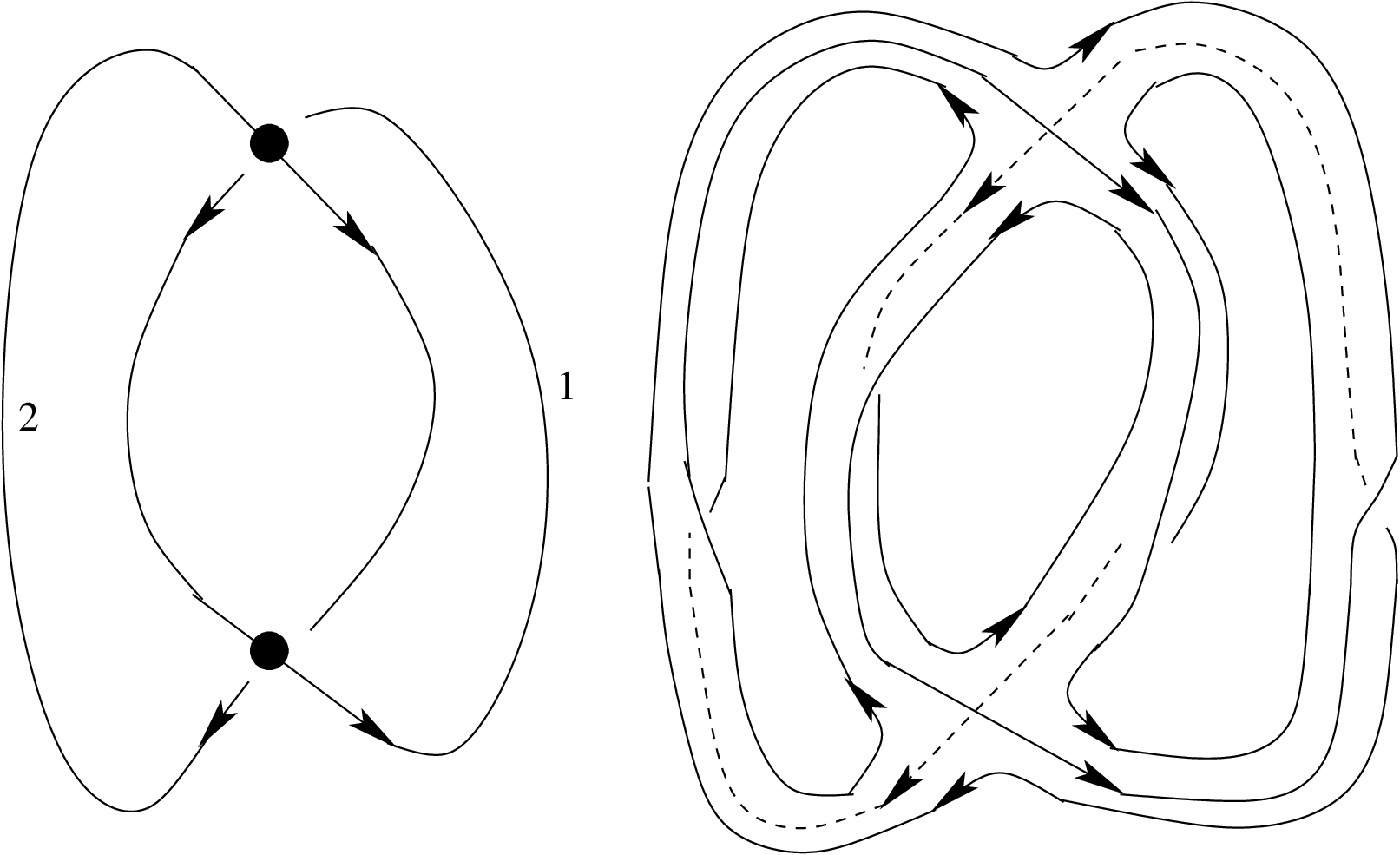}
\caption{\label{I-F8S} A $\Nn$-graph of $(T_1,\tilde b)$ and its decoding.} 
\end{center}
\end{figure}

We stress that this is no longer true for genuine branchings (see an example in Figure \ref{I-F8S}). If $T$ is quasi-regular, then every total ordering of the vertices induces a branching of $T$.

\medskip

\begin{figure}[ht]
\begin{center}
 \includegraphics[width=7cm]{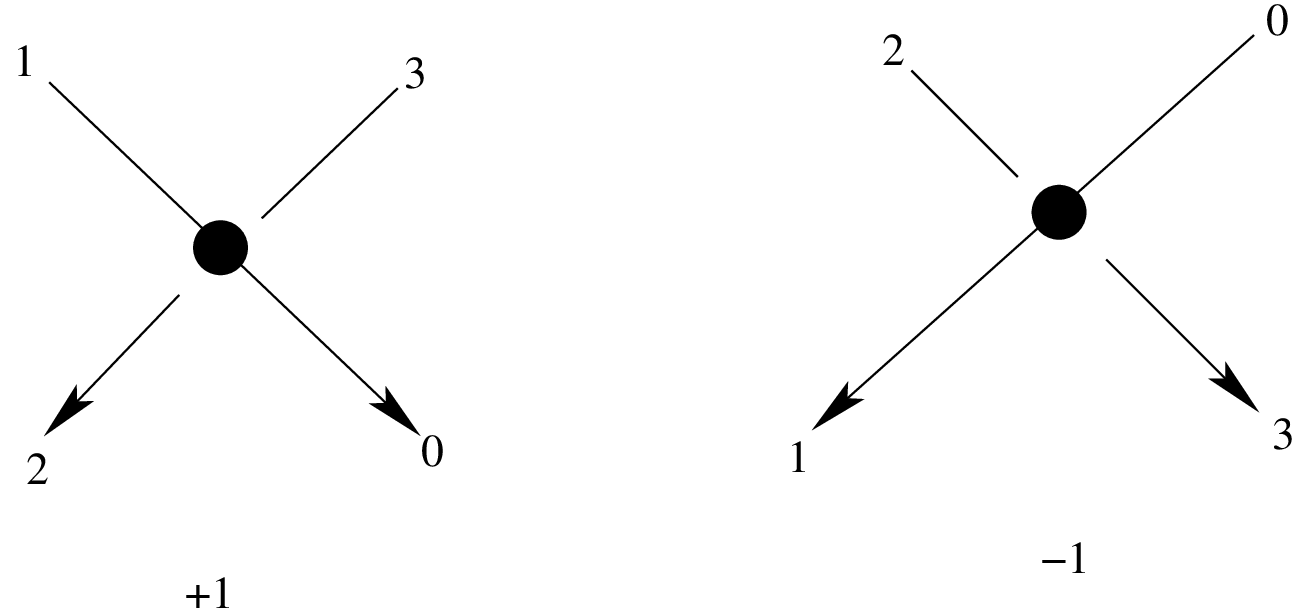}
\caption{\label{N-crossing} $\Nn$-graph crossings for $*_b=\pm 1$.}
\end{center}
\end{figure}

\begin{figure}[ht]
\begin{center}
 \includegraphics[width=9cm]{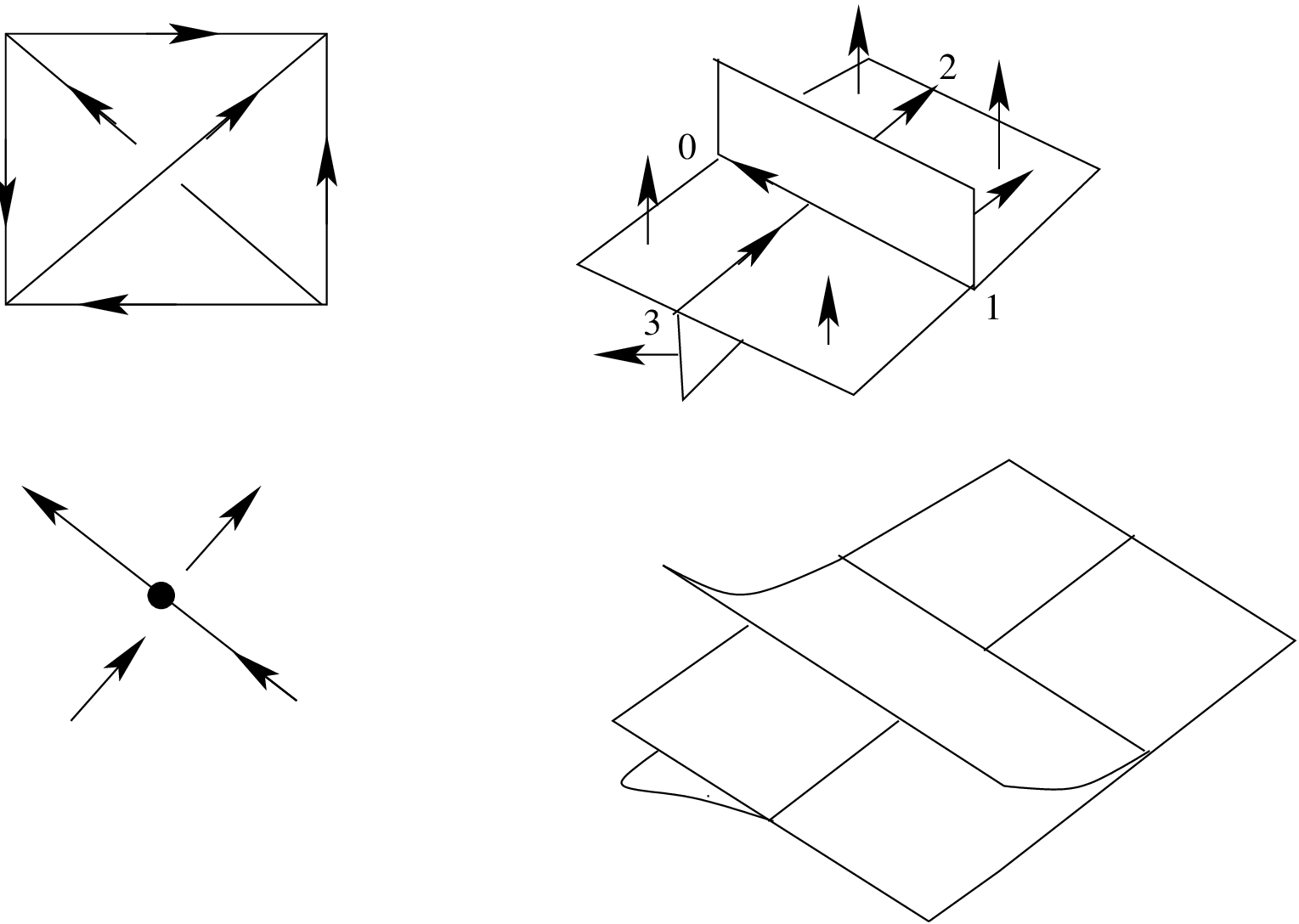}
\caption{\label{b+} Decoding of an $\Nn$-graph crossing ($*_b=1$).}
\end{center}
\end{figure}

It is useful to represent $\Nn(T,\tilde b)$ by a $\Nn$-{\it graph} $\Gamma$, defined as follows. Basically $\Gamma$ is a planar immersion of $({\rm Sing}(P),\bar\sigma)$ with normal crossings. It has two kinds of vertices: {\it essential crossings}, represented by a solid dot, corresponding to the vertices of $P$, and connected by the oriented edges of $\Gamma$, and non essential (``immaterial'') crossings due to the immersion. An essential crossing encodes a branched tetrahedron $(\Delta_v,b_v)$
as shown in Figure \ref{N-crossing}; the labels on the arc endpoints correspond to the $b_v$-ordering of the dual $2$-faces of $(\Delta_v,b_v)$. A full decoding is provided by Figure \ref{b+} in the case $*_b=1$, showing a branched tetrahedron, the dual ``butterfly'' with co-oriented (hence oriented) wings, and its conversion into a portion of a oriented branched surface. An arc of $\Gamma$ corresponding to an edge $e$ of
$({\rm Sing}(P),\bar\sigma)$ inherits the $\Z/3\Z$-color $r(e)$. If
$r(e)=0$ we omit it.

\begin{figure}[ht]
\begin{center}
 \includegraphics[width=8cm]{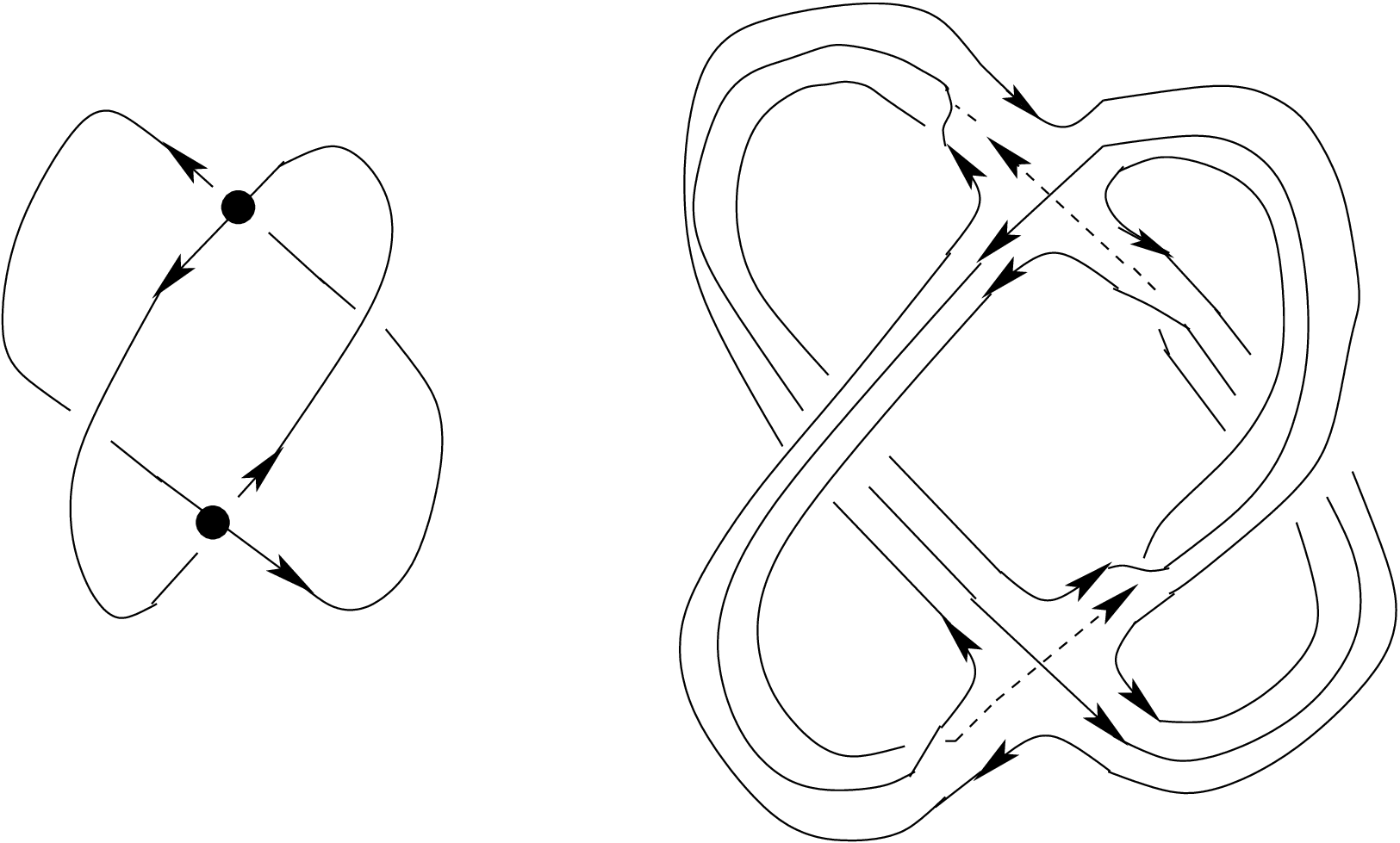}
\caption{\label{I-F8} A $\Nn$-graph of $(T_0,b)$ and its decoding.}
\end{center}
\end{figure}

Given $\Gamma$, an easy decoding procedure (extending
the one of Figure \ref{b+}) produces an embedding in $\mr^3$ of
a closed regular neighbourhood $N(P)$ of ${\rm Sing}(P)$ in the
standard spine $P$. We believe that it is enough to show this
in the case of the simplest cusped manifolds: $M_0$, the
figure-eight-knot complement in $S^3$, and its ``sister'' $M_1$, which
is the complement of a knot in the lens space $L_{5,1}$ (see
\cite{MF}). In both cases the Epstein-Penner decomposition is an ideal
triangulation made by two regular hyperbolic ideal tetrahedra. Denote by $T_0$ and $T_1$ the corresponding triangulations of $\hat V_0$ and $\hat V_1$. Figure \ref{I-F8} (resp. Figure \ref{I-F8S}) shows an $\Nn$-graph and its decoding for a branching $b$ of $T_0$ (resp. weak branching $\tilde b$ of $T_1$).  It is not hard to verify that {\it $T_1$ does not carry any branching}.

\begin{figure}[ht]
\begin{center}
 \includegraphics[width=9cm]{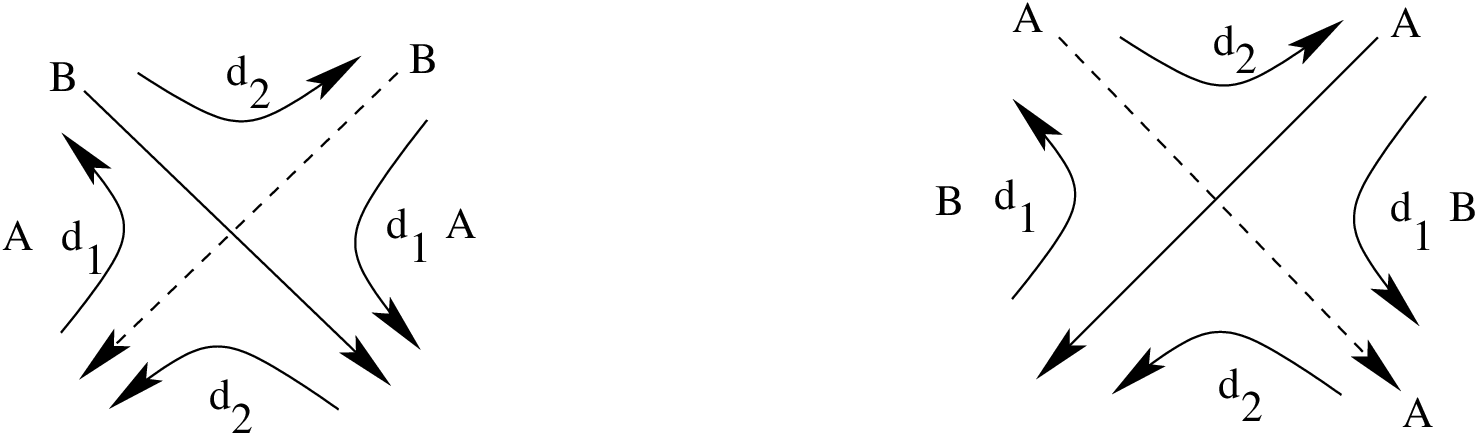}
\caption{\label{I-dec} Edge decorations on $\Nn$-graphs. }
\end{center}
\end{figure}

\medskip

{\bf Edge decorations.} We will use several labellings of the edges $E\in E(\{\Delta_i\})$ of a weakly branched triangulation, called {\it decorations}, such that for every tetrahedron opposite edges have the same label. An example is the decoration of the square edges with $A$ or $B$ according to the pre-branching, and of the diagonal edges with $\emptyset$. Every decoration $d$ is determined on each tetrahedron $\Delta_j$ by the triple $(d_0,d_1,d_2):=(d(E_0^j),d(E_1^j),d(E_2^j))$. In terms of decoded $\Nn$-graphs, they are placed as in Figure \ref{I-dec}, where we understand that
the over/under crossing arcs are labelled by $d_0$, and we show
also $A$, $B$, $\emptyset$.

\begin{remark}\label{o-graph}{\rm In the case of genuine branchings,
    $\Nn$-graphs were used in \cite{LINK} under the name of {\it
      normal o-graphs}, with the opposite
    convention for the sign $*_b$ (in accordance with the usual crossing signs
    of link diagrams). This choice gives the equality
    $H_N(S^3,L)\equiv_{N} J_N(\bar L)$, where $\bar L$ is the mirror
    image of the link $L$. The present convention yields the statement of Theorem \ref{QH=J}.}
\end{remark}

{\bf The model $G(T,\tilde b)$ of the gluing variety.} If $(T,\tilde b)$ is a
weakly branched triangulation, we can use the weak branching $\tilde b$ to fix
the auxiliary choices used in the definition of the
gluing variety $G(T)$. On every branched tetrahedron $(\Delta_j,b_j)$ let $v^j=v^j_3$, and order the edges of
$F^j_3$ as $E^j_0=[v^j_0,v^j_1]$, $E^j_1=[v^j_1,v^j_2]$, and $E^j_2=[v^j_0,v^j_2]$. If $*_{b_j}=+ 1$, this ordering is compatible with the cyclic ordering induced by the opposite of the boundary orientation of $F^j_3$. If $*_{b_j}=- 1$ this is no longer true, and since the coordinates $u^j_r$, $r\in \{0,1,2\}$, of $G(T)$ are the cross-ratio moduli of an isometry class of oriented hyperbolic tetrahedron associated to $\Delta_j$, $u^j_r$ should be replaced by $(u^j_r)^{-1}$ in order to compensate the choice of opposite orientation. Hence we define a new system of coordinates $(w^j_r)\in \Cc^{3s}$ of $G(T)$ by labelling $E^j_r$ with $w^j_r$, setting
\begin{equation}\label{crossw}
(w^j_0,w^j_1,w^j_2)= \left\lbrace\begin{array}{ll} (u^j_0,u^j_1,u^j_2)  
&  {\rm if} \ *_{b_j}=+1 \\ (1/u^j_2,1/u^j_1,1/u^j_0) & {\rm if} \ *_{b_j}=-1 .
\end{array}\right. 
\end{equation}
This gives us the model of $G(T)$ that we are going to use. We denote it by $G(T,\tilde b)$. The defining equations of $G(T,\tilde b)$ are the edge relations, for all edges $e$ of $T$, given by 
\begin{equation}\label{gleq}
\prod_{E\to e} w(E)^{*_E}=1
\end{equation}
where if $E$ is an edge of $(\Delta_j,b_j)$ we have $w(E)=w^j_i$ if
and only if $E$ is $E^j_i$ or the opposite edge, and $*_E:=*_{b_j}$. The
volume function on $G(T,\tilde b)$ takes now the form $\textstyle  {\rm
  Vol}(w)=\sum_j *_{b_j} {\rm D}_2(w^j_0)$. Clearly ${\rm Vol}(u)= {\rm Vol}(w(u))$.
\medskip
\begin{exa}\label{GT1}  {\rm It is easy to recover the edge equations 
of $G(T,\tilde b)$ from a decoded $\Nn$-graph representing it: at each essential crossing one places the
    cross-ratio variables $(w_0,w_1,w_2)$ like the decorations $(d_0,d_1,d_2)$ in Figure \ref{I-dec}, and take the products of cross-ratio variables along the boundary lines of $N(P)$. For example, consider the cusped manifold $M_1$ and $(T_1,\tilde b)$. Assign $(w_0,w_1,w_2)$ to the top crossing of Figure \ref{I-F8S}, and $(W_0,W_1,W_2)$ to the bottom one. Note that $*_b=+1$ for both of them. Then we get the equations
\begin{equation}\label{eqF8S}
w_0w_1^2W_0W_1^2=1, \ \
w_0w_2^2W_0W_2^2 =1 . 
\end{equation}
 Using the relation
$w_{j+1}=1/(1-w_j)$ and the similar one for $W_j$, they reduce to the unique quadratic equation
$$ w_1(w_1-1)W_1(W_1-1)=1 .$$
The parameter space of ``positive solutions'' is the
half plane ${\rm Im}\ w_1>0$ with the ray $0.5 + si$ removed, where $ \sqrt{15}/2 \leq
s< +\infty$ \cite{FM}. The complete hyperbolic structure is realized at $w_1=W_1= \exp(i\pi/3)$ (two regular ideal tetrahedra). Similarly, using $(T_0,b)$ one recovers Thurston's
celebrated treatment of the figure-eight knot complement $M_0$.}
\end{exa}

\section{Analytic configurations}\label{I-CAC}

\subsection{Local analytic configurations}\label{LocalAC} 
Take an oriented $3$-simplex $(\Delta,b)$. As in the previous section the $2$-face $F_3$ is opposite to the vertex $v_3$, and the edges of
$F_3$ are ordered as
$$ E_0=[v_0,v_1], \ E_1=[v_1,v_2], \ E_2=[v_0,v_2] . $$  

\subsubsection{Quantum hyperbolic $3$-simplices} \label{qhsimplices}
Recall that the edge decorations of $(\Delta,b)$ are equal on opposite edges, and hence specified by triples $d=(d_0,d_1,d_2)$ where $d_r=d(E_r)$ ($r\in \{0,1,2\}$).  A {\it quantum hyperbolic $3$-simplex} is a tuple $(\Delta,b,w,f,c)$ where:
\begin{itemize}
\item $w=(w_0,w_1,w_2)$ is the system of cross ratios \eqref{crossw}.

\item $f=(f_0,f_1,f_2)\in \Z^3$ is such that $l_k = \log(w_k)+ if_k\pi$, $k\in\{0,1,2\}$, satisfies $l_0+l_1+l_2 =0$.
 \item $c =(c_0,c_1,c_2)\in \Z^3$ satisfies $ c_0+c_1+c_2 =1$. 
\end{itemize}
Here, $\log$ is the branch of logarithm with imaginary part in $(-\pi,\pi]$. For every $k\in\{0,1,2\}$ set
\begin{equation}\label{ql-b}
l_{k,N,*_b,c} =  \frac{1}{N}(\log(w_k)+i\pi(N+1)(f_k-*_bc_k)) .
\end{equation} 
We call $f$ a {\it flattening}, $c$ a {\it charge}, $l_k$ a {\it classical log-branch}, and $l_{k,N,*_b,c}$  a {\it quantum log-branch}. If Im$(w)>0$ (resp. Im$(w)<0$), then $c$ is a charge if and only if $-c$ (resp. $c$) is a flattening.  

We organize the edge decorations $w$, $f$, $c$ as follows. Denote by 
$p_\infty : \Ww_\infty \ra \Cc$ the maximal abelian covering map. We realize $\Ww_\infty$ as the quotient space $\Ww_\infty = (\Dd \times \Z^2)/\sim$, where $\Dd$ is the result of gluing two copies of $(-\infty,0)$ and $(1,+\infty)$ to the boundary of $\C \setminus ((-\infty,0)\cup (1,+\infty))$, and $\sim$ is the equivalence relation 
\begin{equation}\label{identrel}
\begin{array}{cc} (x+i0;p,q)\sim (x-i0;p+2,q) & {\rm if} \ 
\ x \in (-\infty,0)\\  
(x+i0;p,q)\sim (x-i0;p,q+2) & {\rm if} \ \ x \in (1, +\infty) .
\end{array}\end{equation}
Setting $l(y;m):=  \log(y)+ im\pi$, $m\in \mz$, the bijective map
\begin{equation}\label{identRS}
 (x; \frac{l(x;p)-\log(x)}{i\pi}, 
\frac{l((1-x)^{-1};q)-\log((1-x)^{-1})}{i\pi}) \to \left(l(x;p),l((1-x)^{-1};q)\right)
\end{equation}
identifies $\Ww_\infty$ with the Riemann surface of the maps
\begin{equation}\label{epsmaps}
\phi_{\epsilon,\epsilon'}:\ x\to (\log(x)+ i\epsilon \pi, \log((1-x)^{-1})
+ i\epsilon \pi)\quad ,\ \epsilon, \epsilon' \in \{0,1 \}.
\end{equation}
For every sign
$*=\pm 1$ and couple $c=(c_0,c_1)\in \Z^2$ define the analytic map $ l _{N,*,c} :  \Ww_\infty  \to (\C^*)^2$ by
$$ l _{N,*,c} ([x;p,q])=
(\frac{1}{N}(\log(x)+i\pi(N+1)(p-*c_0)),\frac{1}{N}(\log((1-x)^{-1})
+i\pi(N+1)(q-*c_1))) . $$
Then: 
\begin{itemize}
\item $\Cc$ (resp. $\Ww_\infty$) with coordinate $w_0$ (resp. $[w_0;f_0,f_1]$) is a parametrization of the set of edge decorations $w$ (resp. $(w,f)$) of $(\Delta,b)$.
\item Setting $c_2=1-c_0-c_1$, the components of $l_{N,*_b,c} ([w_0;f_0,f_1])$ determine the three quantum log-branches $l_{k,N,*_b,c}$ of $(\Delta,b,w,f,c)$.
\item Put ${\rm EXP}(z_1,\dots,z_{2s})=
(\exp(z_1),\dots,\exp(z_{2s}))$, $s\geq 1$. The map $ {\rm EXP} \circ  l _{N,*,c}: \Ww_\infty  \to  (\C^*)^2$ provides a distinguished system of
$N$th-roots of the cross ratios $w_k$:
\begin{equation}\label{N_root}
w_k':= \exp(l_{k,N,*_b,c}([w_0;f_0,f_1]). 
\end{equation}
\end{itemize}
The proof of the following Lemma is straightforward (see \cite{GT} and \cite{AGT}, Remark 2.1 and 2.18).
\begin{lem}\label{N-roots} (1) The $N$th-roots $(w_0',w_1',w_2')$ 
  depend only on the residues mod $(N)$ of the flattening $f$ and the
  charge $c$, and verify the relation
\begin{equation}\label{relloc}
w'_0w'_1w'_2 =-\zeta^{*_bm}
\end{equation}
where as usual $N:=2m+1$ and $\zeta := \exp (2i\pi/N)$.
  
(2) For any charge $c$ and $N$-th roots $u_k$ 
of the cross ratios $w_k$ verifying the relation \eqref{relloc}, there
is a flattening $f$ such that $u_k = w_k'$ as in
\eqref{N_root}.

(3) The image $\Ww_N$ of $\Ww_\infty$ via the map ${\rm EXP} \circ l
_{N,*_b,c}$ is the curve in $(\C^*)^2 $ with defining equation
$$u_0^N +(u_1^{-1})^N=1 . \ $$  
There is a natural $N^2$-to-$1$ rational regular map $p_N: \Ww_N \to
\Cc$, $p_N(u_0,u_1):= u_0^N$.
\end{lem}
\subsubsection{Tetrahedral tensors} \label{tettens} 
\begin{figure}[ht]
\begin{center}
 \includegraphics[width=7cm]{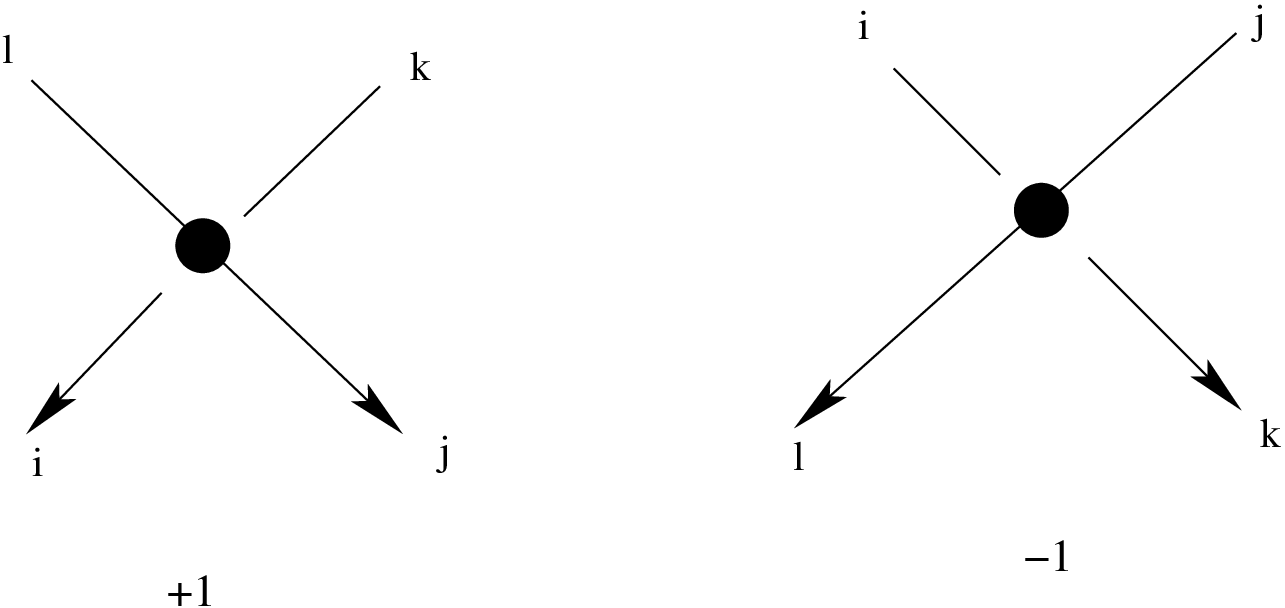}
\caption{\label{tensorindex} Position of tensor indices.}
\end{center}
\end{figure}

Given a branched tetrahedron $(\Delta,b)$ and an operator $A\in {\rm End}(\C^N\otimes \C^N)$, associate a copy $V_j$ of $\C^N$ to
the $2$-face $F_j$, and consider $A$ as a map
\begin{equation}\label{mat}
A = \left\lbrace 
\begin{array}{ll} (A^{i,j}_{k,l}) : V_3 \otimes V_1\rightarrow 
V_2 \otimes V_0 & \mathrm{if}\ *_b=+1 \\ 
(A_{i,j}^{k,l}) : V_2 \otimes V_0 \rightarrow 
V_3 \otimes V_1& \mathrm{if}\ *_b=-1.
\end{array}\right.
\end{equation}
Note that the pair of indices $(i,j)$
(resp. $(k,l)$) corresponds to the space $V_2 \otimes V_0$
(resp. $V_3 \otimes V_1$), and to the outgoing
(resp. ingoing) $2$-face co-orientations induced by $b$ if $*_b=1$, and the converse if $*_b=-1$. The entries $A^{i,j}_{k,l}$ and $A_{i,j}^{k,l}$ are taken in the standard basis of $\C^N\otimes \C^N$. Graphically this is encoded in Figure \ref{tensorindex} (compare with Figure \ref{N-crossing}).

Using \eqref{mat} we define the tetrahedral tensors ${\rm R}_N(\Delta,b,w,f,c)$ of a quantum hyperbolic $3$-simplex by
$$\Rr_1(\Delta,b,w,f):= \Rr_{(1,*_b)}([w_0;f_0,f_1]) \quad \quad {\rm if} \ N=1$$
and
$$\Rr_N(\Delta,b,w,f,c):= \Rr_{(N,*_b,c)}\circ {\rm EXP} \circ l _{N,*_b,c}([w_0;f_0,f_1]) \quad \quad {\rm if} \ N\geq 3 \ {\rm is odd}$$
where $\Rr_{(1,*)}: \Ww_\infty \to {\rm GL}(\C\otimes \C)\cong \C^*$ and $\Rr_{(N,*,c)}: \Ww_N \to {\rm GL}(\C^N\otimes \C^N)$ are defined as follows. The map $\Rr_{(1,*)}$ is analytic. It is given by
 \begin{equation}\label{formR}
 \Rr_{(1,*)} ([x;p,q]):= \exp\left(*\frac{2}{i\pi}\left({\rm L}(x) + 
\frac{i\pi}{2}
(p \log(1-x) + q \log(x))\right)\right)
\end{equation}
where L is the Rogers dilogarithm 
$${\rm L}(x) = -\frac{\pi^2}{6}-\frac{1}{2}\int_0^x
\frac{\log(t)}{1-t}+ \frac{\log(1-t)}{t}\ dt.$$ 
The maps $\Rr_{(N,*,c)}$ are regular rational maps called {\it matrix dilogarithms}. They verify non commutative versions of the five term relations satisfied by L. For every $x\in \C^*$, put
$x^{1/N} := \exp(\log(x)/N)$, extended to $0^{1/N}:=0$ by
continuity. Set
$$[x]:=N^{-1}\frac{1-x^N}{1-x},\quad g(x) := \prod_{j=1}^{N-1}(1 -
x\zeta^{-j})^{j/N},\quad h(x) := g(x)/g(1). $$ 
Let $\Ff_N= \{(u',v') \in \mc^2 \ \vert \ (u')^N +
(v')^N = 1\}$, and $\omega : \Ff_N\times \mz/N\mz\ra \mc$ the regular rational map given by $\omega(u',v'\vert 0) := 1$ and $\textstyle \omega(u',v'\vert n) := \prod_{j=1}^{[n]_N} v'(1-u'\zeta^j)^{-1}$. Define $\Ll_N: \Ff_N \to  {\rm Aut}(\C^N\otimes \C^N)$ by (recall that $m=(N-1)/2$)
\begin{equation}\label{basicf}
 \Ll_N(u',v')_{k,l}^{i,j} =  h(u')\
\zeta^{kj+(m+1)k^2}\ \omega(u',v'\vert i-k) \ \delta_N(i + j - l).
\end{equation} 
We have
\begin{align*}
\bigl( \Ll_N(u',v')^{-1}\bigr)^{k,l}_{i,j}& =  \frac{[u']}{h(u')}\
\zeta^{-kj-(m+1)k^2}\ \frac{\delta_N(i+j-l)}{\omega(u'/\zeta,v'\vert
i-k)}.
\end{align*}
We put
\begin{align}\label{refR2} \Rr_{(N,+1,c)}(u_0,u_1)^{i,j}_{k,l} := &
\bigl(u_0^{-c_1}u_1^{c_0}\bigr)^{\frac{N-1}{2}}
\Ll_N(u_0,u_1^{-1})_{k,l}^{i,j}\\
\Rr_{(N,-1,c)}(u_0,u_1)^{k,l}_{i,j} := &
\bigl(u_0^{-c_1}u_1^{c_0}\bigr)^{\frac{N-1}{2}}
\bigl( \Ll_N(u_0,u_1^{-1})^{-1}\bigr)^{k,l}_{i,j} .\notag
\end{align}
\begin{remark} {\rm The map $\Rr_{(1,*)}$ does not depend on the charge $c$ and $\Rr_{(N,*,c)}$ depends on $c$ only by the factor
$\bigl(u_0^{-c_1}u_1^{c_0}\bigr)^{(N-1)/2}$, but $\Rr_N$ depends on $c$ also via the map $ l _{N,*_b,c}$.}
\end{remark}
\begin{defi} The local analytic configurations over $(\Delta,b,c)$ are the family of spaces and maps indexed by the odd integers $N\geq 3$, given by $$\Aa_N(\Delta,b,c):= \{\Ww_\infty,\Ww_N,p_\infty,p_N, l _{N,*_b,c},\Rr_{(1,*_b)}, \Rr_{(N,*_b,c)}\}.$$
\end{defi}
 
\subsection{Globalization}\label{glob} In addition to the tetrahedral tensors ${\rm R}_N$, when dealing with arbitrary triangulations we need also {\it face tensors} $\Qq_N$, defined as follows. For every odd $N\geq 3$, consider the symmetric $N\times N$ matrices $S$ and $T$ with entries $$S^i_j =N^{-1/2}\zeta^{ij}\ ,\ T^i_j = \zeta^{i^2(m+1)}\delta_N(i+ j).$$ We show in Lemma \ref{projrep} that $S$ and $T$ generate an $N$-dimensional projective representation of $SL(2,\mz)$. We put $\Qq_N:=S\cdot T^{-1}$ (see Section \ref{formal-conv} for the notation ``$\cdot$"). This $N\times N$ matrix has entries 
$$(\Qq_N)^i_j = \sum_{k=0}^{N-1}  S^k_j (T^{-1})^i_k = N^{-\frac{1}{2}}\zeta^{-ij-(m+1)j^2},$$ and is projectively of order $3$, namely $\Qq_N^3= \phi_N^{-1}{\rm I}_{N}$ with $\textstyle \phi_N = \left(\frac{m+1}{N}\right)$ if $N \equiv 1$ mod$(4)$, and $\textstyle \phi_N = \left(\frac{m+1}{N}\right) i$ if $N \equiv 3$ mod$(4)$.

Let $(T,\tilde b)$ be a weakly branched triangulation of $\hat V$ made of $s$ branched tetrahedra $(\Delta_j,b_j)$. For every collection of charges $c^j$ on the $\Delta_j$s and every point $[w;f]:=([w_0^1;f_0^1,f_1^1],\ldots,[w_0^s;f_0^s,f_1^s])\in \Ww_\infty^s$ we call $c$ a {\it rough global charge} and $(T,\tilde b,w,f,c)$ a {\it rough QH triangulation}. Working dually, consider the oriented graph $\Nn$ associated to $(T,\tilde b)$ as before Lemma \ref{r-b}. For every odd $N\geq 1$ define a {\it QH tensor network} over $\Nn$ as follows:
\begin{itemize}
\item associate ${\rm R}_N(\Delta_j,b_j,w^j,f^j,c^j)$ to the structured vertex $v_j$ of $\Nn$ dual to $(\Delta_j,b_j)$;
\item if $N\geq 3$, associate $\Qq_N^{r(e)}$ to each edge $e$ of $\Nn$ with $\mz/3\mz$-color $r(e)$, so that the domain and target of the linear map associated to $\Qq_N^{r(e)}$ correspond respectively to the initial and final endpoints of $e$; if $N=1$, replace $\Qq_N^{r(e)}$ by the scalar $1 \in \mc^*$.
\end{itemize}
Next we define a normalization factor $a_N(T,\tilde b)$. Formally put $a_1(T,\tilde b)=1$. Denote by $q(T,\tilde b)$ the number of edges $e$ of ${\rm Sing}(P)$ such that $r(e) = 2\in \Z/3\Z$. Fix an auxiliary orientation on every edge $e$ of $T$. Let $n_+(e)$
(resp. $n_-(e)$) be the number of abstract {\it diagonal} edges $E\to e$ such that the $\tilde b$-orientation of $E$ agrees (resp. does not agree) with the orientation of $e$. By Lemma \ref{even}, $n_+(e)-n_-(e)$ is an even integer. Set $\epsilon_N := (-1)^{\frac{N-1}{2}}$. For every odd $N\geq 3$ define
\begin{equation}\label{scalfact}
a_N(T,\tilde b):=  N^{-v}\phi_N^{-q(T,\tilde b)}c_N(T,\tilde b)\ ,\ c_N(T,\tilde b) := \epsilon_N^{v+l-\frac{1}{2}\sum_e (n_+(e)-n_-(e))}
\end{equation}
where $\phi_N$ is as above, and $l$, $v$ are respectively the number of edges and the number of vertices of $T$ which are manifold points.
\begin{remark}{\rm Clearly $c_N(T,\tilde b)$ does not depend on the choice of the auxiliary edge orientations, and $ c_N(T,\tilde b)=1$ if $N-1 \equiv 0$ mod$(4)$. If $\tilde b$ is a branching $b$, using the edge $b$-orientation we see that $\textstyle l-\frac{1}{2}\sum_e (n_+(e)-n_-(e))= l- t = \chi(\bar V)=\chi(P)=-v$, where $\chi$ denotes the Euler characteristic and $t$ the number of tetrahedra of $T$. Hence in thise case $c_N(T,b)=1$, and $a_N(T,b)=N^{-v}$. }
\end{remark}
A {\rm $N$-state} of $\Nn$ is an assignment of a label in $\{0,\dots,N-1\}$ to each
edge; every $N$-state $\sigma$ determines a matrix element, denoted by ${\rm R}_N(\Delta_j,b_j,w^j,f^j,c^j)_\sigma$ or $\Qq_{N,\sigma}^{r(e)}$, of each tensor of the QH tensor network. 

\begin{defi}\label{defiSS} {\rm The {\it QH state sum function} $\Hh_N(T,\tilde b,c): \Ww_\infty^{s} \ra \mc$ is the total contraction of the QH tensor network over $\Nn$, normalized by $a_N(T,\tilde b)$:}
\begin{equation}\label{SS}\textstyle 
\Hh_N(T,\tilde b,c)([w;f]):= a_N(T,\tilde b)
\sum_{\sigma}\prod_{j}  {\rm R}_N(\Delta_j,b_j,w^j,f^j,c^j)_\sigma \prod_{e} \Qq_{N,\sigma}^{r(e)}. 
\end{equation}
%J'ai dŽplacŽ le facteur SYM dans la section 9.
\end{defi}
Sometimes we denote $\Hh_N(T,\tilde b,c)([w;f])$ by $\Hh_N(\Tt)$, where $\Tt:=(T,\tilde b,w,f,c)$, and when the context is clear, $\Hh_N$ by $\Hh_N(T,\tilde b,c)$. Note that $\Hh_1$ is a product of analytic scalar functions and does not depend on $c$. When $N\geq 3$, $\Hh_N(T,\tilde b,c)= \Hh'_N(T,\tilde b,c) \circ  ({\rm EXP} \circ  l _{N,*_b,c})$ where $\Hh'_N(T,\tilde b,c)$ is a rational regular function defined on the algebraic variety $\Ww^s_N$.
\smallskip

In order to make \eqref{SS} actual invariant state sums we need to restrict $[w;f]\in \Ww_\infty^s$ to a suitable subspace, that we denote $G_0(T,\tilde b)_\infty$ below (later we will do similarly for charges, see Definition \ref{globalc} and \ref{globalcQHFT}). Consider the model $G(T,\tilde b)$ of the gluing variety, and the products of the maps of Section \ref{LocalAC} over the $\Delta_j$s (keeping, with slight abuse, the same notations):
$$p_\infty: \Ww_\infty^s \rightarrow \Cc^s\ ,\  
p_N: \Ww_N^{s} \to \Cc^s$$
$$l _{N,*,c} :  \Ww_\infty^{s}  \to (\C^*)^{2s}\ ,\ {\rm EXP}\circ l _{N,*,c} :  \Ww_\infty^{s} \to \Ww_N^{s}\subset (\C^*)^{2s}.$$
Set
$$G(T,\tilde b)_\infty := p_\infty^{-1}(G(T,\tilde b)), \ \ \  
G(T,\tilde b,c)_N := {\rm EXP}\circ l _{N,*_b,c}(G(T,\tilde b)_\infty)
. $$ The restrictions of $p_\infty$ and $p_N$ to these spaces are covering maps of $G(T,\tilde b)$, that make a commutative diagram with ${\rm EXP}\circ l
_{N,*_b,c}$; $G(T,\tilde b)_\infty$ is a closed analytic subset
of $\Ww_\infty ^s$, $G(T,\tilde b,c)_N$ an algebraic subset of
$\Ww_N^s$, and ${\rm Sing}(G(T,\tilde b)_*) = p_*^{-1}({\rm Sing}(G(T,\tilde b))$. 

Define the {\it total cross-ratio} of an
edge $e$ of $T$ by (see the notations at the end of
Section \ref{I-triang})
$$W(e):= \prod_{E\to e} w(E)^{*_E}$$
so that $[w;f]\in G(T,\tilde b)_\infty$ if and only if $W(e)=1$ for
all $e$. Similarly, define the {\it total log-branch},
{\it total $N$-th root modulus}, and {\it total charge} at $e$ by
\begin{equation}\label{edgeconst}
L(e):= \sum_{E\to e} *_El(E)\ ,\ W'_N(e) := 
\prod_{E\to e} w'(E)^{*_E}\ ,\ C(e) := \sum_{E\to e} c(E)
\end{equation}
where $l(E)= \log(w^j_k)+ if^j_k\pi$, $w(E)=w^j_k$ and $c(E)=c^j_k$ if
and only if $E$ is an edge of $(\Delta_j,b_j,w^j)$, $E$ is either $E^j_k$ or the opposite edge, and
$*_E=*_{b_j}$.  Fix an auxiliary ordering of the, say $n$, edges of
$(T,\tilde b)$. We have analytic maps
 $$t_L: \Ww_\infty^s \to \C^n, \ \ t_L([w;f]):=(L(e_1),\dots L(e_n))$$
 $$t_{N,W',c}: \Ww^s_\infty \to \C^n, \ \ t_{N,W',c}([w;f]):= 
 (W'_N(e_1),\dots W'_N(e_s)) $$ the dependence of $L(e_k)$
 (resp. $W'_N(e_k)$) on $[w;f]$ (resp. $[w;f]$ and $c$) being
 understood. The following Lemma is evident.
\begin{lem}\label{t|GT} The restriction of $t_L$ to $G(T,\tilde b)_\infty$
  takes values in $(2i\pi\Z)^n$ and the restriction of $t_{N,W',c}$
  takes values in the set $\mu_N$ of
  $N$th-roots of unity.  Both maps are discrete-valued, hence constant
  on the closure of each connected component of the set of
  non-singular points of $G(T,\tilde b)_\infty$. Moreover, the
  restriction of $t_{N,W',c}$ factorizes as
$$ t_{N,W',c}=t'_{N,W}\circ ({\rm EXP}\circ l _{N,*_b,c})$$
where $t'_{N,W}$ is a rational regular map defined on $G(T,\tilde
b,c)_N$ and constant on irreducible components.
\end{lem}
Finally define the analytic subset
$$G_0(T,\tilde b)_\infty = t_L^{-1}(0,\dots,0)\subset  G(T,\tilde b)_\infty 
$$
and, for every $N$, the algebraic subvariety
$$G_0(T,\tilde b,c)_N = {\rm EXP}\circ l _{N,*_b,c}(G_0(T,\tilde b)_\infty)
\subset G (T,\tilde b,c)_N . $$
As usual, keep the same notation for the restriction of any map already
defined. 
\begin{defi}\label{defiAC}
{\rm The $N$-th {\it analytic configuration} over $(T,\tilde b,c)$
is the family of spaces and maps}
$$\Aa_N(T,\tilde b,c):= 
\{G_0(T,\tilde b)_\infty,G_0(T,\tilde b,c)_N,p_\infty,p_N, l _{N,*_b,c},
\Hh_1, \Hh_N\} . $$
\end{defi}
We stress that every map $t_{N,W',c}$ is constant on $G_0(T,\tilde
b)_\infty$, and its value depends on the fixed rough charge $c$,
precisely on the total edge charges $(C(e_1),\dots, C(e_n))$.
\section{Cusped manifolds}\label{I-cusped}
Let $M$ be a cusped manifold, and $V$ be as in Section
\ref{I-QHI-cusped}, that is, a compact oriented connected $3$-manifold
with one torus boundary component such that $M$ is diffeomorphic to
the interior of $V$. We use the notation $\hat V$ as in Section
\ref{I-triang}.

 \subsection {The augmented $PSL(2,\C)$-character variety and
  $A$-polynomial} \label{A-pol}
Fix a geometric basis $(l,m)$ of the fundamental group
$\pi_1(\partial V) \cong \Z \times \Z$. It is given by a couple of
oriented simple closed curves on $\partial V$ which meet at one point,
transversely and positively.  If $M$ is the complement of a hyperbolic
knot $K$ in $S^3$ we can take a canonical longitude $l$ and a
meridian $m$ of $K$.
 
There is a conjugacy class $[\sigma]$ of representations of
$\pi_1(\partial V)$ in $PSL(2,\C)$, each one with images isomorphic to
$ \Z/2\Z \times \Z/2\Z$, acting on the Riemann sphere $\C
\mathbb{P}^1=\partial \mh^3$ without a common fixed point. The class
$[\sigma]$ is obtained as follows: pick two geodesic lines in $\mh^3$,
say $\gamma_l$ and $\gamma_m$, that meet at one point forming a right
angle; then define $\sigma$ by $\sigma(x)=r_x$, where $x\in\{l,m\}$ and
$r_x$ is the rotation by $\pi$ around $\gamma_x$.

By definition, the $PSL(2,\C)$-character variety of $M$ is the
algebro--geometric quotient
$$ X'(M) := 
{\rm Hom}' (\pi_1(M),PSL(2,\C))/\!/ PSL(2,\C) $$ where $PSL(2,\C)$ acts by 
conjugation, and $$ {\rm Hom}' (\pi_1(M),PSL(2,\C))\subset
{\rm Hom}(\pi_1(M),PSL(2,\C))$$ is the subset of representations that do
not restrict to an element of $[\sigma]$ on $\pi_1(\partial V)$. 
Given a peripheral subgroup $\pi_1(\partial V)$ of
$\pi_1(M)\cong \pi_1(V)$, denote by
$$R(M)\subset {\rm Hom}(\pi_1(V),PSL(2,\mc)) \times \mc \mathbb{P}^1$$ 
the set of couples $(r,z)$ such that $z$ is fixed by $r(\pi_1(\partial
V))$. In particular, couples such that $[r]$ restricts to $[\sigma]$
on $\partial V$ are excluded. The {\it augmented $PSL(2,\C)$-character
  variety} of $M$ is the algebro-geometric
quotient \begin{displaymath} X(M) := R(M) /\!/ PSL(2,\mc)
\end{displaymath} 
where $PSL(2,\mc)$ acts by conjugation on $R(M)$ and by M\"obius
transformations on $\mc \mathbb{P}^1$. Hence, any {\it augmented
  character} of $M$, $\rho=[(\rho,z)] \in X(M)$, is a character of
representations $r : \pi_1(M) \rightarrow PSL(2,\mc)$ together with a
choice of fixed points of the peripheral subgroups, invariant by
conjugation. The hyperbolic holonomy of $M$ defines a point $\rho_{\rm
  hyp}\in X(M)$.  In a similar way we can define the
$PSL(2,\C)$-character variety $X'(\partial V)$ of $\partial V$ and its
augmented version $X(\partial V)$. The inclusion $i:\partial V \to V$
induces regular maps $i^*: X'(M)\to X'(\partial V)$ and $i^*: X(M)\to
X(\partial V)$.  
\begin{teo}\label{XM_property}
  (1) (\cite{CS,BZ}; see \cite{Shalen}) Both $X'(M)$ and $X(M)$ are
  complex algebraic affine varieties, and the natural projection $q:
  X(M)\rightarrow X'(M)$ is a regular map.

(2) (\cite{Kapovich}) 
$\rho_{\rm hyp}$ (resp. $q(\rho_{\rm hyp})$) is a regular point of
$X(M)$ (resp. $X'(M)$). Hence it belongs to a unique irreducible
component $X_0(M)$ of $X(M)$ (resp. $X'_0(M)$ of $X'(M)$), which is a
complex algebraic curve.

(3) The restriction of $q$ to $X_0(M)$, $q_0: X_0(M)\to X'_0(M)$, is
generically $2:1$.

(4) (\cite{Dun}) The restricted map $i^*: X'_0(M)\to X'_0(\partial V)$
(resp. $i^*: X_0(M)\to X_0(\partial V)$) is generically $1:1$. Hence
it is a birational isomorphism onto its image.
 \end{teo}
 For any $\rho \in X(\partial V)$ which is non trivial, let
 $\bar{\rho}$ be a representative of $\rho$ such that $\bar
 \rho(\pi_1(\partial V))$ fixes $z=\infty \in \mc \mathbb{P}^1$. For
 any non zero class $\gamma \in \pi_1(\partial V)$, $\bar
 \rho(\gamma)$ acts on $\mc$ as $w\mapsto \gamma_\rho w +
 b$, where $\gamma_\rho \in \mc^*$ and $b\in \mc$. In general
 $\gamma_\rho$ is a squared eigenvalue of $\bar \rho(\gamma)\in
 PSL(2,\C)$; if $\bar \rho(\gamma)$ is loxodromic, then
 $\gamma_\rho\ne 1$, and the two reciprocally inverse eigenvalues are
 distinguished by the augmentation, which selects an endpoint, whence
 an orientation, of $\gamma$. Consider the so called {\it holonomy
   map}
\begin{equation}\label{loghmap1}
\begin{array}{lcll}
{\rm hol}_\gamma : & X(\partial V) & \longrightarrow & \mc^*\\
           & \rho & \longmapsto & \gamma_\rho.
           \end{array}
\end{equation}

By using the above cusp basis $(l,m)$ we get an algebraic isomorphism
$${\rm hol}_m \times {\rm hol}_l : X(\partial V) \to \mc^* \times \mc^* . $$
Define the rational map
$$ \hG: X(M)\to \mc^* \times \mc^*, \ \ \hG= ({\rm hol}_m 
\times {\rm hol}_l)\circ i^* . $$ Following \cite{CCGLS},
\cite{Cha}, and Dunfield's appendix in \cite{BD}, consider the plane
curve $A(M)$ defined as the closure of the $1$-dimensional part of the
image $\hG(X(M))$.  The (suitably normalized) polynomial generating
the ideal of $A(M)$ is by definition the {\it $PSL(2,\mc)$
  $A$-polynomial} of $M$. We denote by $A_0(M)$, and call {\it
  geometric component of $A(M)$}, the closure of $\hG(X_0(M))$. It
follows from Theorem \ref{XM_property} (4) that:
 \begin{cor}\label{2property} 
   The restricted map $\hG : X_0(M) \ra A_0(M)$ is a birational
   isomorphism.
  \end{cor} 
 \subsection{Rich components of $G(T)$}
\label{rich}
Let $T$ be any ideal triangulation of $\hat V$ such that the gluing
variety $G(T)$ is non empty. There is a natural regular map (see eg. \cite{Cha} or the
appendix of \cite{BD})
\begin{equation}\label{rhomap}
\rho: G(T) \to X(M) . 
\end{equation}
As $G(T)$ is a complex algebraic curve, $\hG(\rho(G(T))$ is a union of
irreducible components of the plane curve $A(M)$.  Assume that $G(T)$
contains a {\it non-negative} point $z_h$ such that
$\rho(z_h)=\rho_{\rm hyp}$. Recall that ``non-negative'' means that
for every tetrahedron $\Delta_j$ of $T$ the cross-ratios
$z_h^j$ have non-negative imaginary parts.  The point $z_h$ is not
necessarily a regular point of the gluing variety, hence in general it
could be contained in several irreducible components of $G(T)$.
\begin{defi}\label{rich}{\rm An irreducible component of $G(T)$ is
    {\it rich} if it contains $z_h$ and an infinite sequence of points
    $z_n$ that converge to $z_h$ (in the ``strong'' topology of $G(T)$
    as an analytic space), and correspond to compact closed hyperbolic
    Dehn fillings of $M$ that converge geometrically to $M$.}
\end{defi}

\begin{prop}\label{rich-exist} 
  For every gluing variety $G(T)$ containing a non-negative point
  $z_h$ with $\rho(z_h)=\rho_{\rm hyp}$, the set of rich components of
  $G(T)$ is non empty and finite.
\end{prop}
\noindent {\it Proof.} 
This follows from the proof of the hyperbolic Dehn filling Theorem
given in \cite{PP}, and the fact that the algebraic curve $G(T)$ has a
finite number of irreducible components. \cvd

\begin{remark}\label{onHDST} {\rm In \cite{NZ}, the hyperbolic Dehn
    filling Theorem was proved assuming that $M$ (allowing several
    cusps) has an ideal triangulation $T$ such that $G(T)$ contains a
    {\it strictly positive} point $z_h$ representing the hyperbolic
    holonomy. This proof uses also that $z_h$ is a regular point of
    $G(T)$. By elaborating on the analysis of \cite{NZ}, the paper
    \cite{Choi} has provided a proof of this regularity result,
    extended to every strictly positive point of $G(T)$. In \cite{BP},
    the proof of \cite{NZ} was presented with some mild modifications
    avoiding the assumption that $z_h$ is regular. Unfortunately, it did not give yet a complete proof of the
    hyperbolic Dehn filling Theorem based on gluing varieties, since
    it is not known if every $M$ has strictly positive geometric ideal
    triangulations. The paper \cite{PP} filled this gap; it uses any
    gluing variety $G(T)$ with a non-negative and possibly singular
    point $z_h$, which always exists; the proof uses also the arguments of
    \cite{BP}. In \cite{LUO1} it is proved that
    strictly positive geometric ideal triangulations exist
    ``virtually'', that is, every cusped $M$ has a finite covering
    $\tilde M$ having such an ideal triangulation. On the other hand,
    passing to $\tilde M$ increases the number of cusps.  }
 \end{remark}    
\begin{prop}\label{birational} Let $Z$ be any rich component
  of $G(T)$. Then the restricted maps $\rho: Z\to X_0(M)$ and
  $\hG\circ \rho: Z\to A_0(M)$ are generically $1:1$ and hence are
  both birational isomorphisms. More precisely, there are maximal non
  empty Zariski open subsets $\Omega^0_Z$ and $\Omega_Z$ of $Z$, with
  $\Omega^0_Z\subset \Omega_Z$, such that:
 \begin{enumerate}
 \item $\rho(\Omega_Z)$ and $\rho(\Omega^0_Z)$ are Zariski open subets
   of $X_0(M)$ (resp. $\hG\circ \rho(\Omega_Z)$ and
   $\hG\circ\rho(\Omega^0_Z)$ are Zariski open subsets of $A_0(M)$).
 
 \item The restriction of $\rho$ and $\hG \circ \rho$ to $\Omega_Z$ is
   $1:1$ onto its image, and $z_h \in \Omega_Z$. The restriction of
   $\rho$ and $\hG \circ \rho$ to $\Omega^0_Z$ is a regular rational
   isomorphism onto its image.
 \end{enumerate}
 \end{prop}
 \noindent {\it Proof.} This goes step by step as the proof of Theorem
 3.1 in \cite{Dun}. The key ingredients are the Gromov-Thurston
 ``volume rigidity'' for closed hyperbolic manifolds, the fact that
 the volume of representations yields a well defined function on the
 normalization of the smooth projective model of $A_0(M)$, and the
 existence on $X_0(M)$ of infinite sequences of points corresponding
 to the holonomies of compact closed hyperbolic Dehn fillings of $M$
 that converge geometrically to $M$. In the present situation this last fact is ensured by the definition of a rich component, and the volume of representations lifts to the function \eqref{volf} on
 $G(T)$. \cvd \medskip
 
 If $T$ has a weak branching $\tilde b$ the
 above discussion can be rephrased in terms of the gluing variety
 $G(T,\tilde b)$ of Section \ref{I-triang}.  Recall that
 a system of cross ratios $w^j$ on the branched tetrahedra
 $(\Delta_j,b_j)$ of $(T,\tilde b)$ is non negative if $z^j=
 (w^j)^{*_j}$ is non negative.
 
 \subsection{Refined analytic configurations for one-cusped
   manifolds}\label{cusped-pattern}
 
Let $(T,\tilde b)$ be a weakly branched ideal triangulation of $\hat V$.  

\begin{defi} \label{globalc} {\rm A {\it global charge} $c$ on
    $T$ is a rough global charge satisfying the following
    additional global constraint:} For every edge $e$ of $T$, the
  total edge charge $C(e)$ is constant and equal to $2$.
\end{defi}

For every global charge $c$ the variety $G_0(T,b,c)_N$ has at every
edge $e$ of $T$ the defining equation
\begin{equation}\label{edge-equation} 
W'(e)=\zeta^{-1}
\end{equation}
where $\zeta := \exp(2\pi i/N)$. Note that $\zeta^{-1}=\zeta^{2m}$. We have:

\begin{prop}\label{top-support}
  (1) Every global charge $c$ on $(T,\tilde b)$ determines a $c$-weight
  $(h_c,k_c)$ of $V$.

  (2) For every $c$-weight $(h_c,k_c)$ of $V$ and every ideal
  triangulation $(T,\tilde b)$ of $\hat V$, there exists a global
  charge $c$ on $(T,\tilde b)$ with $c$-weight equal to $(h_c,k_c)$.
  \end{prop}

  \begin{prop}\label{top-support2} Assume that the gluing variety
    $G(T,\tilde b)$ is non empty, and let $\rho : G(T,\tilde b)
    \rightarrow X(M)$ be as in \eqref{rhomap}.

(1) Every point $[w;f]\in G_0(T,\tilde b)_\infty$ 
determines an $f$-weight $(h_f,k_f)$ of $(V,\rho(w))$.

(2) For every point $\rho\in \rho(G(T,\tilde b))$ and every $f$-weight
$(h_f,k_f)$ of $(V,\rho)$, there exists a point $[w;f]$ of
$G_0(T,\tilde b)_\infty$ with $f$-weight equal to $(h_f,k_f)$.
\end{prop}
Proposition \ref{top-support} and \ref{top-support2} are reformulations of results of Neumann \cite{N0,N}.  We recall below the constructions underlying the
statements (1), which claim the existence of the pairs
$(\gamma(c),\gamma_2'(c))$ and $(L([w;f]),\gamma_2'(f))$ in Definition
\ref{defiweight}. The statements (2) are much more subtle. For the sake of completeness, we recall the proofs at the end of the section (part of the material will be used in Section \ref{qCS}). 

Denote by $T_0$ the cell complex obtained from $T$ by removing an open
cone neighborhood of its vertex. Hence $T_0$ is the result of gluing
tetrahedra with truncated vertices, and the polyhedron underlying
$T_0$ is homeomorphic to $V$. Denote by $\partial T_0$ the induced
triangulation of $\partial V$. Represent any non zero class in
$H_1(\partial V;\mz)$ by {\it normal loops}, that is, a disjoint union
of oriented essential simple closed curves in $\partial V$, transverse
to the edges of $\partial T_0$ and such that no curve enters and exits
a $2$-simplex by a same edge. The intersection of a normal loop, say
$C$, with a $2$-simplex $F$ consists of a disjoint union of arcs, each
of which turns around a vertex of $F$; if $F$ is a cusp section of the
tetrahedron $\Delta$ of $T$, for every vertex $v$ of $F$ we denote by
$E_v $ the edge of $\Delta$ containing $v$, and by $*_v$ the branching
sign of $\Delta$. We write $C \rightarrow E_v$ to mean that some
subarcs of $C$ turn around $v$. We count them algebraically, by using
the orientation of $C$: if there are $s_+$ (resp. $s_-$) such subarcs
whose orientation is compatible with (resp. opposite to) the
orientation of $\partial V$ as viewed from $v$, then we set $ind(C,v)
:=s_+-s_-$. For every point $[w;f]\in G_0(T,\tilde b)_\infty$ and
every global charge $c$ on $(T,\tilde b)$, one defines cohomology
classes $L([w;f])\in H^1(\partial V;\mc)$ and $\gamma(c)\in
H^1(\partial V;\mz)$ by ($C$ is a normal loop on $\partial V$):
\begin{align}
  L([w;f])([C]) & := \sum_{C \ra E_v} *_v \ ind(C,v) l(E_v) \label{weightwf}\\
  & = 
\sum_{C \ra E_v} *_v \ ind(C,v)\left(\log(w(E_v)) +\pi \sqrt{-1}f(E_v)\right) 
\notag\\
  \gamma(c)([C]) & := \sum_{C \ra E_v} ind(C,v) c(E_v).
\end{align} 
Additional classes $\gamma_2'(c)$, $\gamma_2'(f) \in H^1(V;\mz/2\mz)$
are defined similarly, by using normal loops in $T$ and taking the sum
mod$(2)$ of the charges or flattenings, respectively, of the edges we
face along the loops. The reduction mod$(2)$ of $\gamma(c)$ coincides
with the image of $\gamma_2'(c)$ under the map induced on cohomology by
the inclusion $i:\partial V \ra V$. Hence
$r(\gamma(c))=i^*(\gamma_2'(c))$. Also, denoting by $d_w\in
H^1(\partial V; \C/2i\pi \Z)$ the log of the linear part of
the restriction of $\rho(w)$ to $\pi_1(\partial V)$, for all $a\in
H_1(\partial V; \Z)$ we have
\begin{equation}\label{fconstraint2}
L([w;f])(a)= d_w(a) \ {\rm mod}(i\pi)\ ,\
(L([w;f])(a)-d_w(a))/i\pi = i^*(\gamma_2'(f))(a) \ \ {\rm mod} (2) . 
\end{equation}
This proves that the following definition agrees with the defining
constraints of weights (see \eqref{cconstraint} and
\eqref{fconstraint} in the Introduction).
\begin{defi} \label{defiweight}{\rm The $f$-weight $(k_f,h_f)$ of a
    point $[w;f]\in G_0(T,\tilde b)_\infty$, and the $c$-weight
    $(k_c,h_c)$ of a global charge $c$ on $(T,\tilde b)$, are defined
    respectively as the pairs of cohomology classes}
    $$(L([w;f]),\gamma_2'(f)) \in H^1(\partial V;\mc)\times
    H^1(V;\mz/2\mz)$$
    and
    $$(\gamma(c),\gamma_2'(c)) \in H^1(\partial
    V;\mz)\times H^1(V;\mz/2\mz).$$
\end{defi}
\begin{exa}\label{moreF8S} 
  {\rm ({\it Example \ref{GT1} continued}). Every point $w$ of the gluing variety $G(T_1,\tilde b)$ gives rise to a class $ \exp(d_w)\in H^1(\partial V; \C^*)$. For a suitable choice of geometric basis $(\lambda,\mu)$ of $H_1(\partial V; \Z)$ we have
$$ \exp(d_w)(\lambda)=w_0w_1W_0^{-1}W_1^{-1}\ , \  \exp(d_w)(\mu)=w_1^{-2}W_2^2 . $$

Denote by $f_j$ and $l_j:=\log(w_j) + i\pi f_j$ (resp. $F_j$ and $L_j:= \log(W_j) + i\pi F_j$) the flattening and log-branch variables at the top
(resp. bottom) crossing of Figure \ref{I-F8S}. The equations for the log-branch variables are (see \eqref{eqF8S})
\begin{equation}\label{eqlb} 
l_0+l_1+l_2=L_0+L_1+L_2= 0, \ \ l_0+2l_1+L_0+ 2L_1= 0, \ \  
l_0+2l_2+ L_0+ 2L_2 = 0 . 
\end{equation}
For simplicity let us deal with the positive points $w$ of the
(unique) rich component of $G(T_1,\tilde b)$ which
contains the hyperbolic solution, $w_1=W_1= \exp(i\pi/3)$. For every such a
positive point the equations for the flattening variables are
\begin{equation}\label{eqf} 
f_0+f_1+f_2=F_0+F_1+F_2= -1, \ \ f_0+2f_1+F_0+ 2F_1= -2, \ \   
f_0+2f_2+ F_0+ 2F_2 = -2 . 
\end{equation}
Indeed,  under positivity of $w$ the pairs of first two equations in the
systems \eqref{eqf} and \eqref{eqlb} are clearly equivalent one to
each other. Although a bit subtler, this is true also for the
remaining pairs of equations, because the consistency relations
\eqref{eqF8S} imply that the sum of the arguments of the cross ratios
around edges is {\it exactly} equal to $2\pi$ (see \cite{BP},
Lemma E.6.1).

The equations for the charge variables are formally obtained from those for the flattening
ones by replacing in \eqref{eqf} each $f_j$ with $-c_j$ and each 
$F_j$ with $-C_j$.  Solving the systems we get
$$\begin{array}{c} f_0=-(2f_1+F_1-F_2+1)\ , \ F_0=-(F_1+F_2+1)\ , \ f_2=f_1+F_1-F_2\\ c_0=-(2c_1+C_1-C_2-1)\ , \ C_0=-(C_1+C_2-1)\ , \ c_2=c_1+C_1-C_2\\ l_0=-(2l_1+L_1-L_2)\ ,\ L_0=-(L_1+L_2), \ \ l_2=l_1+L_1-L_2.\end{array}$$
Let us turn to the boundary $c$-weight $k_c\in H^1(\partial V;\Z)$ and, for every positive $w$ as above, the 
boundary  $f$-weight $k_f=k_f(w)\in H^1(\partial V; \C)$. They are given by
$$\begin{array}{c} k_c(\lambda)=c_0+c_1-C_0-C_1\ , \ k_c(\mu)= 2C_2-2c_1\\ k_f(\lambda)=l_0+l_1-L_0-L_1\ , \ k_f(\mu)= 2L_2-2l_1 .\end{array}$$
At the complete solution  we have 
$$k_f(\lambda)/i\pi =f_0+f_1-F_0-F_1\ , \ k_f(\mu)/i\pi= 2F_2-2f_1 .$$
Hence $k_f/i\pi\in H^1(\partial V;\Z)$, and $k_f(\mu)/i\pi=0$
mod$(2)$. Moreover
$$f_0=k_f(\lambda)/i\pi-k_f(\mu)/2i\pi-2f_1-1\ , \ f_2=f_1+k_f(\mu)/2i\pi-k_f(\lambda)/i\pi$$ 
$$F_0=k_f(\lambda)/i\pi-3k_f(\mu)/2i\pi -2f_1-1\ ,\  
F_1= k_f(\mu)/i\pi-k_f(\lambda)/i\pi+f_1$$
$$F_2=f_1+k_f(\mu)/2i\pi .$$ 
Hence for any given boundary $f$-weight $k_f$, there is a one-parameter family of flattenings that realize $k_f$, depending on the free parameter $f_1$. Similarly for the charges we have
$$c_0=k_c(\lambda)-k_c(\mu)/2-2c_1+1\ , \ c_2=c_1+k_c(\mu)/2-k_c(\lambda)$$
$$C_0=k_c(\lambda)-3k_c(\mu)/2 - 2c_1+1\ , \ 
C_1=k_c(\mu)-k_c(\lambda)+c_1 \ ,\ C_2=c_1+k_c(\mu)/2$$ 
where the free parameter is now $c_1$. By taking $f_1=c_1=0$, the relevant pairs of parameters that
enter the formulas of the tensors $\Rr_N$ are
$$(f_0,f_1)=(k_f(\lambda)/i\pi-1-k_f(\mu)/2i\pi,0)\ ,\  (c_0,c_1)=(k_c(\lambda)+1-k_c(\mu)/2,0)$$ 
$$(F_0,F_1)=(-1+k_f(\lambda)/i\pi-3k_f(\mu)/2i\pi, k_f(\mu)/i\pi - k_f(\lambda)/i\pi)$$
$$(C_0,C_1)= (1+k_c(\lambda)-3k_c(\mu)/2, k_c(\mu) - k_c(\lambda))$$ 
stressing in this way a pure dependence on the weights. Finally, the state sums are
$$\Hh_N(T_1,\tilde b, w,W,f,F,c,C)=
\sum_{i,j,l,k,I,J=0}^{N-1}
\Rr_N(1,w,f,c)^{i,j}_{k,l}\Rr_N(1,W,F,C)_{j,i}^{I,J} (\Qq^2)^l_I
(\Qq)^k_J$$ where we quoted that $*_b=1$ for both
tetrahedra.\hfill$\Box$}
\end{exa}

Let now $Z$ be a rich component of $G(T,\tilde b)$. Denote $Z_\infty
:= p_\infty^{-1}(Z)$ (an analytic subspace of $G_0(T,\tilde
b)_\infty$), and let $Z_{\infty,0}$ be the union of connected
components of $Z_\infty$ made of points $[w;f]$ such that
\begin{equation}\label{mod2condition}
i^*(\gamma_2'(f)) =0 \in H^1(\partial V;\mz/2\mz).
\end{equation}
Recall the basis $(l,m)$ of $\pi_1(\partial V)$. Because of
\eqref{fconstraint2}, for every point $[w;f]\in Z_{\infty,0}$ we have
\begin{equation}\label{property1} 
(e^{L([w;f])(l)},e^{L([w;f])(m)}) = \mathfrak{h}\circ \rho (w),
\end{equation}
and for every global charge $c$ on $(T,\tilde b)$, the formula 
\begin{align} L_{N,c}([w;f])(C) & := 
\sum_{C \ra E_v}  *_v\ ind(C,v) l_{N,*_v,c}(E_v) \notag\\ &  =  
\sum_{C \ra E_v}  \frac{*_v}{N}\ ind(C,v)\left(\log(w(E_v)) +
\pi \sqrt{-1}(N+1)(f(E_v)-*_vc(E_v))\right)\notag\\
&  =  \frac{1}{N}\ (L([w;f])([C]) - \pi \sqrt{-1}(N+1) 
\gamma(c)([C])) \label{moreclass2}\\
& \hspace*{5cm} + \pi \sqrt{-1} \underbrace{\sum_{C \ra E_v}  
*_v ind(C,v)f(E_v)}_{\in 2\mz} \notag
\end{align}
yields a cohomology class ${L_{N,c}([w;f])} \in H^1(\partial V;\mc/2\pi i \mz)$ such that
\begin{equation}\label{property1bis} 
((e^{L_{N,c}([w;f])(l)})^N,(e^{L_{N,c}([w;f])(m)})^N) =  
\mathfrak{h}\circ \rho (w).
 \end{equation}
If we allow $[w;f]$ to live in $Z_\infty\setminus Z_{\infty,0}$, the
sum in \eqref{moreclass2} can be odd, and the equalities
\eqref{property1} $\&$ \eqref{property1bis} hold true only up to signs
in each component. More precisely, for all $[w;f]\in Z_\infty$ the sum
in \eqref{moreclass2} is an integer expression of the residue class
$\gamma_2'(f)([C])\in \mz/2\mz$, so that we have a well-defined
cohomology class
\begin{equation}\label{ultimateclass} {L_{N,c}([w;f])} \in
H^1(\partial V;\mc/2\pi i \mz)
\end{equation} 
such that
$${L_{N,c}([w;f])} = \frac{1}{N}\ (L([w;f]) - \pi i(N+1)\gamma(c)) 
+ \pi\sqrt{-1}\  i^*(\gamma_2'(f)) \quad {\rm mod}(2\pi i\Z)$$
where, as usual, $i^*:H^1(V;\mz/2\mz) \ra H^1(\partial V;\mz/2\mz)$ is
induced by the inclusion map $i:\partial V \ra V$.  The spaces
$Z_{\infty,0}$ and $A_0(M)$ are related as follows (see Proposition
\ref{birational}). Let
$$A_0(M)_\infty := \pi_\infty^{-1}(A_0(M))\ ,\  
A_0(M)_N := {\rm EXP} \circ l_{N,k_c} (A_0(M)_\infty)$$
where
\begin{equation}\label{infcov}\pi_{\infty} : \tilde{\mc}^2 \ra (\mc^*)^2\ ,
\ ([\lambda;p],[\lambda;q]) \mapsto (\lambda,\mu)
\end{equation}
is the universal covering map, and for every
integral class $k_c\in H^1(\partial V;\mz)$, define $l_{N,k_c}:\tilde{\mc}^2 \ra \mc^2$ by
\begin{multline}\label{lNkc} 
l_{N,k_c}([\lambda;p],[\mu;q]) := 
\left( \frac{1}{N}(\log(\lambda)+2\pi\sqrt{-1}p -\pi\sqrt{-1}(N+1)k_c(l)), 
\right.\\ \left. \frac{1}{N}(\log(\mu)+2\pi\sqrt{-1}q 
-\pi\sqrt{-1}(N+1)k_c(m))\right).
\end{multline}
The restricted map $\pi_\infty : A_0(M)_\infty \rightarrow A_0(M)$ and
the map $\pi_N : A_0(M)_N \ra A_0(M)$, $(u,v)\mapsto (u^N,v^N)$,
define respectively a $\mz\times \mz$- and a $\mz/N\times
\mz/N$-covering of $A_0(M)$.
\begin{defi} \label{defilogholcl} {\rm The {\it classical} and {\it quantum log-holonomies}
    are the maps defined on $Z_{\infty,0}$ by the components of the
    classes $L([w;f])\in H^1(\partial V;\mc)$ and $\exp \circ
    L_{N,c}([w;f])\in H^1(\partial V;\mc^*)$, $[w;f]\in Z_{\infty,0}$,
    in the basis $(l,m)$ of $\pi_1(\partial V)$:}
\begin{align*}
& \fonc{loghol_\infty}{Z_{\infty,0}}{A_0(M)_\infty}
{[w;f]}{(L([w;f])(l),L([w;f])(m)).}\\
& \fonc{loghol_{N,c}}{Z_{\infty,0}}{A_0(M)_N}{[w;f]}
{(e^{L_{N,c}([w;f])(l)},e^{L_{N,c}([w;f])(m)}).}
\end{align*}
\end{defi}
Note that \eqref{property1} $\&$ \eqref{property1bis} imply that
loghol$_\infty$ (resp. loghol$_{N,c}$) maps into $A_0(M)_\infty$
(resp. $A_0(M)_N$).  The identification \eqref{identRS} of
$\Ww_\infty^s$ as a product of Riemann surfaces, and the fact that
$Z_{\infty}\subset G_0(T,\tilde b)_\infty$ is an analytic subspace,
imply that both loghol$_\infty$ and loghol$_{N,c}$ are analytic maps.
\begin{remark}\label{extclass} {\rm If one removes the assumption
    \eqref{mod2condition}, the maps loghol$_\infty$ and loghol$_{N,c}$
    extend to the whole of $Z_\infty$. Then loghol$_\infty$ maps $Z_\infty$ to the union of four connected components $A_0(M)_\infty^{\varepsilon,\varepsilon'} :=
    (\pi_\infty^{\varepsilon,\varepsilon'})^{-1}(A_0(M))$, where
    $\pi_\infty^{\varepsilon,\varepsilon'} : \tilde{\mc}^2 \ra (\mc^*)^2$,
    $([\lambda;p],[\mu;q])\mapsto (\varepsilon \lambda, \varepsilon'
    \mu)$, $\varepsilon, \varepsilon' \in\{-1,+1\}$. The target of loghol$_{N,c}$ becomes the union of the spaces ${\rm EXP} \circ
    l_{N,k_c}(A_0(M)_\infty^{\varepsilon,\varepsilon'})$. The
    component $A_0(M)_\infty = A_0(M)_\infty^{0,0}$ is a natural one:
    we explain in Section \ref{qCS} its relation with the Chern-Simons
    line bundle over $\partial V$.}
\end{remark}
We deduce immediately from the formulas \eqref{weightwf} and \eqref{lNkc}:
\begin{lem} \label{diagram1} For every $N$ and every global charge $c$
  on $(T,\tilde b)$ we have a commutative diagram:
$$\xymatrix{A_0(M) & & A_0(M)_N \ar[ll]_{\pi_N} \\  
& Z \ar[lu]^{\mathfrak{h}\circ \rho} \ar[ld]_{\mathfrak{h}\circ \rho} &  
& Z_{\infty,0} \ar[ll]_{\hspace*{1.5cm} p_\infty}^{\hspace*{-1cm} e^{l_{N,k_c}}} 
\ar[lu]_{loghol_{N,c}} \ar[ld]^{loghol_\infty} & \\ A_0(M) & & A_0(M)_\infty 
\ar[uu] |!{[ur];[ul]}\hole \ar[ll]_{\pi_\infty}}$$
where l$_{N,k_c}$ is defined by \eqref{lNkc}, using the charge boundary weight 
$k_c$, and we denote $EXP \circ l_{N,k_c}$ by $e^{l_{N,k_c}}$.
\end{lem}
Let now $\Ff$ be any finite family of rich components, possibly
contained in different gluing varieties associated to different weakly
branched triangulations of $\hat V$. Recall Proposition
\ref{birational}. Define the following non empty Zariski open subsets
of $X_0(M)$:
$$ \Omega_\Ff= \cap_{Z\in \Ff} \ 
\rho(\Omega_Z)\ ,\ \Omega^0_\Ff= \cap_{Z\in \Ff} \  \rho(\Omega^0_Z).$$
If $Z\in \Ff$ is a rich component  of $G(T,\tilde b)$, 
define
$$\Omega_\Ff(Z):=
  \rho_{\vert Z}^{-1}(\Omega_\Ff)\ ,\ \Omega^{0}_\Ff(Z) := 
\rho_{\vert Z}^{-1}(\Omega^0_\Ff)$$
and
$$\Omega_\Ff(Z)_{\infty,0} 
:= p_\infty^{-1}(\Omega_\Ff(Z))\cap Z_{\infty,0}\ ,\ \Omega_\Ff(Z)_N := {\rm EXP} \circ l_{N,*_b,c}(\Omega_\Ff(Z)_{\infty,0}).$$
We have:
\begin{itemize} 
\item  $\Omega^0_\Ff \subset \Omega_\Ff$. 
\item The restriction of $\rho$ to $\Omega_\Ff(Z)$ is a rational
  regular map which is $1:1$ onto $\Omega_\Ff$. Moreover, the non
  negative point $w^Z_h\in Z$ such that $\rho(w^Z_h)=\rho_{\rm hyp}$
  belongs to $\Omega_\Ff(Z)$.
  \item For every component $Z\in \Ff$, the restriction of $\rho$ to 
$\Omega^{0}_\Ff(Z)$ is a rational regular isomorphism onto $\Omega^0_\Ff$.

\end{itemize}
The following Lemma is an immediate consequence of Proposition
\ref{top-support2}.
\begin{lem}\label{prelemma} Let $\Ff$ be any finite family of rich
  components of gluing varieties associated to weakly branched ideal
  triangulations of $\hat V$. For every component $Z$ in $\Ff$, every
  point $\rho \in \Omega_\Ff$, and every $f$-weight $(k_f,h_f)$ of
  $(V,\rho)$ such that $i^*(h_f)=0$, there is a point $[w;f]_Z \in
  \Omega_\Ff(Z)_{\infty,0}$ such that $(L([w;f]_Z),\gamma_2'(f)) =
  (k_f,h_f)$.
\end{lem}
Put
$$A_{0,\Ff}(M) := \mathfrak{h}(\Omega_\Ff)$$
and
$$A_{0,\Ff}(M)_\infty := \pi_\infty^{-1}(A_{0,\Ff}(M))\ ,\  
A_{0,\Ff}(M)_N := {\rm EXP} \circ l_{N,k_c} (A_{0,\Ff}(M)_\infty). $$
Replacing $G(T,\tilde b)$ with $\Omega_\Ff(Z) \subset Z \subset
G(T,\tilde b)$, and $G_0(T,\tilde b)_\infty$ with
$\Omega_\Ff(Z)_{\infty,0} \subset Z_{\infty,0} \subset G_0(T,\tilde
b)_\infty$, and taking the restrictions of each of the spaces and maps
that form the sequence of analytic configurations
$\{\Aa_N(T,\tilde b,c)\}$, one obtains a sequence of
sub-configurations
$$\Aa_N(c,Z,\Ff) :=\{\Omega_\Ff(Z)_{\infty,0}, 
\Omega_\Ff(Z)_N,p_\infty,p_N,l _{N,*_b,c},\Hh_1, \Hh_N\}.$$ 
By Lemma \ref{diagram1} it fits into a commutative diagram:
$$\xymatrix{ & & & \Omega_\Ff(Z)_N \ar[rd]^{\Hh_N^{'}} & 
  \\ A_{0,\Ff}(M) & & A_{0,\Ff}(M)_N \ar[ll]_{\pi_N} & &
  \mc\\ & \Omega_\Ff(Z) \ar[lu]^{\mathfrak{h}\circ \rho}
  \ar[ld]_{\mathfrak{h}\circ \rho} & & \Omega_\Ff(Z)_{\infty,0}
  \ar[ll]_{\hspace*{1.5cm} p_\infty}^{\hspace*{-1cm} e^{l_{N,k_c}}}
  \ar[lu]_{loghol_{N,c}} \ar[ld]^{loghol_\infty} \ar[ru]_{\Hh_N}
  \ar[rd]^{\Hh_1}
  \ar[uu]_{e^{l_{N,*_b,c}}}& \\
  A_{0,\Ff}(M) & & A_{0,\Ff}(M)_\infty \ar[uu] |!{[ur];[ul]}\hole
  \ar[ll]_{\pi_\infty} & & \mc}$$ 
The following two results
show that if we fix the bulk weight $h_f$ together with the $c$-weight
$(h_c,k_c)$, then $\Hh_1$ and $\Hh_N$ factor through loghol$_\infty$
and loghol$_{N,c}$, and hence induce maps on $A_{0,\Ff}(M)_\infty$ and
$A_{0,\Ff}(M)_N$, respectively.
\begin{teo}\label{cusped-invariance} 
  Let $\Ff$ be any finite family of rich components $Z$ of gluing
  varieties $G(T_Z,\tilde b_Z)$ associated to weakly branched ideal
  triangulations $(T_Z,\tilde b_Z)$ of $\hat V$. Let $\rho\in \Omega_\Ff$, the $f$-weight
  $(h_f,k_f)$ and $[w;f]_Z \in \Omega_\Ff(Z)_{\infty,0}$ be as in
  Lemma \ref{prelemma}, and $(h_c,k_c)$ be any $c$-weight realized on the triangulation $T_Z$ by a global charge $c_Z$. Then:

  (1) For every odd $N\geq 3$ the scalar $\Hh^\Ff_N(\Pp):= \Hh_N(T_Z,\tilde b_Z,c_Z)([w;f]_Z)$ is an invariant of the pattern $\Pp=(V,\rho,h_f,k_f,h_c,k_c)$ up to multiplication by $2N$-th roots of unity, and depends on $k_f$ mod$(\pi i N)$ only.
  
  (2) The scalar $\Hh^\Ff_1(\Pp):= \Hh_1(T_Z,\tilde b_Z)([w;f]_Z)$ is an invariant of $\Pp=(V,\rho,h_f,k_f)$ up to multiplication by $6$th roots of unity, and is defined exactly if $\tilde b_Z$ is a branching.
  \end{teo}
\Dim We have to show that, up to the phase ambiguities, $\Hh^\Ff_N(\Pp)$ and $\Hh^\Ff_1(\Pp)$ do not depend on the encoding of $\Pp$, ie. the choice of $(T_Z,\tilde b_Z, c_Z)$, $Z$ and $[w;f]_Z$. For that we use the {\it QH transits} of \cite{Top, GT, AGT}, that we adapt to weakly branched triangulations in Section \ref{I-enhanced}, in the case $N\geq 3$; for $N=1$ this adaptation follows easily from Section 4.2-4.3 of \cite{GT} and the argument of Proposition \ref{no-sign-H}. As a result, given two realizations of the pattern $\Pp$,
say: \smallskip

- $[w;f]_Z \in \Omega_\Ff(Z)_{\infty,0}$, 
included into some configuration $\Aa_N(c_Z,Z,\Ff)$, and
\smallskip

- $[w';f']_{Z'} \in \Omega_\Ff(Z')_{\infty,0}$, included 
into another configuration $\Aa_N(c_{Z'},Z',\Ff)$,
\smallskip

\noindent there exists a finite sequence of QH transits relating $[w;f]_Z$ to $[w'';f'']_{Z'}\in \Omega_\Ff(Z')_{\infty,0}$ included into a configuration $\Aa_N(c_{Z'}',Z',\Ff)$  differing from $\Aa_N(c_{Z'},Z',\Ff)$ at most by the global charge, and such that
$$ \Hh_N(T_Z,\tilde b_Z,c_Z)([w;f]_Z) \equiv_{2N} 
\Hh_N(T_{Z'},\tilde b_{Z'},c_{Z'}')([w'';f'']_{Z'}). $$
It follows from the definition of $\Omega_\Ff \subset X_0(M)$ that $w'' = w'$. A further transit argument for flattenings and charges shows finally that 
$$\Hh_N(T_{Z'},\tilde b_{Z'},c_{Z'}')([w';f'']_{Z'})\equiv_{2N} 
\Hh_N(T_{Z'},\tilde b_{Z'},c_{Z'})([w';f']_{Z'}).$$
Since the QH state sums \eqref{SS} depend only on the 
reduction mod$(N)$ of flattenings, the dependence of $\Hh_N^{\Ff}$, $N\geq 3$, on $k_f$ mod$(\pi i N)$ follows immediately.  \cvd 
\medskip

The above constructions depend on the choice of $\Ff$. A canonical choice is the family $\Ff_{\rm EP}$
of rich components $Z$ contained in any weakly branched
EP-triangulation of $\hat V$ (see Proposition \ref{EP-triang}).  By
definition, the analytic configuration $\Aa_N(Y)$ of $Y=(V,(h_c,k_c))$
is the family of analytic sub-configurations $\Aa_N(c,Z,\Ff_{\rm EP})$
for varying charges and rich components $Z\in \Ff_{\rm EP}$.
\begin{cor}\label{ratinv} The invariants $\Hh^{\Ff_{\rm EP}}_1(\Pp)$ 
  and $\Hh^{\Ff_{\rm EP}}_N(\Pp)$ of patterns $\Pp$ over $M$, such that
  $\Pp=(V,\rho,h_f,k_f,h_c,k_c)$ with $\rho \in \Omega_{\Ff_{\rm
      EP}}$, fixed $c$-weight $(h_c,k_c)$, and
  fixed bulk $f$-weight $h_f\in H^1(V;\mz/2\mz)$ such that $i^*(h_f)
  =0$, define respectively an analytic function
$$\Hh^{\Ff_{\rm EP},h_f}_1 :  A_{0,\Ff_{\rm EP}}(M)_\infty \ra \mc$$
and a regular rational function   
$$\Hh_N^{\Ff_{\rm EP},h_f,h_c,k_c} :
A_{0,\Ff_{\rm EP}}(M)_N\ra \ \mc/\mu_{2N}$$ that give rise to a
commutative diagram
$$\xymatrix{ & & & \Omega_{\Ff_{\rm EP}}(Z)_N \ar[rd]^{\Hh_N^{'}} 
& \\ A_{0,\Ff_{\rm EP}}(M) & & A_{0,\Ff_{\rm EP}}(M)_N \ar[ll]_{\pi_N} 
\ar[rr]_{\hspace*{1,5cm} 
e^{l_{N,*_b,c}}}^{\hspace*{-1,4cm} \Hh^{\Ff_{\rm EP},h_f,h_c,k_c}_N} |
!{[ur];[dr]}\hole & & \mc/\mu_{2N}\\ &
\Omega_{\Ff_{\rm EP}}(Z) \ar[lu]^{\mathfrak{h}\circ \rho_Z} 
\ar[ld]_{\mathfrak{h}\circ \rho_Z} & & 
\Omega_{\Ff_{\rm EP}}(Z)_{\infty,0} 
\ar[ll]_{\hspace*{1.5cm} p_\infty}^{\hspace*{-1cm} e^{l_{N,k_c}}} 
\ar[lu]_{loghol_{N,c}} \ar[ld]^{loghol_\infty} 
\ar[ru]_{\Hh_N} \ar[rd]^{\Hh_1} \ar[uu] & \\
A_{0,\Ff_{\rm EP}}(M) & & A_{0,\Ff_{\rm EP}}(M)_\infty 
\ar[uu] |!{[ur];[ul]}\hole \ar[ll]_{\pi_\infty} 
\ar[rr]^{\hspace*{1cm}\Hh^{\Ff_{\rm EP},h_f}_1} & & \mc}$$

\end{cor}
\noindent {\it Proof.}  Any point $x$ in the image of loghol$_\infty$
or loghol$_{N,c}$ determines a unique point $\rho =
\mathfrak{h}^{-1}(\pi_\infty(x))$ or $\mathfrak{h}^{-1}(\pi_N(x))$ in
$\Omega_{\Ff_{\rm EP}}$. If $x\in {\rm Image}($loghol$_\infty$), by
Definition \ref{defiweight} it determines also a boundary $f$-weight
$k_f$, and by Lemma \ref{prelemma} any boundary $f$-weight $k_f$ for
$\rho$ such that $(k_f-d_\rho)/i\pi$ is the zero class in
$H^1(\partial V;\mz/2\mz)$ is realized by a point of
$\pi_\infty^{-1}(\mathfrak{h}(\rho))$ in the image of loghol$_\infty$
(see \eqref{fconstraint2}). Hence the tuples $(x,h_f)$ and
$(\rho,h_f,k_f)$ with $i^*(h_f)=0$ are in one-to-one
correspondence. If $x\in {\rm Image}($loghol$_{N,c}$) and a boundary
$c$-weight $k_c$ is given, $x$ determines the reduction mod$(\pi i N)$ of
$k_f$. Hence the tuples $(x,h_f,h_c,k_c)$ and $(\rho,h_f,k_f \ {\rm
  mod}(\pi i N),h_c,k_c)$ with $i^*(h_f)=0$ are in one-to-one
correspondence. From this and Theorem \ref{cusped-invariance} it
follows that $\Hh_1$ and $\Hh_N$ descend to maps $\Hh^{\Ff_{\rm
    EP},h_f}_1$ and $\Hh_N^{\Ff_{\rm EP},h_f,h_c,k_c}$ defined on
the image of loghol$_\infty$ and loghol$_{N,c}$ respectively. Since
$\Hh_1$ is analytic and $\pi_\infty$ and $p_\infty$ are covering maps,
$\Hh^{\Ff_{\rm EP},h_f}_1$ is analytic too and it extends uniquely to
an analytic function on $A_{0,\Ff_{\rm EP}}(M)_\infty$. Similarly $\Hh_N^{\Ff_{\rm EP},h_f,h_c,k_c}$ can be extended in
an unique way to a regular rational function on $A_{0,\Ff_{\rm
    EP}}(M)_N$. \cvd 
        
\begin{remark}\label{canonical}{\rm The Epstein-Penner family $\Ff_{\rm EP}$ of rich components is only one of the possible choices.  For instance, for every $M$ and $m\geq 2$ we can consider the family $\Ff_m$ of rich components of gluing varieties of weakly branched triangulations of $\hat V$ with no more than $m$ tetrahedra, or the family $\Ff_{m_0}$ where $m_0$ is the minimum $m$ such that $\Ff_m \neq \emptyset$. For example if $M$ is the complement of the knot $5_2$ in $S^3$, then $m_0=3$, while the Epstein-Penner subdivision is a triangulation made by $4$ hyperbolic ideal tetrahedra. In practice, any rich component is suited to study the asymptotical behaviour of the QHI at small deformations of the point $z_{hyp}$ representing the complete hyperbolic structure.}  
\end{remark}
    
\noindent {\bf Proof of Theorem \ref{more_cusped} (1)-(2).} 
For the first two statements, apply Theorem \ref{cusped-invariance}
to the family $\Ff = \Ff_{\rm EP}$, set $\Omega(M):=
\Omega_{\Ff_{\rm EP}}$, and use Corollary
\ref{ratinv} to pull back the covering $\pi_N$ via $\mathfrak{h}$.
As explained in Remark \ref{extclass}, the assumption $i^*(h_f) =0$ is not necessary; in general one gets functions $A_{0,\Ff_{\rm
    EP}}(M)_N^{\varepsilon,\varepsilon'}\ra \ \mc/\mu_{2N}$ where $\varepsilon,\varepsilon'\in \{\pm 1\}$. \cvd \medskip
        
%The organization here is different with respect to LIM_9 : the above remark was at the end of Section 4.3, and all mentions of the SYM factor and the reduced invariants is dropped now.  Some adaptations are done at the end to make more explicit the connection with Proposition \ref{top-support2} (2)

\noindent {\bf Proof of Proposition \ref{top-support} (2) and \ref{top-support2} (2).}   Let $(T,\tilde
b)$ be a weakly branched ideal triangulation of a pattern $(V,\rho,(h,k))$ over $M$. Denote by $(\Delta_j,b_j)$, $j\in \{1,\ldots,s\}$, the $3$-simplices of $(T,\tilde b)$. We need to review the Neumann-Zagier theory (see \cite{N0,N}), that we will apply to the weakly branched triangulation $(T,\tilde b)$. Denote :
\begin{itemize}
\item by $C_1(\Delta_j)$ the $\mz$-module freely generated by the edges $e_k^j$ and $(e_k')^j$, $k\in \{0,1,2\}$, of $\Delta_j$;
\item by $\bar{J}_j$ the quotient module of $C_1(\Delta_j)$ by the relations $e_k^j =  (e_k')^j$, $k\in \{0,1,2\}$;
\item by $J_j$ the quotient module of $\bar{J}_j$ by the relation $e_0^j+e_1^j+e_2^j = 0$.
\end{itemize}
Let $C_1$ be the $\mz$-module freely generated by the edges of $T$, and put $E(T) := \oplus_j C_1(\Delta_j)$ (the $\mz$-module of ``abstract'' edges of $T$), $\bar{J} := \oplus_j \bar{J}_j$, and $J := \oplus_j J_j$. Consider on theses spaces the following $2$-forms: the standard inner product $(\ ,\ )$ on $C_1$, $E(T)$ and $\bar{J}$ (each defined with respect to the natural basis given by (cosets of) edges), and the {\it signed} antisymmetric bilinear form $\langle \ ,\ \rangle$ on $J$ given by
\begin{equation}\label{pairing0}
\langle \ ,\ \rangle = \oplus_j  *_j \langle \ ,\ \rangle_j
\end{equation}
where $\langle \ ,\ \rangle_j$ is the standard antisymmetric bilinear form on $J_j$, defined in the basis $\{e_0^j,e_1^j\}_j$ by 
\begin{equation}\label{defform}
\langle e_0^j,e_1^j \rangle_j = 1\quad ,\quad \langle e_0^j,e_1^m \rangle =0\  {\rm for}\  m\ne j.
\end{equation}
Hence $\{e_0^j,e_1^j\}$ is a symplectic basis of $\langle \ ,\ \rangle$ restricted to $J_j$ when $\Delta_j$ has positive branching orientation $*_j=+1$, and is an anti-symplectic basis of $J_j$ otherwise. Equivalently, $\{e_0^j,e_1^j\}$ is  symplectic if the edges $e_0^j$, $(e_1')^j$, $e_2^j$ have positive cyclic ordering as viewed from their common endpoint $v_0$ in $\Delta_j$ (in \cite{N0} the branchings are not needed and not used; one fixes positively cyclically ordered edges $e_0^j$, $(e_1')^j$, $e_2^j$ independently for all $\Delta_j$, so that the signs $*_j$ occuring here are systematically replaced by $+1$). Define linear maps $\beta : C_1 \ra E(T)$ and $\beta^* : E(T) \ra C_1$ by 
$$\beta(e) = \sum_{E\ra e} E$$
and 
$$\beta^*(e_k^j) = *_j\left( p(e_{k+1}^j)-p(e_{k+2}^j)+p((e_{k+1}')^j)-p((e_{k+2}')^j)\right) $$ 
and the same for $\beta^*((e_k')^j)$, where indices are regarded mod$(3)$, and $p : E(T) \ra C_1$ is the identification map, assigning to an abstract edge its coset in $T$. They induce maps (we keep the same notations)
$$\beta : C_1 \ra J\quad ,\quad \beta^* : J \ra C_1.$$
Clearly, for any edge $e$ of $T$ we have
\begin{equation}\label{signedad}
(\beta^*(e^j_k),e) = \langle e^j_k,\beta(e)\rangle.
\end{equation}
Also, $\beta^*$ splits as $ \beta^* : J \stackrel{\beta_1}{\lra} \bar{J} \stackrel{\beta_2}{\lra} C_1$, where $\beta_1(e_k^j) := *_j(e_{k+1}^j - e_{k+2}^j)$ and $\beta_2(e_k^j) := p(e_k^j) + p((e_k')^j)$. Similarly to \eqref{signedad}, $\beta_1$ and the natural projection $q:\bar{J} \ra J$ are adjoint maps: for all $x\in \bar{J}$ and $a\in J$ we have
\begin{equation}\label{adjointq}
(\beta_1(a),x) = \langle a,q(x)\rangle.
\end{equation}
Denote by $C_0$ the $\mz$-module freely generated by the vertices of $T$, and by $\alpha^*:C_1\rightarrow C_0$ the map which assigns to an edge  the sum of its vertices. One has ${\rm Im}(\beta^*) \subset {\rm Ker}(\alpha^*)$, and interpreting the elements of ${\rm Ker}(\alpha^*)$ as unoriented $1$-cycles in $T$, there is an isomorphism
\begin{equation}\label{1exactseq}
{\rm Ker}(\alpha^*)/{\rm Im}(\beta^*) \cong H_1(T;\mz/2\mz).
\end{equation}
It is not hard to see that ${\rm Im}(\beta) \subset {\rm Ker}(\beta^*)$. Put $H := {\rm Ker}(\beta^*)/{\rm Im}(\beta)$. The $2$-form $\langle \ ,\ \rangle$ on $J$ descends to $H$. Since ${\rm Ker}(\beta^*) = {\rm Im}(\beta)^{\perp}$, $\langle \ ,\ \rangle$ is non degenerate on $H/{\rm Tors}(H)$. In fact, it is shown in \cite{N0} that there is an exact sequence
\begin{equation}\label{themainseq}
\xymatrix{0 \ar[r] & H \ar[rr]^{\hspace*{-1.8cm}(\gamma',\gamma'_2)} & & H^1(\partial V;\mz) \oplus H^1(V;\mz/2\mz) \ar[r]^{\hspace*{0,9cm} r - i^*}  & H^1(\partial V;\mz/2\mz) \ar[r] & 0.} \end{equation}
%Hence $H$ has only $2$-torsion elements, and taking coefficients in $\mc/2\pi i \mz$, the exact sequence shows that $\gamma'$ yields an isomorphism
%$$\gamma' :H \otimes  (\mc/2\pi i \mz) \lra H^1(\partial V; \mc/2\pi i \mz).$$
Moreover, denoting by $\cdot$ the intersection product on $H_1(\partial V;\mz)$, and by $\gamma = PD \circ \gamma'$ the map $\gamma'$ followed by the Poincar\'e duality isomorphism, for all $a,b\in H$ we have
\begin{equation}\label{intersec}
\gamma(a) \cdot \gamma(b) = 2\langle a,b\rangle.
\end{equation}
The map $\gamma$ is defined as follows. Recall the notations introduced before Definition \ref{defiweight}. For any $2$-face $F$ of $\partial T_0$ which is a boundary section of the truncated tetrahedron of $T_0$ corresponding to $\Delta_j$, let us write $F\ra \Delta_j$, and denote by $s_k^{j,F}$ the 
edge of $F$ which is opposite to the vertex of $F$ that belongs to the edge $e_k^j$ or $(e_k')^j$, with the positive orientation as viewed from that vertex. Then, the linear map 
\begin{equation}\label{mapgamma}
\fonc{\bar{\gamma}}{\bar{J}}{C_1(\partial T_0)}{e_k^j}{\sum_{F\ra \Delta_j}s_k^{j,F}.}
\end{equation}
descends to 
\begin{equation}\label{gammadefmap}\gamma:H \ra H_1(\partial V;\mz).
\end{equation}
Conversely, let $q:\bar{J}\ra J$ be the natural projection. Denote by $\partial T_0'$ the cellulation of $\partial V$ dual to $\partial T_0$. Represent classes in $H_1(\partial V; \mz)$ by simplicial loops in $\partial T_0'$. For all $F\ra \Delta_j$, denote by $a_k^{j,F}$ the  simplicial arc in $\partial T_0' \cap F$ which faces the vertex of $F$ that belongs to the edge $e_k^j$ or $(e_k')^j$, with the positive orientation as viewed from that vertex. Define a linear map
$$\fonc{\bar{\delta}}{C_1(\partial T_0')}{\bar{J}}{a_k^{j,F}}{e_k^j.}$$
Then, $q\bar{\delta}$ descends to a map $\delta : H_1(\partial V; \mz) \ra H$. Moreover, for all $x\in H_1(\partial V; \mz)$ we have
\begin{equation}\label{twice}
\gamma \circ \delta (x) = 2x.
\end{equation}
The identity \eqref{intersec} follows from \eqref{twice}, together with the simple fact that for any simplicial loop $C$ in $\partial T_0'$, we have
\begin{equation}\label{intersecbis}
C \cdot \bar{\gamma}(e_k^j)  = *_j\ \langle \bar{\delta}(C),e_k^j\rangle_j\ ,
\end{equation}
where the sign $*_j$ comes from the orientation induced by the basis $\{ e_0^j, e_1^j\}$ of $J_j$, as discussed after \eqref{defform}. Hence $\gamma$ and $\delta$ are adjoint maps with respect to the form $\langle \ ,\ \rangle$ on $J$ and the intersection form on $H_1(\partial V;\mz)$: for all $x\in H_1(\partial V;\mz)$ and $a\in H$ we have
\begin{equation}\label{intersec3}
\gamma(a) \cdot x  = \langle a,\delta(x)\rangle.
\end{equation}
Taking coefficients in $\mz/2\mz$ and simplicial loops in $T_0'$ in normal position with respect to $T_0$, one defines similarly to $\delta$ a map $\delta_2:H_1(V;\mz/2\mz)\ra H\otimes (\mz/2\mz)$. The map $\gamma_2'$ in the above exact sequence is defined for all $a\in H\otimes (\mz/2\mz)$ and $c\in H_1(V;\mz/2\mz)$ by   
$$\gamma'_2(a)(c) = \langle a , \delta_2(c)\rangle .$$
With this material at hands, on can show that the maps $\gamma'$ and $\gamma'_2$ produce the boundary and the bulk weights of $V$. This is done for the charge weights in Section 6 of \cite{N0}, and for the flattening weights in Section 9 of \cite{N}, with the further constraint that the weights are taken equal to $0$. For the sake of completeness, let us explain the case of the flattening weights, the adaptation to the charge weights being not difficult. 

Consider the weakly branched ideal triangulation $(T,\tilde
b)$ of $(V,\rho,(h,k))$. By assumption, there is a point $w\in G(T,\tilde b)$ having $\rho$ as holonomy.  Take for each $3$-simplex $(\Delta_j,b_j)$ the (local) flattening of the form $(\underline{l}_0^j,\underline{l}_1^j,\underline{l}_2^j) := (\log(w^j_0),\log(w^j_1),\log(w^j_2) +\varepsilon i\pi)$, $\varepsilon \in \{-1,+1\}$, and set $\underline{l}:=(\underline{l}_0^j,\underline{l}_1^j,\underline{l}_2^j)_j$. Consider the vector
\begin{equation}\label{vl}
v_{\underline{l}} := \sum_{j=1}^s  \underline{l}_{1}^j e_0^j - \underline{l}_{0}^j e_1^j \quad \in J\otimes \mc.
\end{equation}
We have
$$\beta_1(v_{\underline{l}}) =  \sum_{j=1}^s *_j( \underline{l}_{0}^j e_0^j + \underline{l}_{1}^j e_1^j - (\underline{l}_{0}^j+\underline{l}_{1}^j)e_2^j) = \sum_{j=1}^s *_j( \underline{l}_{0}^j e_0^j + \underline{l}_{1}^j e_1^j + \underline{l}_{2}^je_2^j) \quad \in \bar{J}\otimes \mc.$$
So $\beta_1(v_{\underline{l}})$ represents the map assigning the value $*_j\underline{l}_k^j$ to the abstract edges $e_k^j$ and $(e_k')^j$ of $T$. Recall that for every $\bar {a}\in \bar{J}$, $\beta_2(\bar {a})$ computes at every edge $e\in C_1$ the sum of the coefficients of $\bar {a}$ attached to the abstract edges having coset $e$. Then, $\beta^* = \beta_2\circ \beta_1$, the gluing equation \eqref{gleq}, and the argument of Lemma \ref{even} imply that $\beta^*(v_{\underline{l}})\in 2\pi i C_1$ (note that in the Lemma the label $\emptyset$ is attached to edges of the type $e_2^j$ and $(e_2')^j$). By \eqref{1exactseq} we deduce that $\beta^*(v_{\underline{l}})/\pi i$ represents the $0$ class in $H_1(T;\mz/2\mz)$, and so there exists $\underline{f} \in J$ such that $\beta^*(\underline{f}) = \beta^*(v_{\underline{l}})/\pi i$. Put $v_l := v_{\underline{l}}-\pi i \underline{f}$, and write its coefficients as 
$$\textstyle v_l = \sum_{j=1}^s  l_{1}^j e_0^j - l_{0}^j e_1^j = \sum_{j=1}^s  (\log(w^j_1) + \pi i f^j_1) e_0^j - (\log(w^j_0) + \pi i f^j_0) e_1^j.$$
Since $\beta^*(v_l)=0$, the collection $l:=(l_0^j,l_1^j,l_2^j)_j$ satisfies the edge relations $L(e)=1$ (see \eqref{edgeconst}). The tetrahedral relations are satisfied because any element of ${\rm Im}(\beta_1)$ has coefficient sum equal to $0$ in each module $\bar{J}_j$. So $l$ defines a flattening of $w$. Consider the class $[v_l]\in H\otimes \mc$. One can see that $\gamma'([v_l]) = L([w;f])$ as follows. Consider the map 
$$\fonc{\bar{\gamma}''}{\bar{J}\otimes \mc}{C^1(\partial T_0';\mc)}{x}{(C \mapsto (x,\bar{\delta}(C))).}$$ 
Denote by $A\in J\otimes \mc$ a representative of $a\in H\otimes \mc$. Then $\bar{\gamma}''$ induces a map 
$$\fonc{\gamma''}{H\otimes \mc}{H^1(\partial V; \mc)}{a}{\bar{\gamma}''(\beta_1(A)).}$$
Using successively \eqref{adjointq} and \eqref{intersec3}, we get
\begin{equation} \label{dualP} \gamma''(a)([C]) = (\beta_1(A), \bar{\delta}(C)) = \langle a,q\bar{\delta}([C])\rangle = \gamma(a)\cdot [C].\end{equation}
Hence $\gamma'' = \gamma'$. Clearly, $\gamma''([v_l]) = \bar{\gamma}''(\beta_1(v_l))$ is given by \eqref{weightwf}.  This proves $\gamma'([v_l]) = L([w;f])$. The collection $f:=(f_0^j,f_1^j,f_2^j)_j$ satisfies also the edge relations mod$(2)$, hence it defines a class $\gamma_2'(f)\in H^1(V;\mz/2\mz)$  as explained after \eqref{weightwf}. By construction we have the compatibility relations \eqref{fconstraint2}, so the pair $(L([w;f]),\gamma_2'(f))$ is a flattening weight of $(V,\rho)$. Finally, take any other flattening weight $(k_f',h_f')$ of $(V,\rho)$. Comparing \eqref{fconstraint} and \eqref{fconstraint2}  we see that $k_f'-L([w;f]) \in H^1(V;\pi i\mz)$ and $(k_f'-L([w;f]))/\pi i  = i^*(h_f'-\gamma_2'(f))$ mod$(2)$. By the exact sequence  \eqref{themainseq}, there exists $a\in H$ such that $\gamma'(a) = (k_f'-L([w;f]))/\pi i$ and $\gamma'_2(a) = h_f'-\gamma_2'(f)$. Then the vector $\beta_1(v_l +\pi i a)\in \bar{J}\otimes \mc$ represents a flattening of $w$ with weight $(h_f',k_f')$.\hfill$\Box$

\subsection{Relation with Chern-Simons theory}\label{qCS}
Recall that $V$ is the compact $3$-manifold obtained by
completing the cusp of $M$ with a torus. Chern-Simons theory with
gauge group $PSL(2,\mc)$ associates to $\partial V$ the {\it Chern-Simons
  bundle} $\Ll_{\partial V} \ra X(\partial V)$, which is a $\mc^*$-bundle with
canonical connection $1$-form and canonical inner product, and to $V$ the parallel {\it Chern-Simons section} $s_V : X(V) \ra i^*\Ll_{\partial V}$, where $i:X(V) \rightarrow X(\partial V)$ is the restriction map (\cite{F,KK}). 

We are going to show that when the bulk $f$-weight $h_f$ is $0$ the bottom row of the diagram of Theorem \ref{ratinv} encodes the restriction of the pair $(i^*\Ll_{\partial V},s_V)$ to $\Omega_{\Ff_{\rm EP}} \subset X_0(M)$, identified as a component of $X(V)$. For arbitrary $h_f$ one gets a twisted version of $(i^*\Ll_{\partial V},s_V)$. The proof relies on a certain number of results that we need to recall. Denote by CS$(M') \in \mr/\mz$ the Chern-Simons invariant of the Levi-Civita connection of a closed riemannian $3$-manifold $M'$, and by CS$(M)$ Meyerhoff's extension to cusped hyperbolic manifolds $M$. By definition, CS$(M)$ is the limit of CS$(M'_n)$ for any sequence $(M'_n)$ of closed hyperbolic
Dehn fillings of $M$ converging to $M$ in Dehn surgery space. In \cite{N}, Corollary 14.6 $\&$ Theorem 14.7, Neumann proved (using the notations of the present paper): 
\begin{teo} \label{geomH0} Let $Y$ be $M$ or a closed hyperbolic Dehn
  filling $M'$ of $M$. Denote by $\rho_Y$ the hyperbolic holonomy of $Y$, identified as a character of $V \subset M'$ in the case $Y=M'$. If $Y=M$, put $k_Y=0$, and if $Y=M'$, let $k_Y$ be any boundary $f$-weight of $V$ relative to $\rho_Y$ such that $k_Y$ vanishes on the meridian of the added solid torus. Then
 $$\Hh_1(V,\rho_Y,(0,k_Y))= \exp\left(\frac{2}{\pi}{\rm Vol}(Y)+
2\pi i{\rm CS}(Y)\right).$$
\end{teo}
The expression on the right hand side is related to the Chern-Simons section $s_V$ as follows. To any {\it compact closed} oriented $3$-manifold $Y$, $PSL(2,\mc)$-Chern-Simons theory associates a {\it function} $s_Y : X'(Y) \ra \mc^*$ defined on the variety $X'(Y)$ of $PSL(2,\mc)$-characters of $Y$. If $Y$ is hyperbolic with holonomy $\rho_{Y}$, a classical result of Yoshida \cite{Yo} gives 
\begin{equation}\label{firstexp0}
s_Y(\rho_{Y}) = \exp\left(\frac{2}{\pi}{\rm Vol}(Y)+2\pi i{\rm CS}(Y)\right).
\end{equation}
Assume that $Y = V \cup_{\varphi} (D^2\times S^1)$ is a closed hyperbolic Dehn filling of $M$ whose holonomy $\rho_Y$ lies in a sufficiently small neighborhood $D$ of $\rho_{hyp}$ in $X(V)=X(M)$ (here we identify $\rho(Y)$ with an augmented character of $V$). Then Kirk-Klassen (\cite{KK}, pages 554--556) showed that
\begin{equation}\label{firstexp}
s_Y(\rho_Y) = \langle s_V(\rho_Y),s_{D^2\times S^1}(\rho_Y)\rangle
\end{equation}
where $\langle \ ,\ \rangle$ is the canonical inner product of $i^*\Ll_{\partial V}$, and on the right hand side $\rho_Y$ denotes also the augmented character induced on the glued solid torus $D^2\times S^1$. Explicitely, if we fix the gauge on $i^*\Ll_{\partial V}$ by taking as coordinates on $D$ the ones corresponding under the map $\mathfrak{h}$ to the {\it standard} logarithms $\log(\lambda)$ and $\log(\mu)$, which are equal to $0$ at $\rho_{hyp}$, we have
\begin{equation}\label{firstexp2}
s_V(\rho_{hyp})  = \exp \left(\frac{2}{\pi}{\rm Vol}(M)+ 2\pi i{\rm CS}(M)\right).
\end{equation}
Then, if the Dehn filling instruction $\varphi$ maps the meridian $\partial(D^2\times *)$ to $l^qm^p$ and the longitude $*\times S^1$ to $l^sm^r$ (so $ps-qr=1$), the formula \eqref{firstexp} splits as
\begin{multline}\label{firstexp1}
s_Y(\rho_Y) = s_V(\rho_{hyp})\\ \times \exp \left(- \frac{1}{2\pi i} \int_{\rho_{hyp}}^{\rho_Y} \left( \log(\lambda) d\log(\mu) - \log(\mu)d\log(\lambda)\right)\right) \exp\left(-(s\log(\lambda)+r\log(\mu))\right).
\end{multline}
The first exponential gives the variation of $s_V$ between $\rho_{hyp}$ and $\rho_Y$ ($s_V$ being a parallel section of $i^*\Ll_{\partial V}$, this exponential is the holonomy of the connection $1$-form of $i^*\Ll_{\partial V}$ between the two points); the product with $s_V(\rho_{hyp})$ is thus the value of $s_V(\rho_Y)$ at the chosen gauge. The second exponential is the value of $s_{D^2\times S^1}(\rho_Y)$ in the gauge corresponding to the Dehn filling: $\log({\rm hol}_{\partial(D^2\times *)}) = q\log(\lambda)+p\log(\mu) = 2\pi i$, $\log({\rm hol}_{*\times S^1}) = s\log(\lambda)+r\log(\mu)$. Its argument is $-$(length $+ i$ rotation angle) of the geodesic core of the surgery torus. The arguments of both exponential can be collected in a single integral by changing the coordinates $\log(\lambda)$, $\log(\mu)$ on $D$ to $\log(\lambda)+2\pi i r$, $\log(\mu)-2\pi i s$, that is, by working on a different leaf of the Riemann surface of these maps. Then the gauge $\log({\rm hol}_{\partial(D^2\times *)})$ vanishes, like $k_Y$ in Theorem \ref{geomH0}.  

We are going to show that $(i^*\Ll_{\partial V},s_V)_{\vert \Omega_{\Ff_{\rm EP}}}$ is determined by $\Hh^{\Ff_{\rm EP},0}_1: A_{0,\Ff_{\rm EP}}(M)_\infty \ra \mc$. Because of Theorem \ref{geomH0} and \eqref{firstexp2}, it is enough to identify the variation of $\Hh^{\Ff_{\rm EP},0}_1$ with the one of $s_V$, given by \eqref{firstexp1}. This is the content of the following result.
\begin{prop}  The function $\Hh^{\Ff_{\rm EP},0}_1$ descends to a parallel section $\Ss^0$ of a flat trivial $\mc^*$-bundle $\mathcal{L}(M)\ra A_{0,\Ff_{\rm EP}}(M)$ with canonical connection $1$-form and inner product, such that $\mathfrak{h}^*(\mathcal{L}(M),\Ss^0)$ is isomorphic to $(i^*\Ll_{\partial V},s_V)_{\vert \Omega_{\Ff_{\rm EP}}}$.
\end{prop} 

\noindent {\it Proof.} We define $\mathcal{L}(M)$ from the variation of the invariant $\Hh^{\Ff_{\rm EP}}_1$ of Theorem \ref{cusped-invariance}, like the Chern-Simons line bundle is defined from the Chern-Simons action in \cite{F,KK}. Hence $\mathcal{L}(M)$ is the quotient of $A_{0,\Ff_{\rm EP}}(M)_\infty
\times \mc^*$ under the action of $\mz^2$ given by
\begin{equation}\label{bundledef}
\forall a,b\in \mz\ ,\  (a,b)\cdot ([\lambda;p],[\mu;q],z) 
:= ([\lambda;p+a],[\mu;q+b],ze^{b \log(\lambda)-a \log(\mu)}).
\end{equation}
We have a $\mc^*$-bundle projection
$$\fonc{\pi_\infty'}{\mathcal{L}(M)}{A_{0,\Ff_{\rm EP}}(M)}
{[[\lambda;p],[\mu;q],z]}{(\lambda,\mu).} $$ Denote by
$\mathcal{L}(M)^{-1}$ the inverse $\mc^*$-bundle, which is the same as $\mathcal{L}_{\Ff_{\rm EP}}(-M)$, where $-M$ denotes $M$ with reversed orientation. An inner
product is defined by the bundle map
$$\fonc{\langle \ ,\ \rangle}{\mathcal{L}(M) 
  \times \mathcal{L}(M)^{-1}}{\mc^*}{(([\lambda;p],
  [\mu;q],z_1),([\lambda;p],[\mu;q],z_2))}{z_1z_2}$$ where $\mc^*$ is regarded as the bundle over a point.  Also, a flat analytic connection $1$-form is defined by the restriction to $A_{0,\Ff_{\rm EP}}(M)_\infty$ of the $1$-form on $\tilde{\mc}\times \tilde{\mc}$ given by (we put $l_a(\lambda) := \log(\lambda)+2\pi i
a$)
\begin{equation}\label{conn}
\eta := -\frac{1}{2\pi i } \left( l_{p}(\lambda) dl_{q}(\mu) 
- l_{q}(\mu)dl_{p}(\lambda)\right).
\end{equation}
By Lemma \ref{prelemma}, given any even class $\theta \in H^1(\partial
V;2\mz)$ and any point $[w;f]_Z\in \Omega_{\Ff_{\rm
    EP}}(Z)_{\infty,0}$, the flattening $f$ can be modified to a
flattening $f'$ such that $[w;f']_Z\in \Omega_{\Ff_{\rm
    EP}}(Z)_{\infty,0}$ and the corresponding boundary weights satisfy
$$L([w;f']_Z) - L([w;f]_Z) = \pi i \theta.$$ 
As in the proof of Proposition \ref{top-support} (2) and \ref{top-support2} (2), denote by $v_l$ (resp. $v_{l'}$) the vector in $H\otimes \mc$ associated to the classical log-branch at $[w;f]_Z$ (resp. $[w;f']_Z$), and similarly denote by $v_d\in H\otimes (\mc/2\pi i \mz)$ the vector associated to $(\log(w_0^j), \log(w_1^j),\log(w_2^j))_j$. Recall the class $d_w\in H^1(\partial V;\mc/2\pi i \mz)$ in \eqref{fconstraint2}. Then $d_w = \gamma'(v_d)$ and
$$\pi i \theta =  \gamma'(v_{l'} - v_{l}) = \gamma'\left(\pi i\sum_{j=1}^s  ((f')^j_1-f^j_1)e_0^j -  ((f')^j_0-f^j_0)e_1^j\right).$$
Hence, denoting $\Hh_1^{\Ff_{\rm EP}} := \Hh_1^{\Ff_{\rm EP}}(T_Z,\tilde b_Z)$ we have
$$\Hh_1^{\Ff_{\rm EP}}([w;f']_Z) \Hh_1^{\Ff_{\rm EP}}([w;f]_Z)^{-1} =\hspace*{7cm}$$
\begin{align}
 \hspace*{2.5cm} &=  
\exp \left(\sum_{j=1}^s *_j \left( ((f')^j_1-f^j_1)\log(w^j_0) - 
((f')^j_0-f^j_0)\log(w^j_1)\right) \right) \notag \\
& = \exp \left(-\frac{1}{\pi i} \langle v_{l'} - v_{l} , v_d\rangle \right)\notag \\ & = \exp \left(-\frac{1}{2\pi i} \langle\! \langle \gamma'(v_{l'} - v_{l}), \gamma'(v_d)\rangle\! \rangle \right)\notag \\
& = \exp \left(\frac{1}{2}\left( \theta(m)\log({\rm hol}_l(\rho(w))) - \theta(l)\log({\rm hol}_m(\rho(w)))\right)\right)
\label{variation}
\end{align}
where the map hol$_C$, $C$ a curve in $\partial V$, is
defined in \eqref{loghmap1}, and $\langle\! \langle \ , \ \rangle\! \rangle$ denotes the cup product followed by the evaluation on the fondamental class. The first equality follows from
\eqref{formR} and \eqref{SS}, the second is by the definition \eqref{pairing0} of the symplectic pairing on $H\otimes \mc$, the third follows from \eqref{intersec} and Poincar\'e duality, and the last one uses that the symplectic basis of $H_1(\partial V;\mz)$ is $(l,m)$. By
setting $\theta(m):=2b$ and $\theta(l):=2a$, we see from
\eqref{bundledef} and \eqref{variation} that the
function $\Hh^{\Ff_{\rm EP},h_f}_1$ (with arbitrary bulk weight $h_f$ for the moment)
descends from $A_{0,\Ff_{\rm EP}}(M)_\infty$ to a section of $\mathcal{L}(M)$. Denote it
$$\Ss^{h_f}: A_{0,\Ff_{\rm EP}}(M) \ra \mathcal{L}(M). $$
We show that this section is parallel with respect to the connection $\eta$ by a computation similar to \eqref{variation}. As remarked in \cite{N0}, Lemma 10.2, for any $\mathbb{Q}$-vector space $E$, the skew symmetric bilinear form $\langle\ ,\ \rangle$ on $J$ induces a symmetric bilinear map
$$\fonc{B}{(H\otimes E)\otimes (H\otimes E)}{E\wedge E}{(a\otimes v)\otimes (b\otimes w)}{\langle a,b\rangle \ v \wedge w.}$$
Consider the map $\cdot_E$ on $H_1(\partial V;E)  = H_1(\partial V;\mz) \otimes E$ induced by the intersection product on $H_1(\partial V;\mz)$:
$$\fonc{\cdot_E}{H_1(\partial V;E)\otimes H_1(\partial V;E)}{E\wedge E}{(x\otimes v)\otimes (y\otimes w)}{(x\cdot y) \ v \wedge w.}$$
Denote again by $\gamma$ the map $\gamma \otimes {\rm id}: H \otimes E \ra H_1(\partial V;E)$. By \eqref{intersec} we have a commutative diagram
$$\xymatrix{(H\otimes E)\otimes (H\otimes E) \ar[r]^{\hspace*{1cm} B} \ar[d]_{\gamma \otimes \gamma}& E\wedge E \ar[d]^{2\times}\\ H_1(\partial V;E) \otimes H_1(\partial V;E) \ar[r]^{\hspace*{1.3cm} \cdot_E} & E\wedge E}$$
In particular, let $E = C^\infty(\mathcal{W}_\infty)$ and denote by $l_0^j$, $l_1^j\in E$ the classical log-branch functions at the edges $e_0^j$, $e_1^j$ of $\Delta_j$. Put $\textstyle a = \sum_{j=1}^s e_0^j \otimes l_1^j - e_1^j \otimes l_0^j \in H\otimes E$. For any point $[w;f]_Z\in \Omega_{\Ff_{\rm EP}}(Z)_{\infty,0}$ we have $\gamma(a)([w;f]_Z) = l\otimes L([w;f]_Z)(m) - m \otimes L([w;f]_Z)(l)$; dually the class $\gamma'(a)\in H^1(\partial V;E)$,  $\gamma'(a) = PD \circ \gamma(a)$, is given by $\gamma'(a)([w;f]_Z) = L([w;f]_Z)\in H^1(\partial V;\mc)$. Moreover, in $E\wedge E$ we have
$$ 2 \ \sum_{j=1}^s  *_j \ (l_0^j \wedge l_1^j) =  B( a,a) = \frac{1}{2}\  \gamma(a) \cdot_E \gamma(a)  = L(.)(l) \wedge L(.)(m).$$
Applying the homomorphism $E \wedge E \rightarrow \Omega^1(\mathcal{W}_\infty)$, $f\wedge g \mapsto fdg -gdf$, we deduce
$$2\ \sum_{j=1}^s  *_j \ (l_0^jdl_1^j - l_1^jdl_0^j) =  L(.)(l) dL(.)(m)  - L(.)(m) dL(.)(l).$$
Hence, for any two points $[w;f]_Z$, $[w';f']_Z$ connected by a continuous path $(\gamma^1,\ldots,\gamma^s)$ in $\Omega_{\Ff_{\rm EP}}(Z)_{\infty,0}$, we have
\begin{equation}\label{varconn2}
\Hh_1^{\Ff_{\rm EP}}([w';f']_Z) \Hh_1^{\Ff_{\rm EP}}([w;f]_Z)^{-1} = \exp  \left(2 \sum_{j=1}^s *_j \int_{\gamma^j}  \eta\right) = \exp\left(\int_{\gamma} \eta\right)
\end{equation}
where $\gamma \subset A_{0,\mathcal{F}_{\rm EP}}(M)_\infty$ is the
image of $(\gamma^1,\ldots,\gamma^s)$ under the map
loghol$_\infty$.  Hence $\Ss^{h_f}$ is a parallel section of $\mathcal{L}(M)$. Multiplying coordinates by $4\pi i$ identifies $\mathfrak{h}^*\mathcal{L}(M)$ with the restriction of $i^*\Ll_{\partial V}$ to $\Omega_{\Ff_{\rm EP}}$ (see \cite{KK}, end of Section 3 and
pages 555-556). Finally, as described before the statement, from \eqref{firstexp1} and \eqref{varconn2} we get $\mathfrak{h}^*\Ss^{0} = (s_V)_{\vert \Omega_{\Ff_{\rm EP}}}$.\hfill  $\Box$

\begin{remark}\label{CSvsCS}{\rm 1) In \cite{KK} the Chern-Simons
    bundle is defined over the variety of characters $X'(\partial V)$
    instead of the variety of augmented characters $X(\partial V)$;
    the coordinates on $X'(\partial V)$ are $1/2\pi i$ times the
    logarithms of the eigenvalues at a pair of meridian and longitude
    curve on $\partial V$. A fundamental domain for the action of
    $\mz\times \mz$ on coordinates is $([0;1/2]\times \mr)\times
    ([0;1]\times \mr)$. The rescaling factor $4\pi i$ at the end of
    the proof replaces these coordinates by the logarithms of the
    squared eigenvalues, so that a fundamental domain for the action
    of $\mz\times \mz$ on $A_{0,\Ff_{\rm EP}}(M)_\infty$ is
    $([0;2i\pi]\times \mr)^2$.

    2) The bundle $\mathcal{L}(M)$ is the restriction
    to $A_{0,\Ff_{\rm EP}}(M)$ of a bundle $\mathcal{L}(T^2)$ over
    $X(T^2)=X(\partial V)$ defined by the action \eqref{bundledef}. Of
    course $\mathcal{L}(T^2)$ is isomorphic to $\Ll_{\partial V}$. Its
    connection $1$-form \eqref{conn} has curvature $F_{\eta} =
    (-1/i\pi) d\log(\lambda)\wedge d\log(\mu)$, its Euler class is
    represented by $ (-1/2(i\pi)^2) d\log(\lambda)\wedge d\log(\mu)$,
    and its Euler number is $-2$. }
\end{remark}

\begin{remarks}\label{more-comments}{\rm {\it Problems and perspectives :}
\medskip

%Ici j'ai simplifiŽ :

\noindent  {\bf a)}   The $\mz^2$-equivariant analytic function $\Hh_1^{\Ff_{\rm EP},0}$ on $A_{0,\Ff_{\rm EP}}(M)_\infty$ is the classical analog of the regular rational functions $\Hh_N^{\Ff_{\rm EP},0,h_c,k_c}$ on $A_{0,\Ff_{\rm EP}}(M)_N$. One expects that a refinement of the $\mz^2$-action appears in the quantum world, involving the quantum torus algebra generated by $\hat l$, $\hat m$ such that $\hat l\hat m = \zeta\hat m\hat l$, where $\zeta$ is a primitive $N$-th root of $1$ and $(l,m)$ a geometric basis of $H_1(\partial V;\mz)$:
        
%. Indeed, take the $2N$-th power to kill the phase anomaly, and consider the complex vector space $\mathbb{L}_N$ freely generated by the functions $(\Hh_N^{0,h_c,k_c})^{2N}$ as $k_c$ varies (we drop the upperscript ``$\Ff_{\rm EP}$" to simplify notations; what follows extends obviously to the functions $(\Hh_N^{h_f,h_c,k_c})^{2N}$ with $h_f\ne 0$ -- see Remark \ref{extclass}). Define a representation of $\mz \times \mz \cong H^1(\partial V;\mz)$ on $\mathbb{L}_N$ as follows: associate to $\gamma' \in H^1(\partial V;\mz)$ the linear operator $\hat \gamma'$ on $\mathbb{L}_N$ defined by \begin{equation}\label{actionlin} \hat \gamma'. (\Hh_N^{0,h_c,k_c})^{2N}  = (\Hh_N^{0,h_c,k_c+\gamma'})^{2N}. \end{equation} From Theorem \ref{more_cusped} (2) it is easy to see that an element of the subgroup $(N\mz)^2$ acts on $\mathbb{L}_N$ by scalar multiplication: the scalars are given by the action on the functions on $A_{0,\Ff_{\rm EP}}(M)_N$ associated to the invariants $\alpha_N(V,\rho,(h,k))^{2N}$, the action on the reduced function $(\Hh_N^{h_{red}})^{2N}$ of Corollary \ref{on-reduced-ss} being trivial. A general class $\gamma'$ acts on $(\Hh_N^{h_{red}})^{2N}$  by a deck transformation. 
\begin{prob} Lift the functions $\Hh_N^{\Ff_{\rm EP},h_f,h_c,k_c}: A_{0,\Ff_{\rm EP}}(M)_N\ra \mc/\mu_{2N}$
  to $\mc$-valued functions defined on some covering space of $A_{0,\Ff_{\rm EP}}(M)_N$.
\end{prob} 

\noindent {\bf b)} By Theorem \ref{more_cusped} (2), the functions $\Hh_N^{\Ff_{\rm EP},h_f,h_c,k_c}$ yield $\mc/\mu_{2N}$-valued rational functions on a
covering space $\tilde X_0(M)_N$ of $X_0(M)$. On another hand, the field $\mc(\tilde X_0(M)_N)$ is a
finite extension of $\mc(X_0(M))$, and $\mc(X_0(M))$ is generated by the functions
of augmented characters associated to ${\rm hol}_l$ and ${\rm hol}_m$ by pull-back via the map $\mathfrak{h}: X_0(M)\ra A_0(M)$.
\begin{prob} Describe $\Hh_N^{\Ff_{\rm EP},h_f,h_c,k_c}$ in terms of functions of augmented
  characters, or as a meromorphic function on the smooth projective
  model of $\tilde X_0(M)_N$.
\end{prob}
For instance, a property that is obvious from the state sum
formulas is that the poles of $\Hh_N^{\Ff_{\rm EP},h_f,h_c,k_c}$ cover ideal
points of $X_0(M)$.

By a theorem of Bullock \cite{Bu}, the ring of
$SL(2,\mc)$-characters $\mc[X(M)]$ is isomorphic to
$K_{-1}(M)/\sqrt{0}$, where $K_{-1}(M)$ is the Kauffman bracket skein
module of $M$ specialized at $A=-1$, with its natural algebra
structure, and $\sqrt{0}$ its nilradical.
\begin{prob} Find skein theoretic constructions of $\Hh_N^{\Ff_{\rm EP},h_f,h_c,k_c}$.
\end{prob}
}
\end{remarks}
\section{QHFT partition functions}\label{I-PF}
For simplicity we consider only the QHFT patterns with topological 
support $(W,L, h_c)$, where $W$ is a closed oriented 
connected $3$-manifold. In that case only the bulk $c$-weight $h_c$ occurs. In order to build the analytic configurations of $(W,L, h_c)$ we
consider {\it quasi-regular} triangulations $T$ of $W$ (see Proposition \ref{totalGT}), and we require that $L$ is realized by a {\it Hamiltonian} subcomplex $H$ of the $1$-skeleton of $T$. Such a
pair $(T,H)$ is called a {\it distinguished} triangulation of $(W,L)$.

\begin{defi}\label{globalcQHFT} {\rm A {\it global charge} $c$ on $(T,H)$ is a rough global charge on $T$ satisfying the following additional global constraint
on the total edge charges:} {\it $C(e) = 0$ if $e$ is an edge of $H$, and $C(e) = 2$ otherwise.}
\end{defi}
Note that any global charge on $(T,H)$ encodes $H$
(i.e. the link $L$). Let us turn to the combinatorial encoding of $(W,L, h_c)$. Similarly to
Proposition \ref{top-support} and \ref{top-support2} we have (see \cite{Top,GT,AGT}):

\begin{prop}\label{on-charge} 
  (1) Every pair $(W,L)$ has quasi-regular distinguished
  triangulations $(T,H)$.

  (2) Every global charge $c$ on $(T,H)$ determines a bulk $c$-weight
  $h_c$ of $W$.

(3) For every bulk $c$-weight $h_c$ of $W$ and every distinguished
  triangulation $(T,H)$ of $(W,L)$ there is a global charge $c$ on
  $(T,H)$ with bulk $c$-weight equal to $h_c$.
\end{prop} 
Let us fix a quasi-regular distinguished triangulation $(T,H)$ of
$(W,L)$, a global charge $c$ on $(T,H)$, and a branching $b$ on $T$ (it
exists). 
\begin{prop}\label{on-flatt} 
  (1) Every point $[w;f]\in G_0(T,b)_\infty$ determines a character
  $\rho=\rho(w)$ and a bulk $f$-weight $h_f$, so that
  $(W,L,\rho,h_f,h_c)$ is a pattern with topological support
  $(W,L, h_c)$.

(2) For every pattern $(W,L,\rho,h_f,h_c)$ there is a point $[w;f]\in
  G_0(T,b)_\infty$ with holonomy $\rho$ and bulk $f$-weight $h_f$. 
\end{prop}

The construction of the QHI of $(W,L,h_c)$ is achieved by the following
Proposition; it is improved in Section \ref{sign} with respect to the QHI sign anomaly. The proof is similar to that of Theorem \ref{cusped-invariance}. For every odd $N\geq 3$ consider the analytic configurations
$\Aa_N(T,b,c)$. The defining equations of the algebraic variety $G_0(T,b,c)_N$ are:
\begin{itemize}
  \item $W'(e)=\zeta^{-1}$, for every edge $e$ of $T$ not contained in $H$;
\item $W'(e)=1$, for every edge $e$ of $T$ contained in $H$.
\end{itemize}
\begin{prop}\label{QHI-W} 
  For every pattern $\Pp= (W,L,\rho,h_f,h_c)$, every distinguished
  branched quasi-regular triangulation $(T,b,H)$ of $(W,L)$, every
  global charge $c$ on $(T,H)$ with bulk $c$-weight equal to $h_c$,
  and every point $[w;f]\in G_0(T,b)_\infty$ with holonomy $\rho$ and bulk
  $f$-weight $h_f$, the scalars
\begin{itemize}
\item $\Hh_1(\Pp):= \Hh_1(T,b)([w;f])$
\item $\Hh_N(\Pp):\equiv_{2N} \Hh_N(T,b,c)([w;f])$
\end{itemize}
do not depend on the choice of $(T,b,c)$ and $[w;f]$, and hence define invariants of $\Pp$.
\end{prop}
By definition, the analytic configuration $\Aa_N(W,L, h_c)$ is the family of analytic
configurations $\Aa_N(T,b,c)$ over the triples $(T,b,c)$ as in Proposition \ref{QHI-W}.

\begin{remark}\label{norm}{\rm The state sum normalization factor $a_N(T,\tilde b)$ includes a term $N^{-v}$. Other choices are possible. In \cite{LINK} we used instead $N^{2-v}$ in order to get Theorem \ref{QH=J}.}
\end{remark}

%We can compute the symmetrization factor $\alpha(T,b,c)$ of the function $\Hh_N(T, b,c)$ as in Theorem \ref{more_cusped} (3), by replacing $V$ by $\bar{V}$ in the proof. Since $\partial \bar{V}$ is a disjoint union of $2$-spheres the boundary weights vanish. Then we get:
%\begin{cor} We have $\alpha(T,b,c) \equiv_2 1$. Hence $\Hh_N(W,L,\rho)\equiv_{2N} \Hh_{N,red}(W,L,\rho)$.
%\end{cor}

For general QHFT patterns the constructions are more elaborated; in particular, one
requires that $H$ is Hamiltonian in the $1$-skeleton of $T$ from which the non-manifold vertices have been removed, and in the
normalization factor $N^{-v}$, $v$ denotes the number of manifold
(i.e. internal) vertices. Anyway, the results have basically the
same flavour. 
 
\section{$\Nn$-graph calculus}\label{N-calculus}
As every QH state sum is the total contraction of a tensor network supported by a $\Nn$-graph, one needs to describe the main features of a ``calculus'' based on such diagrams. 
%Recall that every $\Nn$-graph $\Gamma$ encodes a weakly branched triangulation $(T,\tilde b)$ of $\hat V$, hence also the induced pre-branched triangulation $(T,\sigma_b)$.
\subsection{$\Nn$-graphs representing the same weakly branched
  triangulation}\label{N-movesW}
Using the $\Nn$-graph decoding of Section \ref{I-triang}
it is easy to see that two $\Nn$-graphs $\Gamma$ and
$\Gamma'$ encode the same weakly branched triangulation $(T,\tilde b)$ if
and only if they are equivalent under the relation
generated by:
\begin{figure}[ht]
\begin{center}
 \includegraphics[width=6cm]{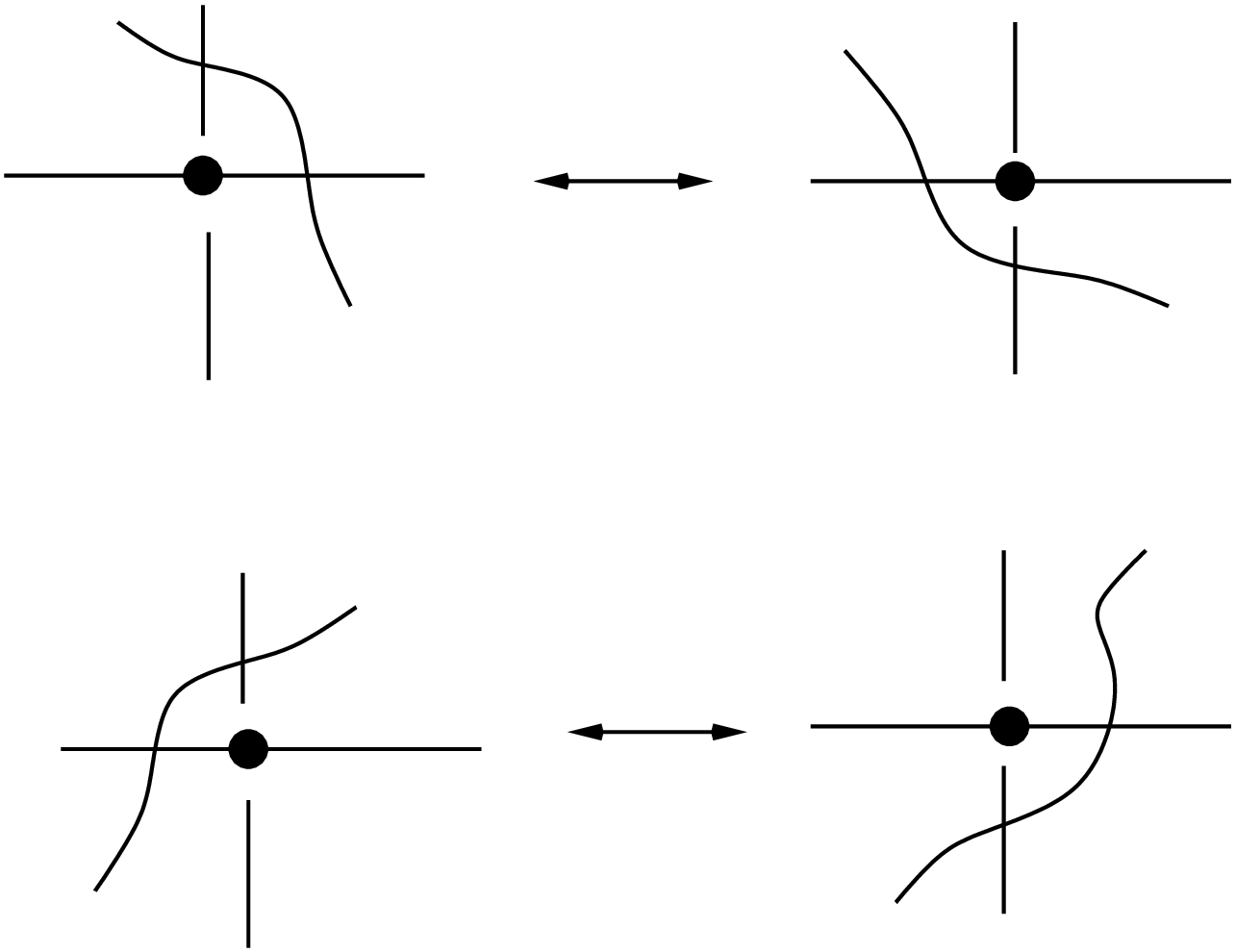}
\caption{\label{I-RIII} Some Reidemeister moves of type III for $\Nn$-graphs}
\end{center}
\end{figure} 

\begin{itemize}
\item Plane isotopy.
\item Switching the over/under arcs at an accidental crossing.   
\item The $\Nn$-graph versions of oriented Reidemeister moves, called
  {\it R-moves} (i.e. formally the usual ones if only accidental
  crossings are involved, or the Reidemeister move of type III shown in Figure \ref{I-RIII} (the orientation being understood) when
  a dotted crossing is involved).
\end{itemize}

\subsection{Changing the weak branching}\label{cwb}
Given a $\Nn$-graph $\Gamma$ representing a weakly
branched triangulation $(T,\tilde b)$, and an arbitrary change of weak branching $(T,\tilde b)\to (T,\tilde b')$, we want to describe a
systematic way to modify $\Gamma$ to a $\Nn$-graph
$\Gamma'$ representing $(T,\tilde b')$.  We begin with the local changes of branching on a given tetrahedron.

\smallskip

\subsubsection {The $S_4$ action $(\Delta,b)\to (\Delta,b_\beta)$} \label{S4action} 
Consider a branched tetrahedron $(\Delta,b)$ with $b$-ordered vertices
$v_0,v_1,v_2,v_3$ and $2$-facets $F_0,F_1,F_2,F_3$, as usual.  Let $S(J_4)$
be the symmetric group on the set $J_4=\{0,1,2,3\}$.  There is a
natural $1-1$ correspondence between the elements $\beta \in S(J_4)$
and the branchings $b_\beta$ on $\Delta$ such that $b=b_{\rm
  Id}$. The branching $b$ induces a branching $b(j)$ on each facet
$F_j$, that is a labelling of the vertices of $F_j$ by elements of $J_3=\{0,1,2\}$. For every $\beta \in S(J_4)$, also the branching
$b_\beta$ induces a branching $b_\beta(j)$ on $F_j$ (note that these
facets are still ordered with respect to $b$, not $b_\beta$), and the
transition from $b(j)$ to $b_\beta(j)$ is encoded by a permutation
$\epsilon_\beta(j)\in S(J_3)$. Clearly we have:
 
\begin{lem}\label{A3} 
(1)  The branching signs $*_{b_\beta}$ and $*_b$ coincide if and only if $\beta \in A(J_4)$.

(2) The $2$-facet $F_j$ has the same $b$- and $b_\beta$-transverse co-orientations if and only if $\epsilon_\beta(j)\in A(J_3)$.
\end{lem}

Let $C(b)$ denote a dotted $\Nn$-graph crossing that encodes
$(\Delta,b)$ as in Figure \ref{N-crossing}. We call $C(b)$ a {\it $2-2$ $\Nn$-tangle}; we consider it as a planar tangle properly embedded in a $2$-disk $D$, with endpoints on
$\partial D$ labelled by $J_4$. One can assume that the crossing is central to $D$, and that the arcs joining it to the endpoints are contained in two oriented diameters
of $D$.  We want to describe an algorithm that produces, for every $\beta\in
S(J_4)$,  a new decorated $2-2$ $\Nn$-tangle $C(b_\beta)$ in
$D$, with the same unoriented arc-germs as $C(b)$ at $\partial D$.  Hence
the arcs of $C(b)$ and $C(b_\beta)$ will be in natural $1-1$
correspondence and oriented by $b$ and $b_\beta$ respectively.  The
tangle $C(b_\beta)$ will verify the following properties:

\begin{itemize}
\item $C(b_\beta)$ contains one dotted crossing. According to the last
  Lemma, $*_{b_\beta}= *_b$ if and only if $\beta \in A(J_4)$. Hence,
  given $\beta$, we know the sign $*_{b_\beta}$.

\item $C(b_\beta)$ can replace $C(b)$ in
  every $\Nn$-graph $\Gamma'$ obtained after a change of
  weak branching $\tilde b \to \tilde b'$ acting as $b\to
  b_\beta$ on the corresponding branched tetrahedron $(\Delta,b)$ of $(T,\tilde b)$.  

\item Each of the four arc-germs of $C(b_\beta)$ at $\partial D$ 
is labelled by an element of $S(J_3)$, which will eventually contribute to 
the edge coloring of $\Gamma'$.
\end{itemize}

%Let us establish at first the following auxiliary construction.
\smallskip

\noindent {\bf An auxiliary construction}. Consider the two {\it basic
  $2-2$ tangles} in a $2$-disk $D'$, made by two {\it oriented} simple
arcs $a$, $a'$ that either are disjoint or intersect transversely at
one point, without any over/under information at the crossing point.
A $J_4$-labelling of the $4$ free endpoints of $a\cup a'$ is {\it
  admissible} if it respects the following conditions:

(a) The endpoints labelled by $1,3$ (resp. $0,2$) belong to
  different arcs and are both either the initial or the final endpoint
  of the corresponding arc.
\smallskip

(b) $3,2$ (resp. $1,0$) label the endpoints of the same arc.
\smallskip

\noindent These conditions are satisfied by the labels of the tangles of Figure
\ref{N-crossing}; recall that in this figure, $*_b=1$ (resp. $*_b=-1$)
if and only if $1,3$ label initial (resp. final) endpoints. By the
same rule we give every basic tangle an admissible
$J_4$-labelling, denoted by $\Bb$, and a sign $*_\Bb =\pm 1$. Note also that in
Figure \ref{N-crossing}, the arc with endpoints labelled
by $0,1$ passes over the other arc.  Now we convert $\Bb$ into a $2-2$
$\Nn$-tangle as follows:

\begin{figure}[ht]
\begin{center}
 \includegraphics[width=8cm]{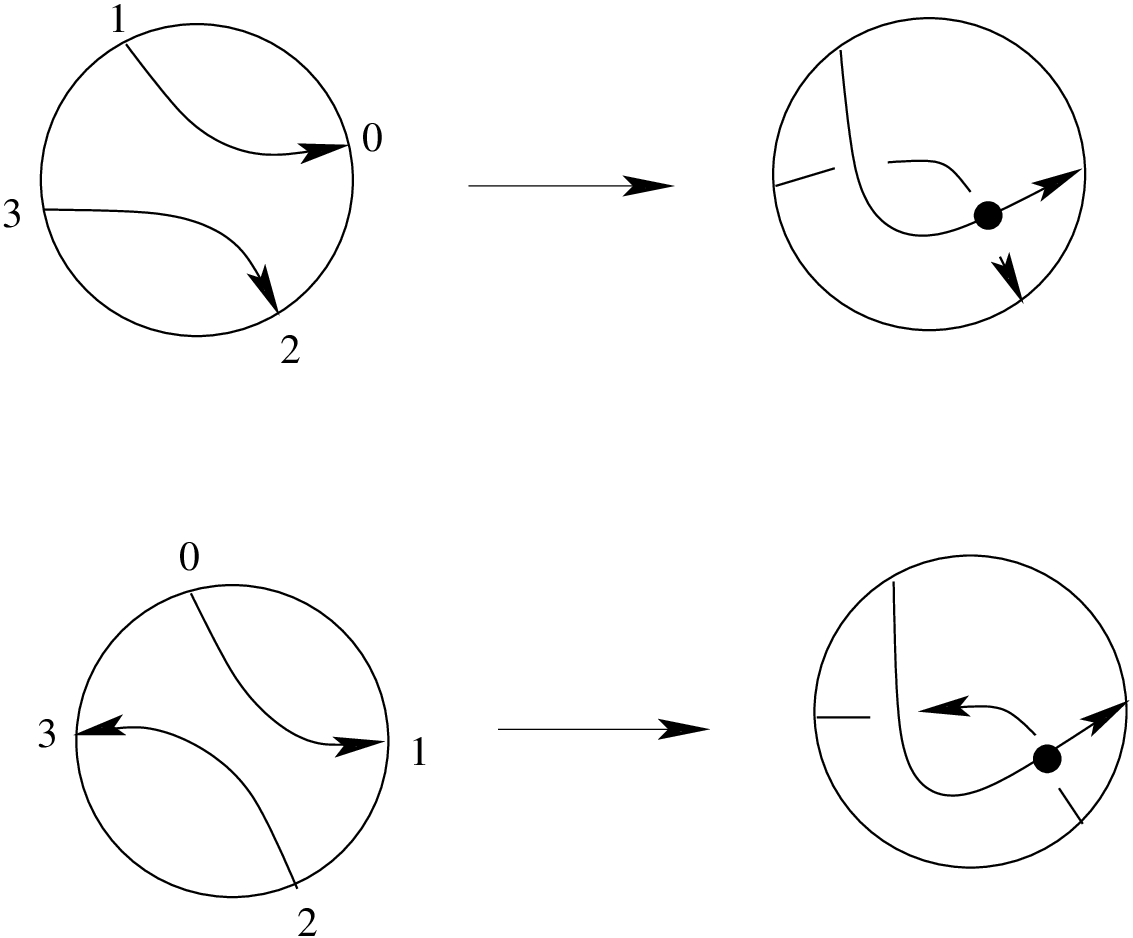}
\caption{\label{I-graph} Basic to $\Nn$-tangles.}
\end{center}
\end{figure}

\begin{figure}[ht]
\begin{center}
 \includegraphics[width=8cm]{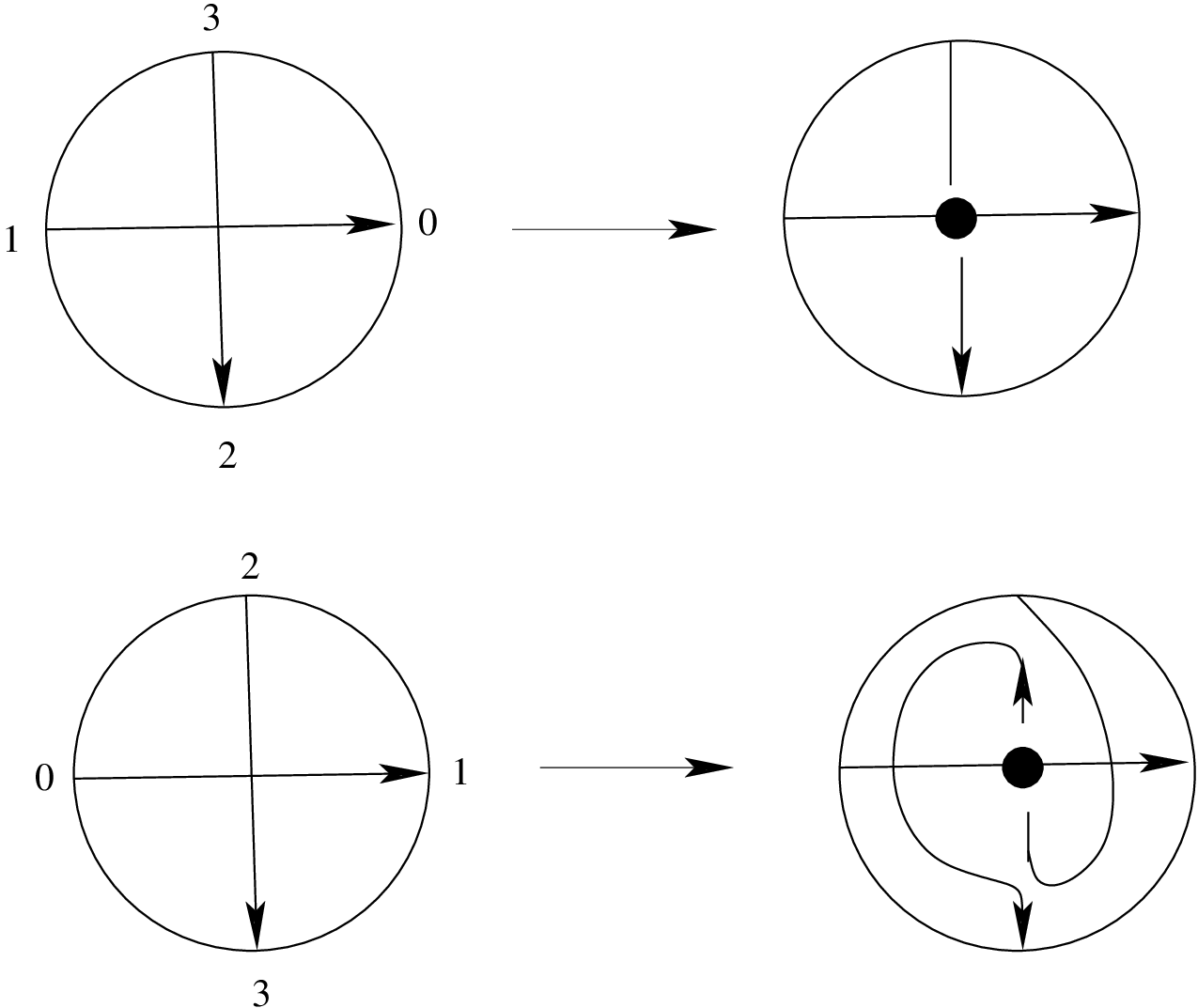}
\caption{\label{I-graph2} More basic to $\Nn$-tangles.}
\end{center}
\end{figure}

(c) If $a\cap a'=\emptyset$, perform an oriented Reidemeister
move of type II of the arc with labels $0,1$ over the other arc, thus
creating two crossings. Only one of them can be made into a dotted crossing so that its sign $*_b$ agrees with $*_\Bb$ (see some examples in Figure
\ref{I-graph}).  

\smallskip

(d) If $a$ and $a'$ cross at one point, let $A$ be the arc with $0,1$
labels, and $I$ a small open interval in the interior of the arc with
$2,3$ labels, such that $I$ contains the crossing point.  Let $A$ pass over
$I$, put a dot at the crossing point, and turn $I$ so that the crossing sign and $*_\Bb$ agree. Finally complete $I\cup A$ to an $\Nn$-tangle with the same endpoints as
$a\cap a'$, by introducing $0$ or $2$ accidental crossings on
opposite sides of $A$ with respect to $I\cap A$ (see some examples in
Figure \ref{I-graph2}).   \cvd 

\smallskip

{\bf The algorithm $C(b) \rightarrow C(b_\beta)$:}
\begin{enumerate}

\item Remove the interior of a smaller concentric sub-disk $D'$ of $D$. Each  
endpoint of $C(b)$, say $x_j$, $j\in J_4$, is connected by an oriented  
sub-arc to a point $x'_{\beta(j)}$ on $\partial D'$.

\item There is only one basic $2-2$ tangle $\Bb$ in $D'$ such that 
the $b_\beta$-labelling of the endpoints $x'_{\beta(j)}$ is admissible and
$*_\Bb=*_{b_\beta}$. As in the auxiliary construction, convert $\Bb$ into an $2-2$ $\Nn$-tangle in $D'$; glue it to the four sub-arcs $[x_j,x'_{\beta(j)}]$, and denote by $C'(b_\beta)$ the resulting $2-2$ $\Nn$-tangle in $D$. 

\item  
The $2-2$ $\Nn$-tangle $C(b_\beta)$ is obtained from $C'(b_\beta)$ by giving a label $\hat \epsilon_\beta(j)\in
  S(J_3) $ to each sub-arc $[x_j,x'_{\beta(j)}]$ as follows. Recall the permutation $\epsilon_\beta(j)\in S(J_3)$ defined before Lemma \ref{A3}; the $b$- and $b_\beta$-orientations of $[x_j,x'_{\beta(j)}]$ coincide one to each other iff $\epsilon_\beta(j) \in A(J_3)$).  Then set:
 \begin{itemize}
 \item $\hat \epsilon_\beta(j)=\epsilon_\beta(j)$, if the $b_\beta$-orientation of $[x_j,x'_{\beta(j)}]$ coincides with the orientation pointing towards the center of $D$;
\item $\hat \epsilon_\beta(j)=\epsilon_\beta(j)^{-1}$ otherwise.\cvd
\end{itemize}
\end{enumerate}
\subsubsection {Iterations}
Assume that  $\beta=\alpha \cdot \gamma$, where by definition
$\alpha \cdot \gamma := \gamma \circ \alpha$. Then we can produce 
$C(b_\beta)$ by iterating two times the preceding algorithm. We get a sequence of nested disks $$D''\subset D'\subset D,$$ 
a subdivision of the segment  $[x_j,x''_{\alpha \cdot \gamma(j)}]$ into sub-arcs
$$[x_j,x''_{\alpha \cdot \gamma(j)}]= 
[x_j,x'_{\alpha(j)}]\cup [x'_{\alpha(j)}, x''_{\alpha\cdot \gamma(j)}],$$ 
and a factorization of the permutation associated to $[x_j,x''_{\alpha \cdot \gamma(j)}]$:
$$\epsilon_\beta(j)= 
\epsilon_\alpha(j) \cdot \epsilon_\gamma(\alpha(j)). $$ 
The ``adjacent'' transpositions $(01), (12), (23)$ form a standard
set of generators of $S_4$. Let $\beta=\tau_1 \cdots \tau_k$ where
each $\tau_j$ is such a generator. Then we can produce again
$C(b_\beta)$ by iterating $k$ times the algorithm. We get a
sequence of nested disks
$$ D_k \subset .... \subset D_1 \subset D,$$ 
a subdivision of the arc $[x_j, x^{(k)}_{\beta(j)}]$ of  
$C(b)\setminus {\rm Int}(D_k)$ into sub-arcs
$$[x_j,x^{(1)}_{\tau_1(j)}]\cup[ \dots ]\cup[x^{(k-1)}_{\tau_1 \cdots \tau_{k-1}(j)}, 
x^{(k)}_{\beta(j)}],$$   
and a factorization of the permutation $\epsilon_\beta(j)\in S_3$ associated to the arc  $[x_j, x^{(k)}_{\beta(j)}]$ :

$$\epsilon_\beta(j)= \epsilon_{\tau_1}(j) \dots \epsilon_{  
 \tau_k}(\tau_1 \cdots \tau_{k-1}(j)) =: f_1 \cdots f_k. $$

\begin{lem} Every transposition $f_i\in S_3$ belongs to the set $\{{\rm Id}, (01),  
(12)\}$.
\end{lem}

\begin{figure}[ht]
\begin{center}
 \includegraphics[width=9cm]{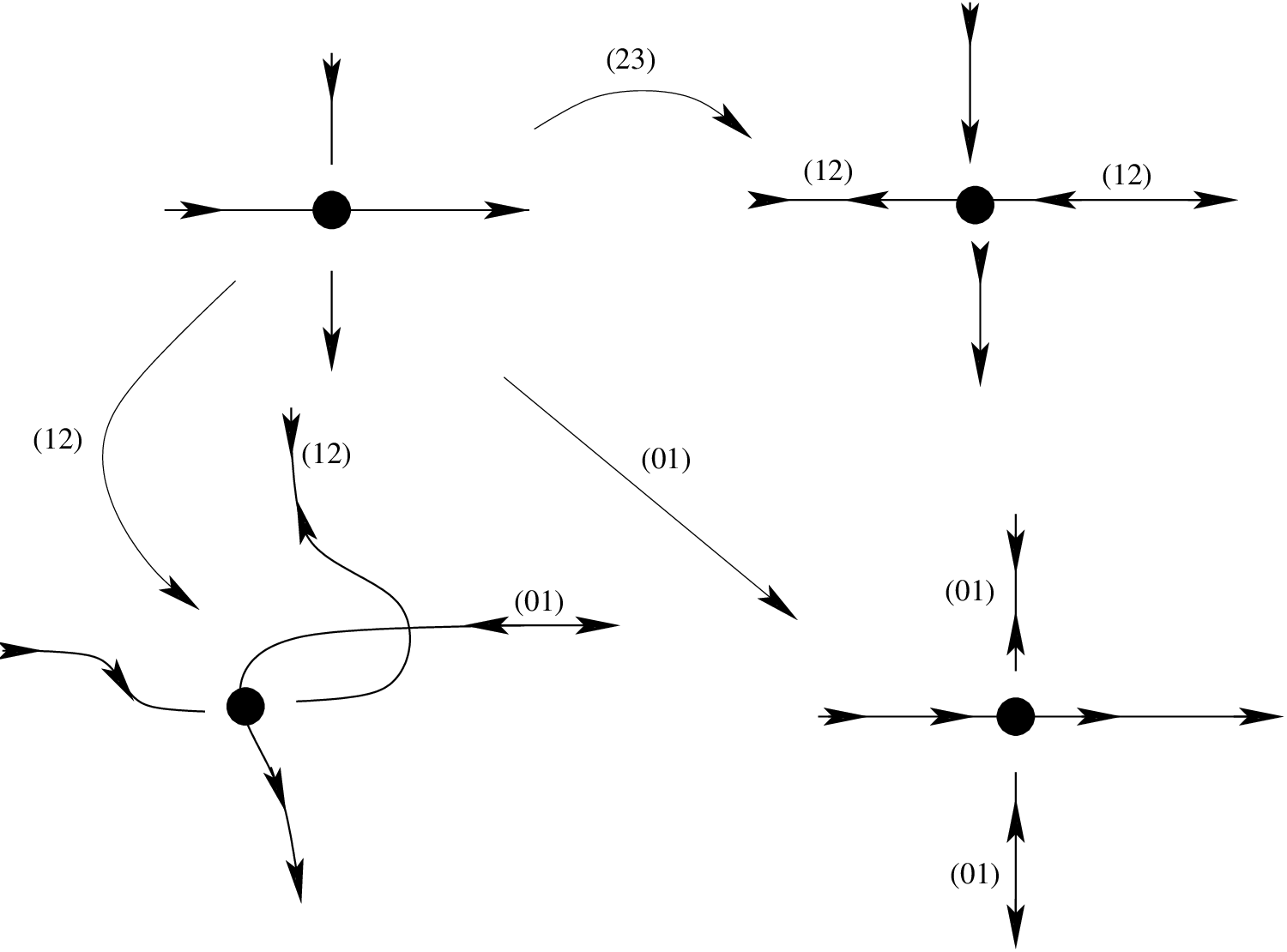}
\caption{\label{(ij)} $C(b_{(ij)})$, $*_b=+1$.}
\end{center}
\end{figure}

\Dim It is enough to check the case $\beta=\tau_j$ (see Figure \ref{(ij)}).
\cvd 
\medskip

For future application, it is convenient to associate a symbol to each sub-arc of $C(b_\beta)$, that indicates its label $f \in \{{\rm Id}, (01), (12)\}$ and depends on its $b$-orientation:
\begin{itemize}
\item If the $b$-orientation is outgoing (i.e. opposite to the one pointing towards  
the center of $D$), then let this symbol be $< {\rm Id} <$, or $<f>$ (a ``source'') if $f$ is a transposition; 
\item if the $b$-orientation is ingoing, then let this symbol be $>{\rm Id}>$, or $>f<$ (a ``pit'') if $f$ is a transposition. 
\end{itemize}
Under concatenation of sub-arcs, adjacent symbols are separated by $<<$ or $>>$. Here are examples: 
$$\cdots <f> >g< \cdots , \  \cdots >f< <g> \cdots , \ \cdots >id> >id> \cdots$$
$$\cdots <id< <id< \cdots , \  \cdots >id> >f< \cdots, \ \cdots <id< <f> \cdots,
\ {\rm etc}. $$  
Letting $?\in \{<<,\ >>\}$, in the case of $\epsilon_\beta(j) = f_1 \cdots f_k$ we get:
\begin{itemize}
\item {\it $ <f_1 ?\cdots  ? f_k<$ or $>f_1? \cdots ?f_k>$ if $b$ and $b_\beta$  induce on $[x_j, x^{(k)}_{\beta(j)}]$ the same orientation;
\item $<f_1? \cdots ?f_k>$ or $>f_1? \cdots ?f_k<$ otherwise.}
\end{itemize}
We can also make simplifications of the form:    
$$ ?s> >f< <f> >t? =\ ?s> >t?$$
and so on, so that the concatenation of symbols of sub-arcs eventually reduces to $> {\rm Id} >$, $< {\rm Id} <$, or alternating pits and sources carrying the transpositions $(01)$ and $(12)$. In Figure \ref{(ij)} we show graphically the result when $*_b=+1$ and $\beta$ varies in $\{(01), (12), (23)\}$; we dropped the notation ${\rm Id}$ for the symbol $> {\rm Id} >$. The case $*_b=-1$ is
similar.

\subsubsection {Oriented $C$-moves}\label{orientedC}
We call $C(b)\to C(b_\beta)$ a $C$-move. If $C(b)$ and  $C(b_\beta)$ define the same pre-branching then we say that
it is an {\it oriented $C$-move}. Set 
$$\sigma=(0123) =  (23)\cdot (12)\cdot (01)\ ,\ \tau = (02) \cdot (13).$$
$\sigma$ generates a subgroup $\CG$ of $S(J_4)$ isomorphic
to $(\Z/4\Z,+)$ via the map $\sigma \mapsto 1$.
\begin{lem}\label{orientedCmove} For both $*_b=\pm 1 $, $\beta \in S_4$ 
  induces an oriented $C$-move $C(b)\to C(b_\beta)$ iff $\beta \in \CG$.
\end{lem}  
\Dim Looking at Figure \ref{Branched_Delta} anti-clockwise from the top left, and denoting by $b$ the first
branching (so that $*_b=+1$), the others are successively $b_\tau$,
$b_{\sigma \cdot \tau}$, $b_{\tau \cdot \sigma \cdot \tau}$,
$b_{\sigma^{-1}\cdot \tau \cdot \sigma \cdot \tau}=b$. Hence $\beta \in
\CG$. One can do similarly when $*_b=-1$. \cvd
\smallskip

\begin{figure}[ht]
\begin{center}
 \includegraphics[width=9cm]{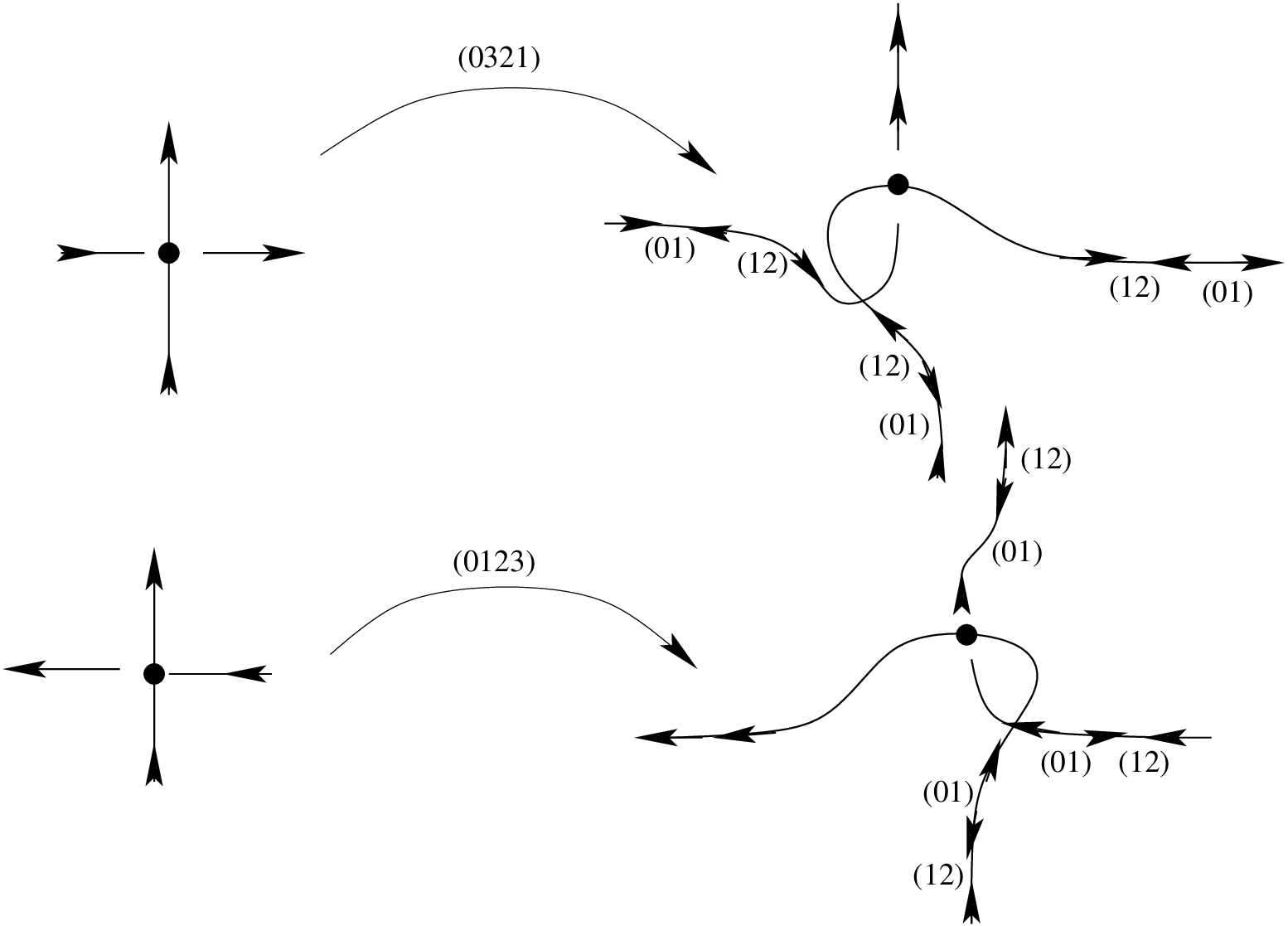}
\caption{\label{or-C-gen} Generating oriented $C$-moves.}
\end{center}
\end{figure}

Figure \ref{or-C-gen} shows the oriented $C$-moves $C(b)\to
C(b_{\sigma^{-1}})$ for $*_b= +1$ and $C(b)\to C(b_{\sigma})$ for $*_b=-1$.

\subsubsection{Globalization} Let $\Gamma$ be an $\Nn$-graph
encoding $(T,\tilde b)$, and $\tilde b'$ another weak branching on $T$. We want to produce an $\Nn$-graph $\Gamma'$ encoding $(T,\tilde b')$.  
For every abstract branched tetrahedron
$(\Delta_j,b_j)$ of $(T,\tilde b)$ there is a permutation $\beta_j\in S(J_4)$ such
that the restriction of the transit $\tilde b \to \tilde b'$ to $\Delta_j$ is $b_j\to (b_j)_{\beta_j}$. Applying at every crossing of
$\Gamma$ the algorithm $C(b_j)\rightarrow C((b_j)_{\beta_j})$,
we get a diagram $\Gamma'$.  Forgetting the orientations, there is a
natural $1-1$ correspondence between the edges of $\Gamma$
and $\Gamma'$. Let
$e$ be an edge of $\Gamma$ with color
$r(e)$. Denote by $\bar e$ the corresponding edge of $\Gamma'$. According to the last step of the algorithm, $\bar e$ is labelled by ``initial'' and ``final'' permutations $\hat
\epsilon^i(\bar e), \hat \epsilon^f(\bar e)\in S(J_3)$. There is an alternative:
\begin{itemize}
\item Both $\hat \epsilon^i(\bar e)$ and $\hat \epsilon^f(\bar e)$
  belong to $A(J_3)$ (in such a case $e$ and $\bar e$ have compatible
  orientations); then the color $r(\bar e)$ is defined by
$$ \hat \epsilon^i(\bar e)\cdot (012)^{r(e)}\cdot \hat \epsilon^f(\bar e)=
(012)^{r(\bar e)} .$$
\item Both $\hat \epsilon^i(\bar e)$ and $\hat \epsilon^f(\bar e)$ do
  not belong to $A(J_3)$ (in such a case $e$ and $\bar e$ have
  opposite orientations); then the color $r(\bar e)$ is defined by
$$ \hat \epsilon^i(\bar e)\cdot (012)^{-r(e)}\cdot \hat \epsilon^f(\bar e)=
(012)^{r(\bar e)} .$$ 
\end{itemize}

\noindent The oriented $C$-moves have the nice feature to be local and
independent one to each other. On the contrary, a family of 
permutations $\beta_j\in S(J_4)$ as above must satisfy non trivial global constraints to induce a change of weak branching $\tilde b \to
\tilde b'$.  Fortunately, the following local/global calculus on
$\Nn$-graphs (introduced in \cite{BP2}) covers the general case.

\begin{lem}\label{local/global} Any change of $\Nn$-graphs $\Gamma\to \Gamma'$ corresponding to a change of weak branching $(T,\tilde b)\to (T,\tilde b')$ is a composition of a finite sequence of moves $\Gamma_j\to \Gamma_{j+1}$ of 
  the following types:
\begin{enumerate}
\item The moves of Subsection \ref{N-movesW};
\item The oriented $C$-moves;
\item (A non-local {\rm ``circuit'' move}). Let $\gamma$ be a simple circuit of $\Gamma_j$, ie. a circuit with one component. Assume that each edge of $\gamma$ is the upper strand at each crossing through which $\gamma$ passes. Then $\Gamma_{j+1}$ is obtained by
  acting on each of these crossings with the transposition $(23)\in
  S_4$. Equivalently, $\Gamma_{j+1}$ is obtained by
  reversing the orientation of $\gamma$ and keeping the $\mz/3\mz$-color $r(e)$ unchanged at every edge $e$ of ${\rm Sing}(P)$.
\end{enumerate}
\end{lem}
\Dim The equivalence between the two descriptions of the non-local
``circuit'' move is clear. Let $\Gamma\to \Gamma'$
be as in the statement. The set $\mathcal{S}$ of edges of ${\rm Sing}(P)$ where the
pre-branchings associated to $\tilde b$ and $\tilde b'$ disagree is
the union of non-overlapping simple circuits $\gamma_1$, $\ldots, \gamma_r$ oriented by $\tilde
b$. Using oriented $C$-moves one can modify $\Gamma$ so that $\gamma_1$ eventually passes over at each crossing where it goes through. Then change its orientation. Doing similarly with $\gamma_2$, $\ldots, \gamma_r$, one ends up with $\mathcal{S}= \emptyset$.\cvd

\smallskip 

\subsection{Transits}\label{MP-bub}
\begin{figure}[ht]
\begin{center}
 \includegraphics[width=9cm]{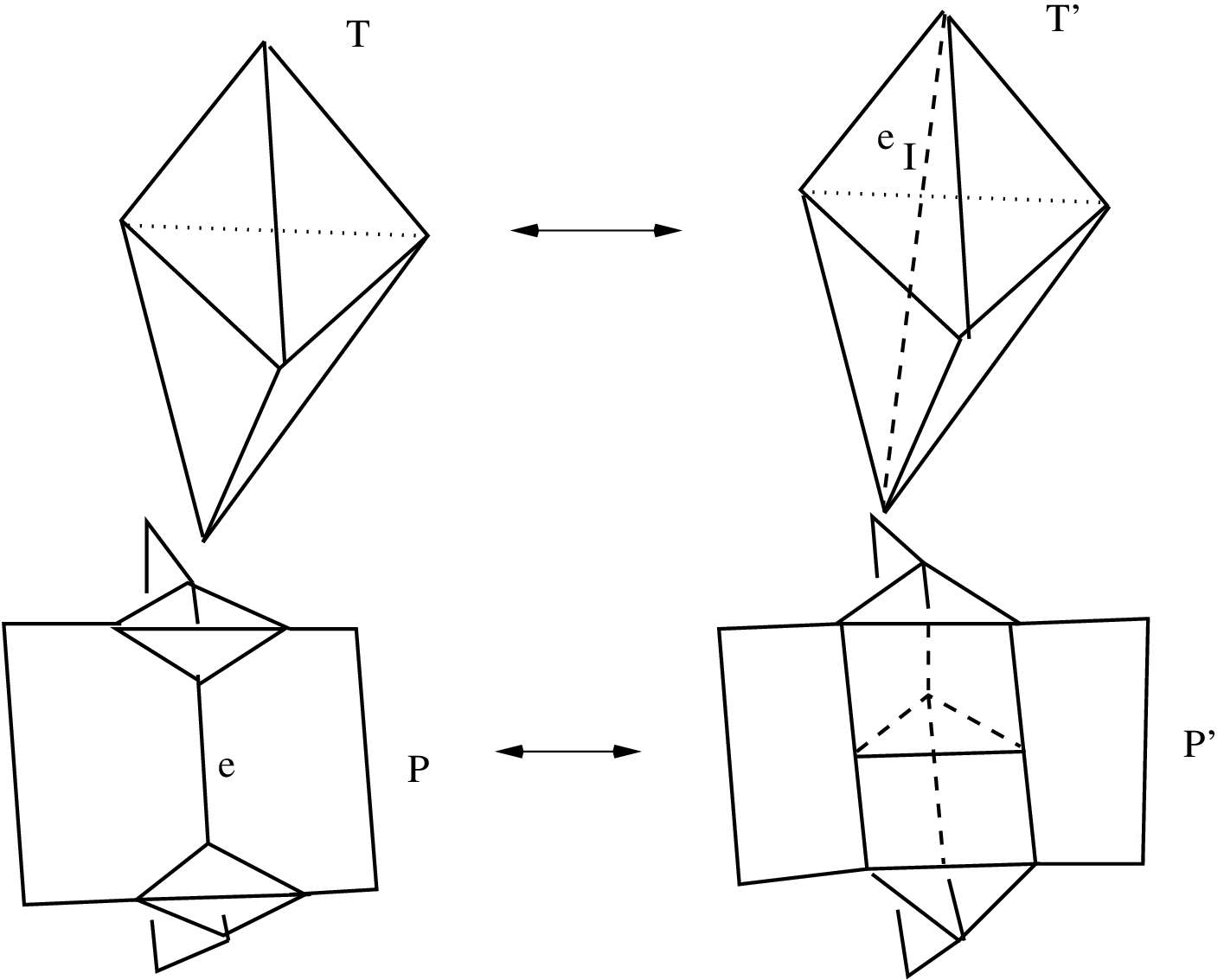}
\caption{\label{IntrinsicMP} The MP move.}
\end{center}
\end{figure}

\begin{figure}[ht]
\begin{center}
 \includegraphics[width=9cm]{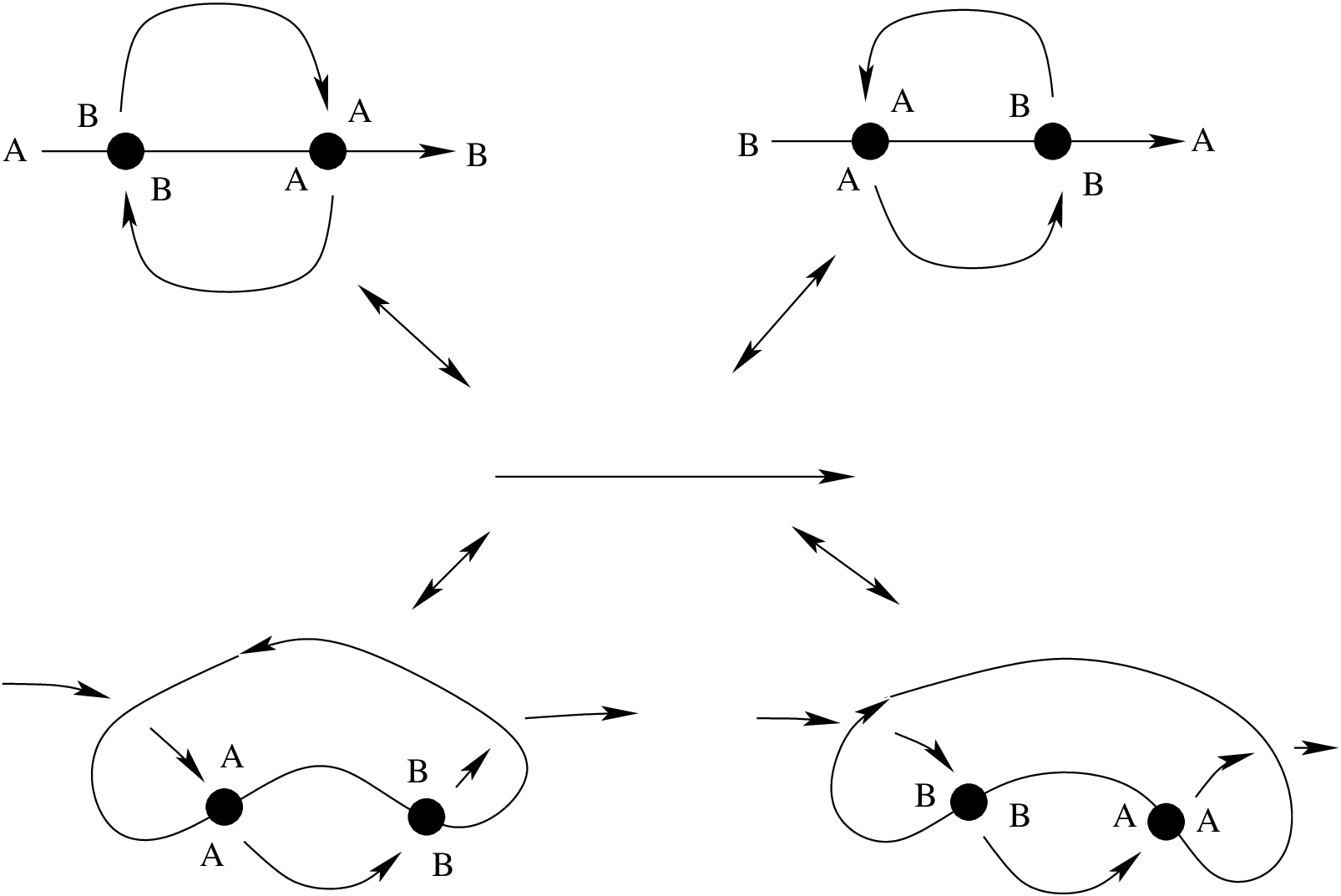}
\caption{\label{bubbles} Bubble $b$-transits.}
\end{center}
\end{figure}

It is well-known that two ``naked'' triangulations $T$ and $T'$ of
$\hat V$ can be connected by a finite sequence of $3$-dimensional
Pachner's moves, also called MP and bubble moves.  The MP
moves,  also called ``$2\leftrightarrow 3$'' moves, are illustrated in Figure \ref{IntrinsicMP} in terms of
triangulations and dual spines. Branched versions of the bubble moves are shown in Figure \ref{bubbles} in terms of $\Nn$-graphs.

%Let us enhance theses moves for structured triangulations.

\begin{defi}\label{p-b-transits} {\rm %Given any move (either bubble or MP, positive or negative) $T\to T'$ there is a maximal subset $\Ff$ of the $2$-faces of $T$ that ``persist'' in $T'$. 
Let $(T,\sigma)$ and $(T',\sigma')$ be pre-branched triangulations of
    $\hat V$. A move $T\to T'$ (either bubble or
    MP, positive or negative) lifts to a {\it
      pre-branching transit} $(T,\sigma)\to (T',\sigma')$ if and only
    if $\sigma$ and $\sigma'$ coincide on every common $2$-face of $T$ and $T'$. We say that $T\to T'$ lifts to a $\tilde b$-{\it transit} $(T,\tilde b) \ra
    (T',\tilde{b}')$ of weakly-branched triangulations if it induces a pre-branching transit
    $(T,\sigma_{\tilde b})\to (T',\sigma_{\tilde b'})$.  If $(T,b)$
    and $(T',b')$ are branched triangulations, then $(T,b) \ra (T', b')$ is
    a $b$-{\it transit} if $b$ and
    $b'$ coincide on every common edge of $T$ and $T'$.}
\end{defi}

\begin{figure}[ht]
\begin{center}
 \includegraphics[width=8cm]{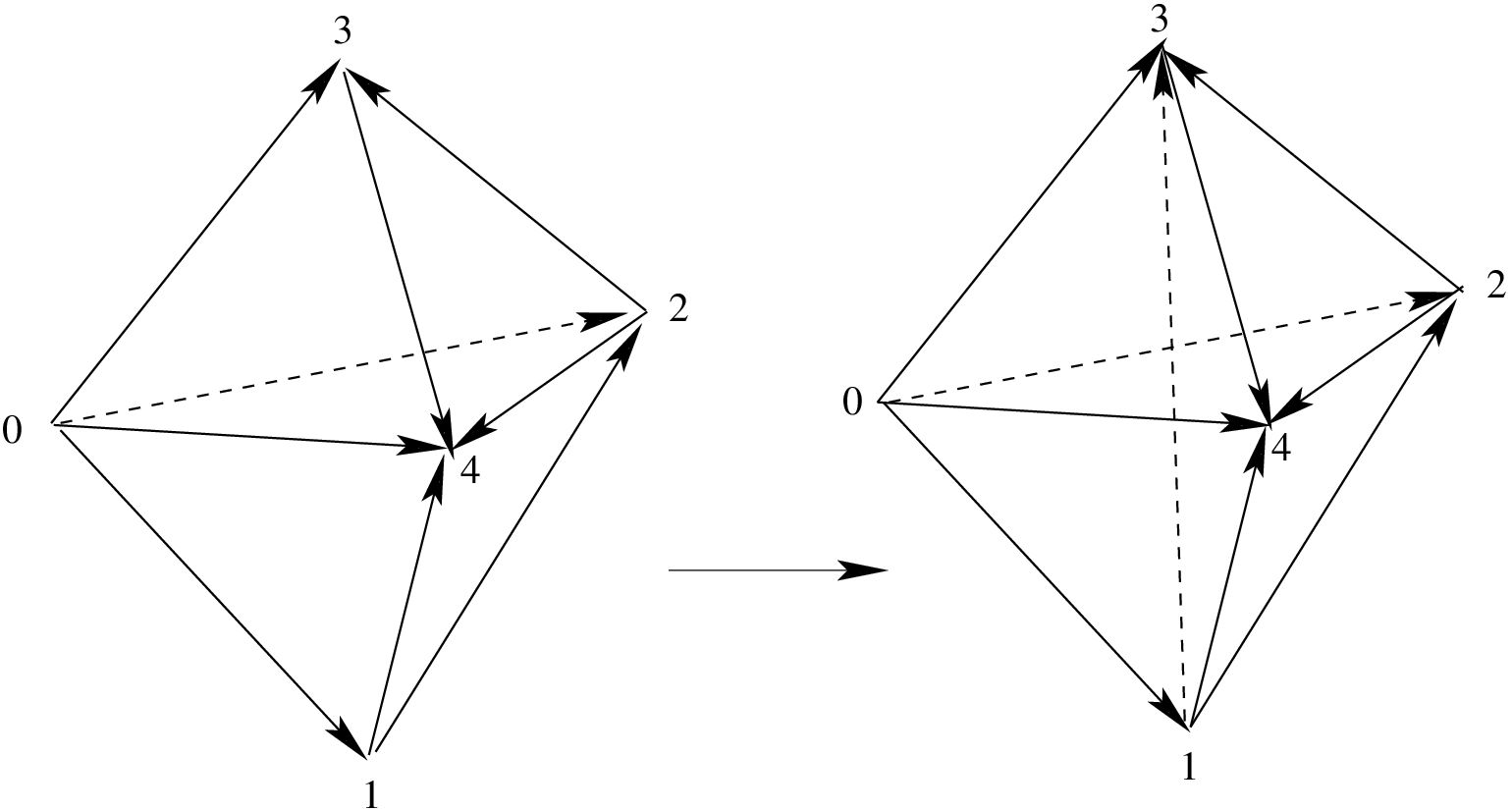}
\caption{\label{5-terms} Schaeffer's $b$-transits.}
\end{center}
\end{figure}

\begin{figure}[ht]
\begin{center}
 \includegraphics[width=8cm]{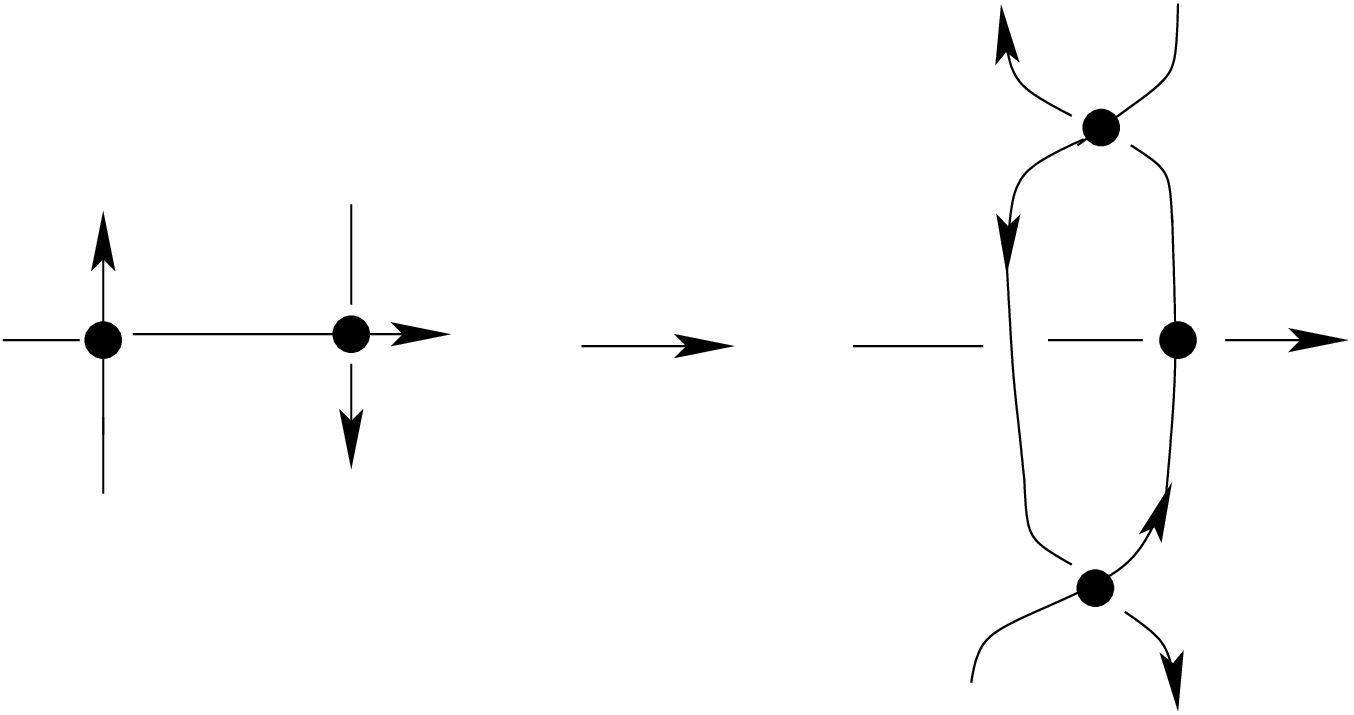}
 \caption{\label{I-MP-b-transit} Normal $\Nn$-graph Schaeffer's
   $b$-transit.}
\end{center}
\end{figure}

Figure \ref{bubbles} shows the different bubble $b$-transits in terms of
normal $\Nn$-graphs (recall Remark \ref{o-graph}), and the decorations $A$, $B$
of the dual square edges according to Figure \ref{I-dec}.  Figure \ref{5-terms} and \ref{I-MP-b-transit} show an example of MP
$b$-transit in terms of branched triangulations and normal $\Nn$-graphs; note that $*_b=+1$ for all $3$-simplices (resp. dotted crossings). It is called the {\it
  Schaeffer} $b$-transit, and plays a distinguished role in the study
of matrix dilogarithms (see Section \ref{sign}).

The following Lemma will be useful in applications to the QH state
sums. Denote by $\cong_{pb}$ the equivalence relation on the set of weakly branched triangulations, defined by $\tilde b$-transits $(T,\tilde b)\to (T,\tilde b')$.

\begin{lem}\label{b-vs-pb} 
  The equivalence relation $\cong_{pb}$ is generated by the following two transformations:
\begin{itemize}
\item Change the weak-branching $(T,\tilde b)\to (T,\tilde b')$ by moves preserving the pre-branching (ie. dually do oriented $C$-moves, see Lemma \ref{local/global}). 
\item Perform $\tilde{b}$-transits $(T,\tilde b)\to (T',\tilde b')$ which look like $b$-transits at the sub-patterns of tetrahedra involved in the move.
\end{itemize}
\end{lem}

\begin{figure}[ht]
\begin{center}
 \includegraphics[width=9cm]{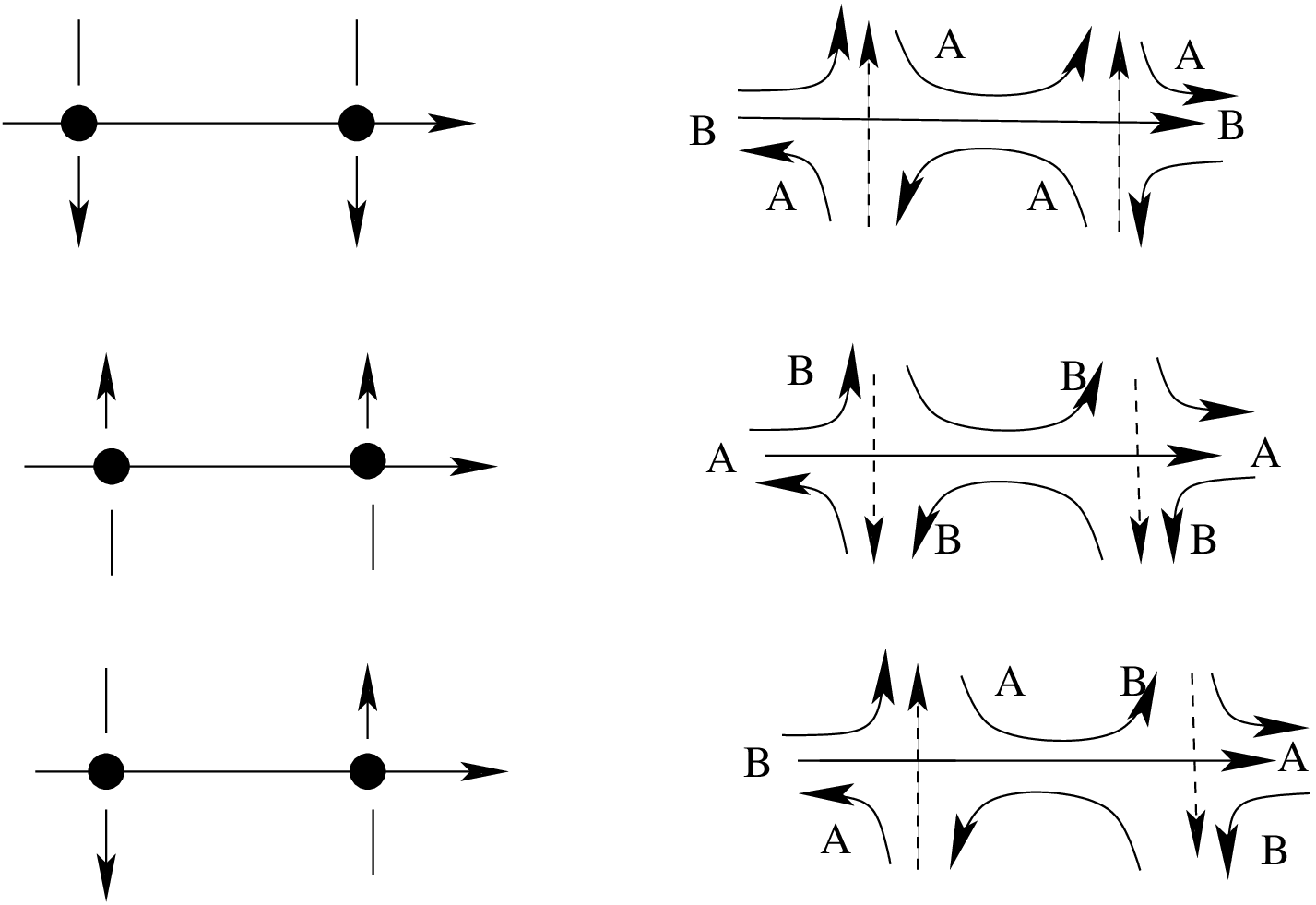}
\caption{\label{pb-vs-b} Branched realization of MP $pb$-transits.}
\end{center}
\end{figure}

\Dim The proof for bubble moves is easy and basically illustrated in Figure
\ref{bubbles} in  terms of $\Nn$-graphs. Consider a positive MP $\tilde{b}$-transit $(T,\tilde
b)\to (T',\tilde b')$, and the two tetrahedra of $T$ modified by the transit.  At their common $2$-face $F$ one of the following possibilities is realized:

(1) there are exactly two square edges which are both {\it not monochromatic}
(i.e. each one inherits {\it different} $A$,$B$ labels from the pre-branchings of the two
tetrahedra);

(2) there is exactly one square edge which is monochromatic (either $A$ or $B$), and the
other two edges are monochromatic (either $B$ or $A$ coupled with the empty
label). 

Moreover, this information at $F$ determines completely the co-orientations of the three $2$-faces produced by the transit, and hence determines the pre-branching transit $(T,\sigma_{\tilde b})\to (T',\sigma_{\tilde b'})$. Then, it is enough to realize the above possibilities by means of $b$-transits (changing the weak branchings before and/or after the transit, if necessary,  by moves preserving the pre-branching).  Figure \ref{pb-vs-b} shows it
in terms of (decoded) normal
$\Nn$-graphs. Clearly this realization is not unique.  \cvd

\section{State sum invariance over weakly branched triangulations}\label{I-enhanced}
Let $(T,\tilde b,w,f,c)$ be a weakly branched QH triangulation of a pattern $\Pp$; $w$ and $f$ verify the defining equations of the spaces in Definition \ref{defiAC}, and $c$ those of Definition \ref{globalc} (resp. \ref{globalcQHFT}) if $\Pp$ is a pattern over a cusped manifold (resp. a QHFT pattern).

\smallskip

We want to prove that $\Hh_N(\Pp):\equiv_{N} \pm \Hh_N(T,\tilde b,w,f,c)$ is an invariant of $\Pp$. In the case when $(T,b)$ is a branched triangulation we achieved this result in \cite{Top, GT, AGT} by introducing the notion of {\it QH transit} $(T,b,w,f,c)\ra (T',b',w',f',c')$. A QH transit is such that: $(T,b)\to (T',b')$ is a $b$-transit, the tuples $(w,f,c)$ and $(w',f',c')$ coincide on the ``common'' tetrahedra of $T$ and $T'$, and for every common edge $e$ we have (see Section \ref{glob} for the notations, where $W_T(e)$ is $W(e)$, and so on)
\begin{equation}\label{eqQHt}
W_T(e) = W_{T'}(e)\ ,\ L_T(e)= L_{T'}(e)\ ,\ C_{T}(e) = C_{T'}(e).
\end{equation}
QH triangulations related by QH transits encode a same pattern $\Pp$. Moreover:\begin{itemize}
\item $\Hh_N(T,b,w,f,c)$ is invariant under any QH transit.
\item Any two branched QH triangulations of $\Pp$ are related by a finite sequence of QH transits and branching changes.
\item $\Hh_N(T,b,w,f,c)$ is invariant under any change of the branching $b$.
\end{itemize}

Clearly these three facts prove the invariance of $\Hh_N(\Pp)$. We follow the same strategy when $\tilde b$ is an arbitrary weak branching. Since the QH transits are local transformations, we define QH transits $(T,\tilde b,w,f,c) \ra (T,\tilde b',w',f',c')$ by requiring that the $\tilde b$-transit $(T,\tilde b) \ra (T',\tilde b')$ restricts to a genuine $b$-transit on the tetrahedra modified by the move. Then, trivially $\Hh_N(T,\tilde b,w,f,c)$ is again invariant under QH transits, and using Lemma \ref{b-vs-pb} and the results of \cite{Top, GT, AGT} it is easy to see that any two weakly branched QH triangulations of $\Pp$ are related by a finite sequence of QH transits and changes of weak branching.  So it remains to show: 
\begin{prop}\label{wbinvariance}
 For any change of weak branching $\tilde b \to \tilde b'$ we have $$\Hh_N(T,\tilde b,w,f,c)\equiv_{N} \pm \Hh_N(T,\tilde b',w,f,c).$$
\end{prop} 
In order to develop the proof we need to fix some conventions regarding the contraction of tensors, and to describe the relations between tetrahedral tensors related by a change of branching. This is the content of the two following subsections.
\subsection{Formal conversion of tensor networks}\label{formal-conv}
Recall the notation $V=\C^N$. We use the canonical isomorphism $V\ra V^*$ sending the standard basis to the dual basis. Let $A\in {\rm End}(\C^N\otimes \C^N)$ be associated to a $2-2$
$\Nn$-tangle as in Figure \ref{tensorindex}, according to the
conventions of Section \ref{LocalAC}. Let $\beta \in
S(J_4)$. Assuming that $*_b=1$, under the natural isomorphisms
\begin{align} {\rm Hom}(V_3\otimes V_1,V_2\otimes V_0) & \cong
  V_2\otimes V_0 \otimes (V_3\otimes V_1)^* \label{tensiso1} \cong
  V_2\otimes V_0 \otimes V_1^*\otimes
  V_3^*
\end{align}
we see that $A$ belongs to $V_2\otimes V_0
\otimes V_1^*\otimes V_3^*$, where as usual $V_j$ is the copy of
$\C^N$ associated to the $j$-th face of $(\Delta,b)$.  Consider
$W_2\otimes W_0 \otimes W_1^*\otimes W_3^*$ if $*_{b_\beta}=1$, and
$W_3\otimes W_1 \otimes W_0^*\otimes W_2^*$ if $*_{b_\beta}=-1$, where
$W_j$ is the copy of $\C^N$ associated to the $j$th-face  according to
$b_\beta$. By using ${\rm Id}:V_j\to W_{\beta(j)}$ or the
canonical isomorphisms $V_j \to W^*_{\beta(j)}$, we get
further canonical isomorphisms
$$V_2\otimes V_0 \otimes V_1^*\otimes V_3^*\cong W_2\otimes W_0 \otimes
W_1^*\otimes W_3^*\quad {\rm and}\quad V_2\otimes V_0 \otimes V_1^*\otimes V_3^*\cong W_1 \otimes W_3
\otimes W_0^*\otimes W_2^* .$$

Denote by $A_\beta$ the operator supported by $(\Delta,b_\beta)$,
defined as the image of $A$ via such an isomorphism. We call $A_\beta$
the {\it formal $\beta$-conversion} of $A$. There are explicit
identities between the matrix elements of $A$ and $A_\beta$. Following
our conventions for the index positions, we get for example
$$A^{p,q}_{s,t} = (A_{(13)(02)})^{q,p}_{t,s} = 
(A_{(23)})^{p,t}_{s,q} =(A_{(02)})^{s,t}_{q,p} = 
(A_{(01)})^{s,q}_{p,t}= (A_{(12)})^{s,p}_{t,q} .$$ Let $\Aa$ be a QH tensor network supported by a weakly branched triangulation $(T,\tilde b)$. For any change of weak branching $\tilde b \to \tilde b'$, the formal conversions of the 
tetrahedral tensors of $\Aa$ match to produce a new tensor network
$\Aa_{\tilde b'}$, supported by $(T,\tilde b')$. The state sums associated to $\Aa$ and $\Aa_{\tilde b'}$ being total contractions, it is tautologically evident
that they take the same values.

\begin{figure}[ht]
\begin{center}
 \includegraphics[width=5cm]{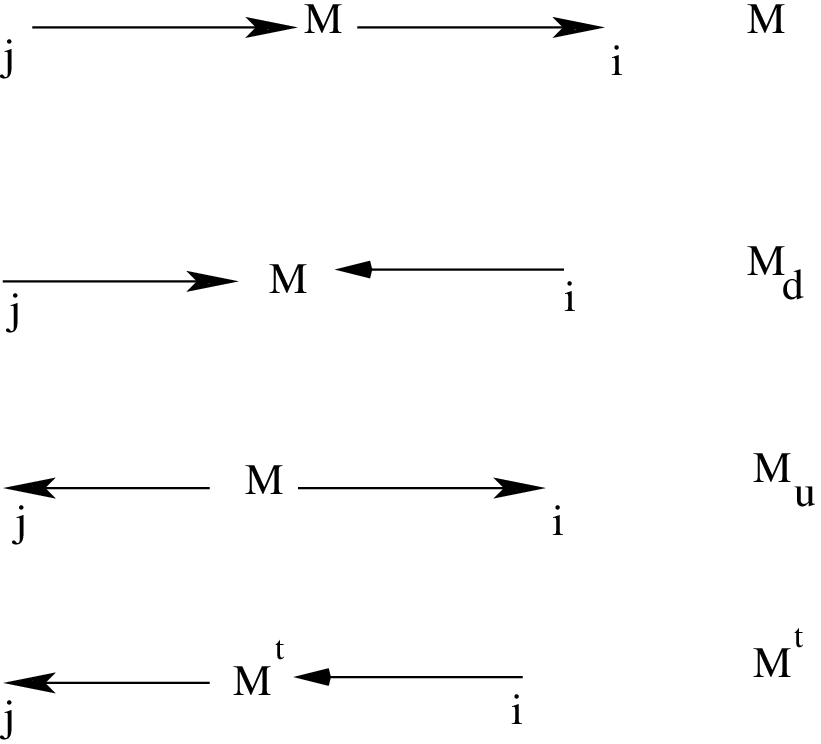}
\caption{\label{convM} Formal conversion of a square matrix.}
\end{center}
\end{figure}

\begin{figure}[ht]
\begin{center}
 \includegraphics[width=9cm]{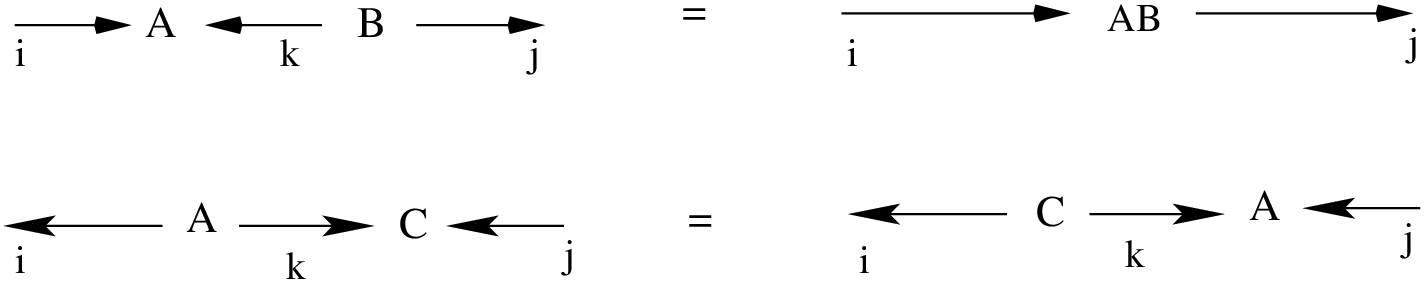}
\caption{\label{convM2} Relations for $M$-conversion.}
\end{center}
\end{figure}

The formal conversion of tensors can be widely applied. Here is
a second example that we use below. Let $M=(M^i_j)$ be any $N\times N$ matrix. We interpret $M$ as an endomorphism of $V=\C^N$.  There are two formal
conversions of $M$, $M_u$ and $M_d$, which are bilinear forms
on $V^*$ and on $V$ respectively, defined by:
\begin{equation}\label{Mconv}
(M_u)^{i,j}=M^i_j=(M_d)_{i,j} .
\end{equation}
Note that we have also $(M_u)^{i,j}=(M^t)^j_i=(M_d)_{i,j}$, where $M^t$ is the transpose matrix of $M$, that is, the adjoint
endomorphism. Figure \ref{convM} shows a graphical encoding of these identities. If $M$ is symmetric, the position
of the indices $i,j$ is immaterial. Such formal conversions satisfy a few
relations shown in Figure \ref{convM2}, where $A$, $B$ and
$C$ denote $N\times N$ matrices and we assume that $AC=CA$.  Note that
in every case the index $k$ is traced out, and the top/right picture
represents the matrix element of a composite endomorphism 
$$A\cdot B := B\circ A.$$
\subsection{The relations between $\Rr_N(\Delta,b,w,f,c)$
and $\Rr_N(\Delta,b_\beta,w',f',c')$} Consider a change of branching $b\ra b_\beta$, $\beta \in S(J_4)$. It is enough to treat the case where $\beta$ is one of the standard generators $(01), (12), (23)\in S(J_4)$. Consider the algorithm $C(b)\to C(b_{\beta})$ of Section \ref{cwb}. Associate a tensor network to $C(b_{\beta})$ as follows: 
\begin{itemize}
\item Associate $\Rr_N(\Delta,b_{\beta},w',f',c')$ to the dotted crossing of $C(b_{\beta})$.
\item Using the matrices $S$ and $T$ of Section \ref{glob}, replace $>(01)<$ with $>T<$, $ >(12)<$ with $>S<$, $<(01)>$ with $<T^{-1}>$, $<(12)>$  with $<S^{-1}>$, and 
$>{\rm Id}>$ with $>I>$.
\end{itemize}
Contracting the tensor network we get a tensor $\hat \Rr_N(T,b_{\beta},w',f',c')$ which has the same type as $\Rr_N(\Delta,b,w,f,c)$, since both have the same oriented
arc-germs at $\partial D$. The next Proposition rephrases the content of Lemma 3.3 of \cite{Top} and Corollary 5.6 of \cite{GT}. Matrix indices are raised or lowered using the standard inner product of $\mc^N$ (see Section \ref{formal-conv}). 
\begin{prop}\label{3.3} 
If $\beta = (01),(12),(23)$, then
$$\Rr_N(T,b,w,f,c) \equiv_{N} \pm \hat \Rr_N(T,b_{\beta},w',f',c').$$
More precisely, exchanging the roles of $b$
and $b_\beta$ and setting $\epsilon_N = (-1)^{\frac{N-1}{2}}$, we have:
$$\begin{array}{l} \Rr_N(\Delta,b_{(01)},w',f',c')_{k,j}^{i,l}
\equiv_N \epsilon_N^{c_0}\ T_{k,k'}
\Rr_N(\Delta,b,w,f,c)_{i',j}^{k',l} \ (T^{-1})^{i',i}\\ \\
\Rr_N(\Delta,b_{(12)},w',f',c')_{j,l}^{i,k} \equiv_N
\epsilon_N^{c_1}\  T_{l,l'}\Rr_N(\Delta,b,w,f,c)_{i',j}^{k,l'}\ (S^{-1})^{i',i}
\\ \\ \Rr_N(\Delta,b_{(23)},w',f',c')^{k,j}_{i,l} \equiv_N
\epsilon_N^{c_0}\  S_{l,l'}\Rr_N(\Delta,b,w,f,c)_{i,j'}^{k,l'}\ (S^{-1})^{j',j}.
\end{array}$$
\end{prop}
The main properties of $S$ and $T$ are described in the following Lemma.
\begin{lem} \label{projrep} We have:
\begin{enumerate}
\item $S^4 = {\rm I}_{N}$. 
\item $ (TS)^3=\phi_N S^2$ where $\phi_N \in \{\pm 1, \pm
  i\}$ is given by
$$\phi_N =  \left\lbrace \begin{array}{ll}   
\left(\frac{m+1}{N}\right) & {\rm if}\ N \equiv 1 \ {\rm mod}(4)\\ 
\left(\frac{m+1}{N}\right) i  & {\rm if}\ N \equiv 3 
\ {\rm mod}(4).\end{array} \right. $$
\item $ (S^{-1}T)^3= \phi_N{\rm I}_{N}$.
\end{enumerate} 
\end{lem} 
{\it Proof.}  (1) is immediate. We have $\textstyle ((TS)^2)_i^j = N^{-1} \zeta^{-(m+1)i^2-ij}\sum_{k=0}^{N-1} \zeta^{(m+1)(k-i-j)^2}$. 
For all non vanishing coprime integers $a$, $b$ with $b>0$ set
$$G(a,b) = \sum_{x\ {\rm mod}(b)} e^{2\pi \sqrt{-1}a x^2/b}.$$
The sum in $((TS)^2)_i^j$ has this form for $a=m+1$, $b=N$.  By
\cite[page 86-87]{L} we have
$$G(a,b) = \left(\frac{a}{b}\right) G(1,b),\qquad b \ {\rm odd}$$
where $G(1,b) = \sqrt{b}$ if $b \equiv 1$ mod$(4)$, and $G(1,b) = i\sqrt{b}$ if $b \equiv 3$ mod$(4)$. Hence $((TS)^2)_i^j = \phi_N N^{-1/2} \zeta^{-(m+1)i^2-ij}$, and (2) follows easily from this. As $TS^{-1}=S(S^{-1}T)S^{-1}$, it is enough to prove the last
statement for $TS^{-1}$. This follows: $(TS)^4=(TS)(TS)^3=\phi_N TS^3=\phi_N TS^{-1}$; $(\phi_N
TS^{-1})^3=((TS)^4)^3=\phi_N^4{\rm I}_N$; 
finally $(TS^{-1})^3=\phi_N {\rm I}_N$.  
\hfill$\Box$

\subsection{Proof of Proposition \ref{wbinvariance}.} By Lemma \ref{local/global} we can assume that $(T,\tilde b)\to (T,\tilde b')$ corresponds to an oriented $C$-move or a circuit move on the associated $\Nn$-graphs. A circuit move preserves the $\mz/3\mz$-color $r(e)$ of every edge $e$, and produces one extremal pit and one extremal source associated to the matrices $S$ and $S^{-1}$ at the endpoints of each edge of the circuit. Hence $\Hh_N(T,\tilde b,w,f,c)$ is preserved, as it is the contraction of a tensor network. Consider the generating oriented $C$-move associated to $\sigma = (0123) \in S(J_4)$, and the tensor $\hat \Rr_N(T,b_{\sigma},w',f',c')$ resulting from the formal conversion of the sequence 
$$ C(b) \lra C(b_{(23)}) \lra C(b_{(23)(12)})\lra C(b_{(23)(12)(01)}) = C(b_{\sigma})$$
as before Proposition \ref{3.3}. From Figure \ref{or-C-gen} we see that an edge $e'$ of $C(b_{\sigma})$ carries the matrix $\Qq_N=S\cdot T^{-1}$ if $r(e')=1+r(e)\in \Z/3\Z$ ($e$ being the arc of $C(b)$ associated to $e'$), and that $e'$ carries $\Qq_N^{-1} = T\cdot S^{-1}$ if $r(e')= -1+ r(e)\in \Z/3\Z$.  Otherwise it carries the identity matrix. By the relation $\Qq_N^3= \phi_N^{-1}{\rm I}_{N}$ proved in Lemma \ref{projrep}, the sum in $\Hh_N(T,\tilde b,w,f,c)$ is changed by a factor $\pm i$ when $r(e')=1+r(e)$ and $r(e)=2$, or when $r(e')=-1+r(e)$ and $r(e)=0$. The term $\phi_N^{-q(T,\tilde b)}$ of $a_N(T,\tilde b)$ compensates exactly such factors. \cvd

\section{Resolution of the sign anomaly}\label{sign}
In this Section we prove Theorem \ref{more_cusped} (3)
and the analogous result for QHFT partition functions:
\begin{prop}\label{nosign} Let $\Pp$ be a QHFT pattern or a pattern over a cusped manifold. If $N-1 \equiv 0 \ {\rm mod} (4)$, or $N-1 \equiv 2 \ {\rm mod} (4)$ and the bulk $c$-weight of $\Pp$ vanishes, then $\Hh_N(\Pp)$ is defined up to multiplication by $N$-th roots of unity.
\end{prop}

\Dim We combine several results described below in Sections \ref{sign-transit} and \ref{sign-wb}. By Corollary \ref{QHtnosign}, for any odd $N\geq 3$ the invariance of $\Hh_N(\Tt)$ under QH transits holds true up to multiplication by $N$-th roots of unity. If $N-1 \equiv 0 \ {\rm mod} (4)$, we have $\epsilon_N=1$, so $\Hh_N(\Tt)$ has no sign anomaly with respect to any change of weak branching (Proposition \ref{3.3}).  If $N-1 \equiv 0 \ {\rm mod} (4)$ and $h_c=0$, the same is true by Lemma \ref{sign-C} and \ref{sign-circuit}. Then the conclusion follows as explained before Proposition \ref{wbinvariance} \cvd
\medskip

Let $\Pp$ be a pattern for which $\Hh_N(\Pp)$ can be defined by using branched triangulations, as in \cite{Top, GT, AGT}. Then $a_N(T,b) = N^{-v}$ or $1$, so we can wonder if the sign anomalies occuring in Proposition \ref{3.3} disappear because of global compensations. This is eventually true.
\begin{prop}\label{no-sign-H} For any pattern $\Pp$ as above, $\Hh_N(\Pp)$ is defined up to multiplication by a $N$-th root of unity, with no assumption on the weights of $\Pp$. 
\end{prop}
\Dim By  a result of Costantino \cite{C}, any change of branching can be achieved by means of a finite sequence of branching transits. The conclusion follows from Corollary \ref{QHtnosign}.  \cvd 
\begin{remark}{\rm One can ask if a weakly branched version of Costantino's result holds true: given weakly branched triangulations $(T,\tilde b)$ and $(T,\tilde b')$ of $\hat V$ that differ only by the weak branchings, are they related by a sequence of oriented $C$-moves and ``local" $b$-transits ? The answer is negative. In fact, we know that the state sums $\Hh_N$ have no sign anomaly with respect to QH transits and oriented $C$-moves, but there are examples of weakly branched triangulations $(T,\tilde b)$ admitting circuit moves along circuits $\gamma$ such that $h_c(\gamma)\ne 0 \in \mz/2\mz$.}
\end{remark}

\subsection {Behaviour of the sign anomaly under QH transits}\label{sign-transit} Recall the formula  \eqref{refR2}. The matrix dilogarithms $\Rr_N$ satisfy Proposition \ref{3.3}, whereas the {\it basic} dilogarithms $\Ll_N$ satisfy the relations (\cite{GT}, Proposition 5.3):
\begin{equation}\label{basicrelations}
\begin{array}{l} \Ll_N(\Delta,b_{(01)},w',f',c')_{k,j}^{i,l}
  \equiv_N (w_0')^{\frac{1-N}{2}}\ T_{k,k'}  \Ll_N(\Delta,b,w,f,c)_{i',j}^{k',l}\ (T^{-1})^{i',i}
\\ \\
  \Ll_N(\Delta,b_{(12)},w',f',c')_{j,l}^{i,k} \equiv_N
  (w_1')^{\frac{N-1}{2}}\ T_{l,l'} \Ll_N(\Delta,b,w,f,c)_{i',j}^
{k,l'}\ (S^{-1})^{i',i} 
  \\ \\ \Ll_N(\Delta,b_{(23)},w',f',c')^{k,j}_{i,l} \equiv_N
  (w_0')^{\frac{1-N}{2}}\ S_{l,l'}
\Ll_N(\Delta,b,w,f,c)_{i,j'}^{k,l'}\ (S^{-1})^{j',j}.
\end{array}
\end{equation}
Define the transits of $N$-th root cross ratio moduli like in \eqref{eqQHt}, by replacing $W(e)$ with the total $N$-th root modulus $W'(e)$ at all edges $e$ of $T$ and $T'$. Clearly any QH transit induces a transit of $N$-th root cross ratio moduli.
\begin{teo}{\rm (\cite{GT}, Theorem 5.2)} \label{bmatdilSchaeffer} The basic matrix
dilogarithm identity supported by the {\it Schaeffer}
$2\leftrightarrow 3$ transit of $N$-th root cross ratio moduli holds true up to an $N$-th root of
unity anomaly.  

\end{teo}
The proof goes as follows. The basic matrix dilogarithm identity supported by the Schaeffer $2\leftrightarrow 3$ QH transit is equivalent to the following equality
in ${\rm Aut}(\mc^N \otimes \mc^N \otimes \mc^N)$:
\begin{equation} \label{FadKas}
\Psi_{23}^1(V)\Psi_{12}^3(U) =
\Psi_{12}^4(U)\Psi_{13}^2(-UV)\Psi_{23}^0(V)
\end{equation}
where $U$ and $V$ are explicit $N^2\times N^2$ matrices
satisfying $U^N=-Id \otimes Id$, $V^N=-Id \otimes Id$, and $VU=\zeta
UV$, $\zeta$ being a primitive $N$-th root of $1$, and for a tensor $A$ satisfying $A^N=-Id$ (on the appropriate
space: $\mc^N$, $\mc^N \otimes \mc^N$, etc.) one sets
\begin{equation}\label{formmat2b}
\Psi^i(A)= \sum_{l=0}^{N-1}A^l\ 
\prod_{s=1}^{l} \frac{((w_0^i)')^{-1}((w_1^i)')^{-1}}{1-((w_0^i)')^{-1}\zeta^{-s}}.
\end{equation}
In \eqref{FadKas}, $\Psi^1_{23}(V)$ means $\Psi^1(V)$ acting on the second and third tensorands, i.e. ${\rm Id} \otimes \Psi^1(V)$, etc. One can check that $\Psi^i$ is uniquely determined up to
multiplication by scalars by the functional relation
\begin{equation}\label{factor2}
\Psi^i(\zeta^{-1}A) =\Psi^i(A)\ \left( \frac{1-((w_0^i)')^{-1}
((w_1^i)')^{-1}A}{((w_0^i)')^{-1}}\ \right) = \Psi^i(A)\ \left((w_0^i)' -
((w_1^i)')^{-1}A \right).
\end{equation}  
The identity \eqref{FadKas} lifts to the algebra generated by $U$ and $V$, so we have to show
\begin{equation} \label{FadKas2}
\Psi^1(V)\Psi^3(U) = \Psi^4(U)\Psi^2(-UV)\Psi^0(V)
\end{equation}
when $U^N=-1$, $V^N=-1$, and $VU=\zeta UV$. To do it, 
one uses crucially the relations between $N$-th root cross-ratio moduli implied by the Schaeffer QH transit. One shows that $UV$ commutes with $P(-UV):=\Psi^4(U)^{-1}\Psi^1(V)\Psi^3(U)\Psi^0(V)^{-1}$; hence $P(-UV)$ is a function of $UV$. Moreover, it satisfies
\begin{equation} \label{eqnouv}
P(-\zeta^{-1}UV) = P(-UV) \ \bigl((w_0')^2 - ((w_1')^2)^{-1}\  (-UV)\bigr).
\end{equation}
Comparing with \eqref{factor2} we see that $P(-UV)$ and $\Psi^2(-UV)$
are equal up to a multiplicative constant $C:= C((w_0')^2,(w_1')^2)$. We introduce an $N$-th root of the reciprocal of $det(\Psi^i(A))$ in each entry of $\Psi^i(A)$ so that $det(\Psi^i(A)) = 1$ (this is the term $h((w_0^i)') = g((w_0^i)')/g(1)$ in \eqref{basicf}). Then $det(P(-UV)) = 1 = C^N det(\Psi^2(-UV))$ and $C^N =1$. As a result the tensor $\Ll_N(w_0',w_1')$ can be expressed in terms of $\Psi(A)$ for a suitable $A$ (see \cite{GT}, formulas (32) and (33)).
\cvd
\medskip

The basic matrix dilogarithm identities hold true up to an
$N$-th root of unity anomaly only for some $2\leftrightarrow 3$ branching transits. From \eqref{basicrelations} one deduces the following table of the anomalies produced by the symmetries of the Schaeffer $2\leftrightarrow 3$ transit (it is enough to consider those induced by the transpositions $(01), \ldots, (34)$ on the vertex ordering used in Figure \ref{5-terms}). The transits of $N$-th root moduli imply that the anomalies are equal for the symmetries induced by $(01)$ and $(34)$, but differ for $(12)$ and $(23)$ because of the occurence of the opposite exponents $(1-N)/2$, $(N-1)/2$ for the $N$-th root moduli
$(w_0^k)'$, $(w_1^k)'$ respectively. For instance $(w_1^0)'(w_0^4)' = (w_1^1)'$, but we see $((w_1^0)')^{\frac{N-1}{2}}((w_0^4)')^{\frac{1-N}{2}} \ne ((w_1^1)')^{\frac{N-1}{2}}$.

$$\begin{tabular}{||c||c|c||c|c|c||} \hline &
  \rule{0cm}{0.4cm}{$\Delta^1$}& $\Delta^3$ & $\Delta^0$ & $\Delta^2$
  & $\Delta^4$ \\ \hline

\rule{0cm}{0.8cm} \raisebox{0.2cm}{$(01)$} & \raisebox{0.2cm}{$1$} &
\raisebox{0.2cm}{$((w_0^3)')^{\frac{1-N}{2}}$} & \raisebox{0.2cm}{$1$} &
\raisebox{0.2cm}{$((w_0^2)')^{\frac{1-N}{2}}$} & 
\raisebox{0.2cm}{$((w_0^4)')^{\frac{1-N}{2}}$}\\ \hline

\rule{0cm}{0.8cm} \raisebox{0.2cm}{$(12)$} & \raisebox{0.2cm}{$1$} &
\raisebox{0.2cm}{$((w_1^3)')^{\frac{N-1}{2}}$} & 
\raisebox{0.2cm}{$((w_0^0)')^{\frac{1-N}{2}}$} &
\raisebox{0.2cm}{$1$} & \raisebox{0.2cm}{$((w_1^4)')^{\frac{N-1}{2}}$} 
\\ \hline

\rule{0cm}{0.8cm} \raisebox{0.2cm}{$(23)$} &
\raisebox{0.2cm}{$((w_1^1)')^{\frac{N-1}{2}}$} & \raisebox{0.2cm}{$1$} &
\raisebox{0.2cm}{$((w_1^0)')^{\frac{N-1}{2}}$} & \raisebox{0.2cm}{$1$} &
\raisebox{0.2cm}{$((w_0^4)')^{\frac{1-N}{2}}$} \\ \hline

\rule{0cm}{0.8cm} \raisebox{0.2cm}{$(34)$} &
\raisebox{0.2cm}{$((w_0^1)')^{\frac{1-N}{2}}$} & \raisebox{0.2cm}{$1$} &
\raisebox{0.2cm}{$((w_0^0)')^{\frac{1-N}{2}}$} 
& \raisebox{0.2cm}{$((w_0^2)')^{\frac{1-N}{2}}$} &
\raisebox{0.2cm}{$1$} \\ \hline
\end{tabular}$$

%\begin{remark}\label{no-local-basic}{\rm The above discussion justifies the comment {\bf c)} after Theorem \ref{more_cusped}: {\it there is no invariance proof of the reduced QH state sums based on the full set of QH transits. }}\end{remark}
\medskip

Replacing the tensors $\Ll_N$ with $\Rr_N$ we have:

\begin{prop}\label{transit-no-sign} The matrix dilogarithm identities hold true up to an $N$-th root of unity anomaly for all QH transits.
\end{prop}
\noindent {\it Proof.} Since we have the identity \eqref{FadKas}, in the case of the Schaeffer QH transit identity it is enough to compare the
symmetrization factors at both sides. They are ($\textstyle \epsilon_N := (-1)^{\frac{N-1}{2}}$):
$$\epsilon_N^{f_0^1c_1^1+f_1^1c_0^1+f_0^3c_1^3+f_1^3c_0^3}\ 
\exp\biggl(\frac{N-1}{2N}\bigl(-c_1^1{\rm l}_0^1+ c_0^1{\rm l}_1^1-c_1^3
{\rm l}_0^3 +c_0^3{\rm l}_1^3\bigr)\biggr)$$
and 
$$\epsilon_N^{f_0^0c_1^0+f_1^0c_0^0+f_0^2c_1^2+
  f_1^2c_0^2+f_0^4c_1^4+f_1^4c_0^4}\
\exp\biggl(\frac{N-1}{2N}\bigl(-c_1^0{\rm l}_0^0+c_0^0{\rm
  l}_1^0-c_1^2{\rm l}_0^2+c_0^2{\rm l}_1^2-c_1^4{\rm l}_0^4+c_0^4{\rm
  l}_1^4\bigr) \biggr).$$ A computation shows that the transits of log-branches and charges
imply that the two exponentials are equal, and with the flattening transit
mod$(2)$ that the signs are equal. To prove the claim
for all others $2\leftrightarrow 3$ QH transits we compare the anomalies resulting from symmetries on the Schaeffer QH transit. Since all QH transits are branched, the actions by the matrices $S$, $T$ due to Proposition \ref{3.3} cancel out along interior faces. The following table shows the sign anomalies. Again, the transit of charges implies that they are equal at both sides.  

$$\begin{tabular}{||c||c|c||c|c|c||} \hline &
  \rule{0cm}{0.4cm}{$\Delta^1$}& $\Delta^3$ & $\Delta^0$ & $\Delta^2$
  & $\Delta^4$ \\ \hline

\rule{0cm}{0.8cm} \raisebox{0.2cm}{$(01)$} & \raisebox{0.2cm}{$1$} &
\raisebox{0.2cm}{$\epsilon_N^{c_0^3}$} & \raisebox{0.2cm}{$1$} &
\raisebox{0.2cm}{$\epsilon_N^{c_0^2}$} & \raisebox{0.2cm}{$\epsilon_N^{c_0^4}$}\\ 
\hline

\rule{0cm}{0.8cm} \raisebox{0.2cm}{$(12)$} & \raisebox{0.2cm}{$1$} &
\raisebox{0.2cm}{$\epsilon_N^{c_1^3}$} & \raisebox{0.2cm}{$\epsilon_N^{c_0^0}$} &
\raisebox{0.2cm}{$1$} & \raisebox{0.2cm}{$\epsilon_N^{c_1^4}$} \\ \hline

\rule{0cm}{0.8cm} \raisebox{0.2cm}{$(23)$} &
\raisebox{0.2cm}{$\epsilon_N^{c_1^1}$} & \raisebox{0.2cm}{$1$} &
\raisebox{0.2cm}{$\epsilon_N^{c_1^0}$} & \raisebox{0.2cm}{$1$} &
\raisebox{0.2cm}{$\epsilon_N^{c_0^4}$} \\ \hline

\rule{0cm}{0.8cm} \raisebox{0.2cm}{$(34)$} &
\raisebox{0.2cm}{$\epsilon_N^{c_0^1}$} & \raisebox{0.2cm}{$1$} &
\raisebox{0.2cm}{$\epsilon_N^{c_0^0}$} & \raisebox{0.2cm}{$\epsilon_N^{c_0^2}$} &
\raisebox{0.2cm}{$1$} \\ \hline
\end{tabular}$$
\medskip  

As usual, the identities for the bubble moves are formal consequence
of those for the $2\leftrightarrow 3$ moves. This concludes the proof. \cvd

\begin{cor} \label{QHtnosign} Let $(T,\tilde b,w,f,c)\to (T',\tilde b', w',f',c')$ be any
QH transit of weakly branched QH triangulations of a pattern $\Pp$. Then
$$\Hh_N(T,\tilde b, w,f,c) \equiv_N  \Hh_N(T,\tilde b',w',f',c').$$
\end{cor}
\Dim If the transit is supported by a MP move it is immediate that $a_N(T,\tilde b)=a_N(T',\tilde b')$, hence the conclusion follows from Corollary \ref{transit-no-sign}. If it is supported by a positive bubble move, spelling the contribution of each term we realize that if 
$$ a_N(T,\tilde b)= N^{-v}\phi_N^{-q(T,\tilde b)}c_N(T,\tilde b)\ , \  c_N(T,\tilde b)=\epsilon_N^{v+l-\frac{1}{2}\sum_e (n_+(e)-n_-(e))}$$ then
$$ a_N(T',\tilde b')= N^{-(v+1)}\phi_N^{-q(T',\tilde b')}c_N(T',\tilde b')\ ,\ c_N(T',\tilde b')= \epsilon_N^{(v+1)+(l+3)-(2+\frac{1}{2}\sum_e (n_+(e)-n_-(e))}$$
and $\phi_N^{-q(T',\tilde b')}=\phi_N^{-q(T,\tilde b)}$. Hence $ a_N(T',\tilde b')= N^{-1} a_N(T,\tilde b)$, as it must be in order
to ensure the invariance under the bubble move, which increases by $1$ the
number of internal vertices.  \cvd   

\subsection{Behaviour of the sign anomaly under a change of weak branching}
\label{sign-wb} A change of weak branching $\tilde b\to \tilde b'$ modifies $a_N(T,\tilde b)$ by its sub-factors $\phi_N^{-q(T,\tilde b)}$ and
$$c^*_N(T,\tilde b):=  \epsilon_N^{l-\frac{1}{2}\sum_e (n_+(e)-n_-(e))}.$$
By the proof of Proposition \ref{wbinvariance},  $\phi_N^{-q(T,\tilde b)}$ and the unnormalized QH state sum have reciprocal variations. So it remains to analyze the behaviour of $c^*_N$ with respect to oriented $C$-moves and the circuit move of Lemma \ref{local/global}. 

It is relevant here to clarify the nature of $c^*_N$.  By extending Chapter 7 of \cite{BP0}, the paper \cite{BP2} provides a combinatorial realization of spin structures on
$3$-manifolds, based on weakly branched triangulations instead of
branched ones. In particular, one associates to every weakly branched
triangulation $(T,\tilde b)$ a framing $\nu_{\tilde b}$ of the manifold defined along
${\rm Sing}(P)$, and computes a cellular cochain 
$\alpha_{\tilde b} \in C^2(P; \Z/2\Z)$ representing the {\it obstruction}
to extend $\nu_{\tilde b}$ over the whole of the spine $P$. Consider the
cellular $2$-chain $\textstyle R(P)=\sum_R R \in C_2(P;\Z/2\Z)$, where $R$ varies
among the $2$-regions of $P$. It is a fact that 
$$l-\frac{1}{2}\sum_e (n_+(e)-n_-(e)) = 
\alpha_{\tilde b}(R(P)) \  {\rm mod} (2). $$
How $\alpha_{\tilde b}$, whence $\alpha_{\tilde b}(R(P))$,  varies with the
weak branching is part of the theory developed in \cite{BP2}. As a consequence 
we get:
\begin{lem}\label{sign-C} For any change of weak branching $\tilde b \to \tilde b'$ preserving the induced pre-branching we have
 $$\Hh_N(T,\tilde b,w,f,c)\equiv_{N} \Hh_N(T,\tilde b',w,f,c).$$
\end{lem}  
\Dim It is enough to show the result for one generating oriented $C$-move $\sigma = (0123)$. By local inspection one checks that $\alpha_{\tilde b'}(R(P))= \alpha_{\tilde b}(R(P))+1$, so $c^*_N(T,\tilde b')=\epsilon_N c^*_N(T,\tilde b)$. On the other hand, using Proposition \ref{3.3}, the factorization $(0123) = (23)\cdot (12)\cdot (01)$, and  the relation $c_0+c_1+c_2=1$, we see that $\sigma$ changes the unnormalized QH state sum by $\epsilon_N^{c_1} \epsilon_N^{c_2}\epsilon_N^{c_0} = \epsilon_N$. Hence the two sign variations are equal.\cvd
  
\begin{lem} \label{sign-circuit} For any change of weak branching $\tilde b \to \tilde b'$ induced by one
circuit move along a circuit $\gamma$ we have
$$\Hh_N(T,\tilde b,w,f,c)\equiv_N  \epsilon_N^{h_c(\gamma)} \Hh_N(T,\tilde b',w,f,c)$$
where $h_c$ is the bulk $c$-weight of $\Pp$. In particular, $\Hh_N(\Pp)$ has no sign
anomaly if $h_c=0$.
\end{lem} 
\Dim Let us realize the circuit move in terms of $\Nn$-graphs. By the third identity
of Proposition \ref{3.3} the move changes the state sum by one factor $\epsilon_N^{c_0}$ for each crossing of $\gamma$. It is easy to check that their product computes $h_c(\gamma)$. On the other hand $\alpha_{\tilde b'}(R(P))= \alpha_{\tilde b}(R(P))$, so that $c^*_N(T,\tilde b')= c^*_N(T,\tilde b)$.  \cvd

\section{Examples}\label{EXAMPLES} % Le texte de toute la section a ŽtŽ remaniŽ - pas les \'equations
In this section we describe the quantum hyperbolic invariants of the figure-eight knot complement, its ``sister", and the complement of the knot $5_2$. 

We give also some samples of numerical computations obtained by using Maple, which are part of our current analytical and numerical exploration of the QHI. At first we quote that they corroborate the invariance of the QH state sums with respect to the choice of flattenings or charges having the same weights. In the case of the sister manifold, they confirm also the invariance of the QH state sums supported by weakly branched triangulations. Moreover, the numerical computations that we present below are proving: 
\begin{itemize}
 \item The actual dependence of the invariants on the weights (and not only on the
characters $\rho$).
\item This dependence persists as $N\rightarrow \infty$. 
\item For some weights the invariants seem to grow exponentially with $N$, and yield instances of the volume conjecture, and for some other weights they do not.
\end{itemize}
Other interesting phenomena happen, that deserve to be understood. For instance, in the case of the complements of the knots $4_1$ and $5_2$, it appears that the asymptotical behaviour of the absolute values of the invariants do not depend on the values on the meridian of the charge weights or the flattening weights. 

\subsection{Notations} \label{NOT} Let $(\Delta,b,w,f,c)$ be a quantum $3$-simplex. Put $d_k := f_k -*_jc_k$, $k\in \{0,1,2\}$. We call $d_k$ an {\it edge color}. In order to simplify the formulas we denote the $N$-th roots of the cross-ratios $w_k$ by bold letters $\bw_k$ (rather than $w_k'$, as we did in \eqref{N_root}), and call them {\it q-shape parameters}. So  
$$\bw_k := \exp\left(\frac{1}{N}(\log(w_k) +\pi i (N+1)d_k)\right).$$
For any cusped manifold $M$ and any QH triangulation $\Tt$ of $M$, we have the tetrahedral and edge relations (see \eqref{relloc}, \eqref{edgeconst} and \eqref{edge-equation})
$$\sum_{k=0,1,2} (\log(w_k) +\pi i d_k) = -*_b\pi i\ ,\ \sum_{E\ra e} *_E(\log(w(E)) +\pi i d(E)) =  -2\pi i$$ 
and 
$$\prod_{k=0,1,2} \bw_k = e^{-*_b\frac{\pi i}{N}}\ ,\ \prod_{E\ra e} \bw(E)^{*_E} =  e^{-\frac{2\pi i}{N}}$$
where, as usual, $d(E)=d_k$ and $\bw(E)=\bw_k$ if $E$ is the edge $E_k$ or the opposite edge of a quantum $3$-simplex $(\Delta,b,w,f,c)$ of $\Tt$, and $*_E:=*_{b}$. We denote by $G_N(T,\tilde b)$ the variety thus defined by the q-shape parameters, and call it the {\it QH gluing variety} (NB: it coincides with the variety $G(T,\tilde b,c)_N$ of Section \ref{glob}; to simplify notations, here we drop the reference to the charge $c$).

As usual, denote by $V$ a compact $3$-manifold with one torus boundary component such that $M$ is diffeomorphic to the interior of $V$. We compute the invariant $\Hh_{N}(V, \rho,k_f, k_c,h_f,h_c)$ by means of the function $\Hh_N'(T,\tilde b,c)$ on $G_N(T,\tilde b)$ introduced after Definition \ref{defiSS}. As $\Hh_N'(T,\tilde b,c)$ is formulated in terms of q-shape parameters and charges, we express equivalently $\Hh_{N}(V, \rho,k_f, k_c,h_f,h_c)$ as $\Hh_{N}(V, \rho,\kappa, k_c,h,h_c)$, using the classes $\kappa$ and $h$ obtained from the flattening weights $k_f$ and $h_f$ by replacing the flattenings $f_k$ by the edge colors $d_k$ in the respective formulas (see Definition \ref{defiweight}, and \eqref{weightwf} for $\kappa$). Clearly $\kappa = k_f - \pi i k_c$ and $h=h_c+h_f$. We call $\kappa$ and $h$ the {\it reduced} boundary and bulk weights. By \eqref{cconstraint} and \eqref{fconstraint2} one has the compatibility relations
\begin{equation}\label{dconstraint} \kappa(a)= d_w(a) \ {\rm mod}(i\pi\mz)\ ,\ (\kappa(a)-d_w(a))/i\pi = i^*(h)(a) \ \ {\rm mod} (2\mz) \end{equation}
for all $a\in H_1(\partial V; \Z)$. 
\smallskip

We consider triangulations with ordered tetrahedra. We denote by $u,v,w\ldots$ (resp. $\bu,\bv,\bw\ldots$, resp. $a,b,c\ldots$, resp. $a',b',c'\ldots$) the triples of cross-ratios (resp. q-shape parameters, resp. edge colors, resp. charges) of the tetrahedra, in the same order. 

\subsection{The figure-eight knot complement}\label{4_1} The Epstein-Penner decomposition of $M$ is an ideal triangulation $T$ having a branching $b$ and two tetrahedra, with opposite $b$-orientations. We showed the dual spine in Figure \ref{I-F8}. We denote by $\Delta^0$ the tetrahedron with positive $b$-orientation. In \cite{AGT} we described the flattenings at positive points of the gluing variety $G(T,b)$, as well as the charges and the QH state sums carried by $(T,b)$. %We are going to use these datas. 

\subsubsection{The QH gluing variety} The gluing variety $G(T,b)$ is the irreducible plane curve with coordinates $(u_0,v_0)$ and defining equation $u_1u_2^2v_0^{-2}v_1^{-1}=1$, which may be written as $u_0^2v_0^2=(1-u_0)(1-v_0)$ (as usual $u_{i+1} = 1/(1-u_i)$ and similarly for the $v_i$'s). The point $(e^{\frac{\pi i}{3}}, e^{-\frac{\pi i}{3}})$ realizes the complete hyperbolic structure. Let $l$ and $m$ be the canonical longitude and the meridian of the figure-eight knot, and let $(u_0,v_0)\in G(T,b)$. The dilation factors of the holonomy $\rho := \rho(u_0,v_0)$ are given by
$${\rm hol}_l(\rho) = u_0^2u_2^{-2} = \frac{u_0^4}{(1-u_0)^2}\ ,\ {\rm hol}_m(\rho) = u_2v_2 = \frac{(u_0-1)(v_0-1)}{u_0v_0}.$$
The tetrahedral and edge relations between q-shape parameters are
$$\bu_0\bu_1\bu_2 =  e^{-\frac{\pi i}{N}}\ ,\ \bv_0\bv_1\bv_2 =  e^{+\frac{\pi i}{N}}$$
$$\bu_1\bu_2^2\bv_0^{-2}\bv_1^{-1} = e^{-\frac{2\pi i}{N}}\ ,\  \bu_1\bu_0^2\bv_2^{-2}\bv_1^{-1}= e^{-\frac{2\pi i}{N}}.$$
Together with the identities $(\bu_{i+1})^N = 1/(1-(\bu_i)^N)$ and $(\bv_{i+1})^N = 1/(1-(\bv_i)^N)$, these relations define $G_N(T,b)$. Using the tetrahedral relations it is easy to see that the two edge relations are equivalent. Recall the reduced weights $h$ and $\kappa$ introduced in Section \ref{NOT}. By \eqref{dconstraint} we have 
\begin{equation}\label{compbord}
\bu_0^{2}\bu_2^{-2}  = e^{\frac{\kappa(l)}{N}+\pi i h(l)}\ ,\ \bu_2\bv_2= e^{\frac{\kappa(m)}{N}+\pi i h(m)}.
\end{equation}
As $M$ is a knot complement, the bulk weight $h$ is determined by $\kappa$, and similarly $h_c$ is determined by $k_c$. Moreover, $h(l)=0$, so $h$ is non zero if and only if $h(m) = a_2+b_2 = 1$ mod$(2)$. 

Over the point $(e^{\frac{\pi i}{3}}, e^{-\frac{\pi i}{3}})$ we have $\kappa(m)\in \pi i \mz$ and $\kappa(l)\in 2\pi i \mz$, and the tetrahedral relations, the edge relations, and the compatibility relations with $\kappa$ imply the following expressions of the edge colors $a_k$ and $b_k$ in terms of $a_0$ and $\kappa$: $a_2=-2-a_0-a_1$, $b_2=2-b_0-b_1$,  and
\begin{equation}\label{standardf3} a_1=-2a_0+\frac{\kappa(l)}{2\pi i}-2\ ,\ b_0 = \frac{\kappa(m)}{\pi i}-a_0\ ,\ b_1 = 2a_0-\frac{2\kappa(m)}{\pi i}-\frac{\kappa(l)}{2\pi i}+2.
\end{equation}
For points $(\bu,\bv)\in G_N(T,b)$ lying over a sufficiently small neighborhood of the point $(e^{\frac{\pi i}{3}}, e^{-\frac{\pi i}{3}}) \in G(T,b)$, the same expressions are valid if one replaces $\kappa(\cdot)/\pi i$ with $(\kappa(\cdot )-\log({\rm hol}_{\cdot}(\rho)))/\pi i$. 

Similarly the charges $a'_k$ and $b'_k$ are given by
\begin{equation}\label{standardf3} a'_1=-2a_0'+\frac{k_c(l)}{2}+1\ ,\ b'_0 = a'_0-k_c(m)\ ,\ b'_1 = -2a_0'+2k_c(m)+\frac{k_c(l)}{2}+1.
\end{equation}
Let $M'$ be a closed hyperbolic $(p,q)$-Dehn filling of $M$ represented by $(u_0,v_0) \in G(T,b)$ in a small neighborhood of $(e^{\frac{\pi i}{3}}, e^{-\frac{\pi i}{3}})$. The holonomy $\rho$ factors through the quotient map $\pi_1(M) \rightarrow \pi_1(M')$ to define the holonomy $\rho'$ of $M'$. One has
\begin{equation}\label{chdfeq}
p\log({\rm hol}_m(\rho') ) +q\log({\rm hol}_l(\rho') ) = 2\pi i,
\end{equation}
and letting $r$, $s \in \mz$ be such that $ps-qr=1$, the reduced weight $\kappa$ satisfies $\kappa(m^pl^q)= 0$ if and only if
\begin{equation}\label{surgf3}
a_1 = r-2-2a_0\ ,\ b_0=-2s-a_0\ ,\ b_1 = 2-r+4s+2a_0.
\end{equation}
Similar formulas express the identity $k_c(m^pl^q)= 0$.

\subsubsection{The state sum formulas} Recall the notations of Section \ref{tettens}.  For any $(\bu,\bv)\in G_N(T,b)$ and any charge $c$ on $T$, the QH state sum 
\begin{equation}\label{form1}
\Hh_{N}'(T,b,c)(\bu,\bv) = (\bu_0^{-a_1'}\bu_1^{a_0'}\bv_0^{-b_1'}\bv_1^{b_0'})^{\frac{N-1}{2}}
\frac{[\bv_0]g(\bu_0)}{g(\bv_0)}\sum_{\alpha,\beta=0}^{N-1}
    \zeta^{\beta^2-\alpha^2} \frac{\omega(\bu_0,\bu_1^{-1}\vert N-\beta)}{
    \omega(\bv_0/\zeta,\bv_1^{-1}\vert N-\alpha)}
\end{equation}
computes the invariant $\Hh_{N}(V, \rho,\kappa, k_c,h,h_c)$, and defines a (non-constant) rational function on $G_N(T,b)$. By a surgery theorem proved in \cite{AGT}, at points $(\bu,\bv)$ where $\kappa$ and $k_c$ satisfy $\kappa(m^pl^q)= k_c(m^pl^q)=0$ as above, $\Hh_{N}(V, \rho,\kappa, k_c,h,h_c)$ coincides with the $N$-th quantum hyperbolic invariant of the closed manifold $(M',L,\rho')$, where $L$ is the surgery core.

\subsubsection{Some numerical results}\label{num1} Denote by $\rho(A,B):=\rho(u_0(A,B),v_0(A,B))$ the holonomy of $M$ given by the shape parameter $u_0(A,B):=e^{\frac{\pi i}{3}}+A+iB$ and the solution $v_0(A,B)$ of the gluing equation $u_0(A,B)^2v_0(A,B)^2-(1-u_0(A,B))(1-v_0(A,B))=0$ such that $u_0(A,B)$ has positive imaginary part and $v_0(A,B)$ has negative imaginary part. So $\rho(0,0)$ is the hyperbolic holonomy. 

Recall that $h$ and $h_c$ are determined by $\kappa$ and $k_c$ respectively. In the following table we present some values of $\vert \Hh_{15}(V, \rho(A,B),\kappa,k_c, h,h_c) \vert$ up to 4 digits, where $(A,B)$ and $k_c(l)$ are as indicated on the top and the left of the corresponding column and row ($(A,B)_{(2,3)}:=(0.06734378...
..., -.42400885...)$ defines the holonomy of the $(2,3)$-Dehn filling of $M$), and we put $\kappa(m)=k_c(m)=0$ and $\kappa(l)=k_f(l) - \pi i k_c(l)$, where $k_f(l)=L([w;f])(l)$ has the collection of shape parameters $w:=w(A,B)$ determined by $(u_0(A,B),v_0(A,B))$, and a (compatible) flattening $f$ such that $L([w(0,0);f])=-2\pi i$ at the complete hyperbolic structure. 

$$\begin{tabular}{||c||c|c|c|c|c||} \hline &
  \rule{0cm}{0.4cm}{$(0,0)$}& $(1,0)$ & $(1,0.5)$ & $(1,1)$ & $(A,B)_{(2,3)}$ \\ \hline

\rule{0cm}{0.8cm} \raisebox{0.2cm}{$k_c(l)=6$} & \raisebox{0.2cm}{$2.5587...$} &
\raisebox{0.2cm}{$2.6504....$} & \raisebox{0.2cm}{$2.0018....$} &
\raisebox{0.2cm}{$1.6118....$} & \raisebox{0.2cm}{$5.0307....$}\\ 
\hline

\rule{0cm}{0.8cm} \raisebox{0.2cm}{$k_c(l)=4$} & \raisebox{0.2cm}{$58.5466...$} &
\raisebox{0.2cm}{$58.3761...$} & \raisebox{0.2cm}{$47.0533...$} &
\raisebox{0.2cm}{$39.9892...$} & \raisebox{0.2cm}{$95.0326...$}\\ 
\hline

\rule{0cm}{0.8cm} \raisebox{0.2cm}{$k_c(l)=2$} & \raisebox{0.2cm}{$2.1356...$} &
\raisebox{0.2cm}{$2.0279...$} & \raisebox{0.2cm}{$1.9058...$} &
\raisebox{0.2cm}{$1.8138...$} & \raisebox{0.2cm}{$2.7491...$}\\ 
\hline

\rule{0cm}{0.8cm} \raisebox{0.2cm}{$k_c(l)=0$} & \raisebox{0.2cm}{$77.4851...$} &
\raisebox{0.2cm}{$77.5401...$} & \raisebox{0.2cm}{$77.5885...$} &
\raisebox{0.2cm}{$77.6997...$} & \raisebox{0.2cm}{$76.4850...$}\\ 
\hline

\rule{0cm}{0.8cm} \raisebox{0.2cm}{$k_c(l)=-2$} & \raisebox{0.2cm}{$0.1118...$} &
\raisebox{0.2cm}{$0.1620...$} & \raisebox{0.2cm}{$0.1672...$} &
\raisebox{0.2cm}{$0.1746...$} & \raisebox{0.2cm}{$0.0650...$}\\ 
\hline

\rule{0cm}{0.8cm} \raisebox{0.2cm}{$k_c(l)=-4$} & \raisebox{0.2cm}{$77.4851...$} &
\raisebox{0.2cm}{$77.5401...$} & \raisebox{0.2cm}{$95.8738...$} &
\raisebox{0.2cm}{$112.7032...$} & \raisebox{0.2cm}{$47.4247...$} \\ \hline

\rule{0cm}{0.8cm} \raisebox{0.2cm}{$k_c(l)=-6$} &
\raisebox{0.2cm}{$2.1356...$} & \raisebox{0.2cm}{$2.0279...$} &
\raisebox{0.2cm}{$2.7647...$} & \raisebox{0.2cm}{$3.5183...$} &
\raisebox{0.2cm}{$1.0549...$} \\ \hline

\rule{0cm}{0.8cm} \raisebox{0.2cm}{$k_c(l)=-8$} &
\raisebox{0.2cm}{$58.5466...$} & \raisebox{0.2cm}{$58.3761...$} &
\raisebox{0.2cm}{$89.7752...$} & \raisebox{0.2cm}{$124.5239...$} &
\raisebox{0.2cm}{$21.8881...$} \\ \hline

\rule{0cm}{0.8cm} \raisebox{0.2cm}{$k_c(l)=-10$} &
\raisebox{0.2cm}{$2.5587...$} & \raisebox{0.2cm}{$2.6504...$} &
\raisebox{0.2cm}{$4.2536...$} & \raisebox{0.2cm}{$6.1199...$} &
\raisebox{0.2cm}{$0.7577...$} \\ \hline

\rule{0cm}{0.8cm} \raisebox{0.2cm}{$k_c(l)=-12$} &
\raisebox{0.2cm}{$33.2019...$} & \raisebox{0.2cm}{$33.1224...$} &
\raisebox{0.2cm}{$63.6250...$} & \raisebox{0.2cm}{$104.4917...$} &
\raisebox{0.2cm}{$7.4274...$} \\ \hline
\end{tabular}$$
\medskip

Note that the table shows that the invariants depend on the character and the weights. Also, the first two columns show a symmetry about $k_c(l)=-2$ (that is, when $(A,B)=(0,0)$, about $\kappa(l)=0$), corresponding to a change of orientation of $l$. The dominant rows are for $k_c(l)=0$ or $-4$ (that is, when $(A,B)=(0,0)$, about $\kappa(l)=-2\pi i$ or $2\pi i$).
\smallskip

Now consider the behaviour of $\vert \Hh_{N}(V, \rho(A,B),\kappa,k_c, h,h_c) \vert$ as $N\rightarrow +\infty$. Put $G_N(\kappa(l)):=\pi\log\left(\vert \Hh_{N+2}(V,\rho(0,0),\kappa,k_c,h,h_c) \vert / \vert \Hh_{N}(V,\rho(0,0),\kappa,k_c,h,h_c) \vert\right)$, where the arguments are as above except that we fix $(A,B)=(0,0)$ (the complete hyperbolic structure). Let $\Hh_{N}''(T,b,c)$ be $\Hh_{N}'(T,b,c)$ divided by $(\bu_0^{-a_1'}\bu_1^{a_0'}\bv_0^{-b_1'}\bv_1^{b_0'})^{\frac{N-1}{2}}$, the product of the matrix dilogarithm symmetrization factors. Note that $G_N:=\pi\log\left(\vert \Hh_{N+2}''(T,b,c)(\bu,\bv)) \vert / \vert \Hh_{N}''(T,b,c)(\bu,\bv) \vert\right)$ and $G_N(\kappa(l))$ are equivalent as $N\rightarrow +\infty$. The following table gives a sample of values of $G_{151}$, where $\kappa(l)$ takes the values indicated in the first row. 
\medskip

%$$\begin{tabular}{||c||c|c|c|c|c|c|c||} \hline &
%  \rule{0cm}{0.4cm}{$-6\pi i$}& $-4\pi i$ & $-2\pi i$ & $0$ & $2\pi i$ & $4\pi i$ & $6\pi i$ \\ \hline

%\rule{0cm}{0.8cm} \raisebox{0.2cm}{$G_N'(\kappa(l))$} & \raisebox{0.2cm}{$2.03069972...$} &
%\raisebox{0.2cm}{$-0.49036336...$} & \raisebox{0.2cm}{$2.02968013...$} &
%\raisebox{0.2cm}{$0.48922984...$} & \raisebox{0.2cm}{$2.02968013...$}& \raisebox{0.2cm}{$-0.49036336...$}& \raisebox{0.2cm}{$2.03069972...$}\\ \hline
%\end{tabular}$$
%\medskip

$$\begin{tabular}{||c||c|c|c|c|c|c|c||} \hline &
  \rule{0cm}{0.4cm}{$-6\pi i$}& $-4\pi i$ & $-2\pi i$ & $0$ & $2\pi i$ %& $4\pi i$ & $6\pi i$ 
  \\ \hline

\rule{0cm}{0.8cm} \raisebox{0.2cm}{$G_{151}$} & \raisebox{0.2cm}{$2.03069...$} &
\raisebox{0.2cm}{$-0.49036...$} & \raisebox{0.2cm}{$2.02968...$} &
\raisebox{0.2cm}{$0.48922...$} & \raisebox{0.2cm}{$2.02968...$}%& \raisebox{0.2cm}{$-0.49036...$}& \raisebox{0.2cm}{$2.03069...$}
\\ \hline
\end{tabular}$$
\medskip

Note that the dependence on the weights persists as	$N>>1$. Further numerical computations show that $G_N(\pm 2\pi i)$ converges to $Vol(M) \approx 2.02988321...$ as $N>\!>1$. 

\subsection{The complement of the knot $5_2$}\label{5_2} The ideal triangulation $T$ of $M $ with smallest complexity has three tetrahedra $\Delta^0$, $\Delta^1$, $\Delta^2$, and it has a branching $b$ that gives each tetrahedron the negative branching orientation. This triangulation is the one provided by Snappea; we keep the same ordering of tetrahedra. The branching $b$ is determined by the ordering it induces on the set of vertices of $\Delta^0$, which is obtained from that provided by Snappea by applying the permutation $(v_0,v_1,v_2,v_3)\mapsto (v_1,v_2,v_0,v_3)$. 

\subsubsection{The QH gluing variety} The edge relations of $(T,b)$ are $u_0u_2v_0^2w_1w_2=1$, $u_1u_2v_2v_1^2w_0w_1=1$, and $u_0u_1v_2w_0w_2=1$. By the tetrahedral relations they reduce to two independent relations, so the gluing variety of $(T,b)$ is the irreducible curve in $\mc_*^3$ with coordinates $(u_0,v_0,w_0)$ and defining equations
$$u_1 v_0^{-2} w_0 = 1\ ,\ u_2v_2^{-1} w_1=1.$$ 
The complete hyperbolic structure is realized by the point $(u_{hyp},v_{hyp},w_{hyp})$ given up to $8$ digits by (the imaginary parts are negative since all tetrahedra have the negative branching orientation)
$$\begin{array}{l} u_{hyp} =  0.21507987..-i\ 1.30714121... \\
v_{hyp} = 0.33764102...-i\ 0.56227951...\\
 w_{hyp} =  0.33764102...-i\ 0.56227951... \end{array}$$
Since $G(T,b)$ is irreducible and $(u_{hyp}^{-1},v_{hyp}^{-1},w_{hyp}^{-1})$ is a positive point, it is a rich variety.  Let $l$ and $m$ be the canonical longitude and the meridian of the knot, and let $(u_0,v_0,w_0)\in G(T,b)$. The dilation factors of the holonomy $\rho := \rho(u_0,v_0,w_0)$ are given by
$${\rm hol}_l(\rho) = u_1^2u_2^{-3}v_2v_0^{-2}w_1 = \frac{u_0^3(v_0-1)}{(1-u_0)^2(u_0-1)^3v_0^3(1-w_0)}$$
$${\rm hol}_m(\rho) = u_1v_0^{-1}= \frac{1}{(1-u_0)v_0}.$$ 
The tetrahedral and edges relations between q-shape parameters are (note that the exponents $-1$ on the left of the edge relations are due to the negative branching orientation of the tetrahedra)
$$\bu_0\bu_1\bu_2 =  e^{+\frac{\pi i}{N}}\ ,\ \bv_0\bv_1\bv_2 =  e^{+\frac{\pi i}{N}}\ ,\ \bw_0\bw_1\bw_2 =  e^{+\frac{\pi i}{N}}.$$
$$(\bu_0\bu_2\bv_0^{2}\bw_1\bw_2)^{-1} = e^{-\frac{2\pi i}{N}}\ ,\  (\bu_1\bu_2\bv_2\bv_1^2\bw_0\bw_1)^{-1}= e^{-\frac{2\pi i}{N}}\ ,\  (\bu_0\bu_1\bv_2\bw_0\bw_2)^{-1} = e^{-\frac{2\pi i}{N}}.$$
Together with the identity $(\bu_{i+1})^N = 1/(1-(\bu_i)^N)$ and the similar ones for the $(\bv_{i})^N$'s and the $(\bw_{i})^N$'s, these relations define $G_N(T,b)$. Using the tetrahedral relations it is easy to see that any of the edge relations is a consequence of the other two. By definition, the weights satisfy 
$$\bu_1^{2}\bu_2^{-3}\bv_2\bv_0^{-2}\bw_1  = e^{\frac{k(l)}{N}+\pi i h(l)}\ ,\ \bu_1\bv_0^{-1}= e^{\frac{k(m)}{N}+\pi i h(m)}.$$
As $M$ is a knot complement, the bulk weight $h$ is determined by $\kappa$, and similarly $h_c$ is determined by $k_c$. Moreover $h(l)=0$, so $h$ is non zero if and only if $h(m) = a_1+b_0= 1$ mod$(2)$. 

At the complete hyperbolic structure we have $\kappa(m)\in \pi i \mz$ and $\kappa(l)\in 2\pi i \mz$, and the tetrahedral relations, the edge relations, and the compatibility relations with $\kappa$ imply the following expressions of the edge colors $a_k$, $b_k$ and $c_k$, in terms of $a_0$, $a_1$ and $\kappa$: $a_2=2-a_0-a_1$, $b_2=2-b_0-b_1$, $c_2=2-c_0-c_1$, and
$$b_0 = a_1 - \frac{\kappa(m)}{\pi i}\ ,\ b_1=-4+2a_0+a_1+2\frac{\kappa(m)}{\pi i}-\frac{\kappa(l)}{2\pi i}$$
$$c_0 = -4+a_1-2\frac{\kappa(m)}{\pi i}\ ,\ c_1 = -a_0 - a_1 -\frac{\kappa(m)}{\pi i}+\frac{\kappa(l)}{2\pi i}.$$
If $(\bu,\bv,\bw) \in G_N(T,b)$ lies over a point in a small neighborhood of $(u_{hyp},v_{hyp},w_{hyp})\in G(T,b)$, then the same expressions are valid if one replaces $\kappa(\cdot)/\pi i$ with $(\kappa(\cdot )-\log({\rm hol}_{\cdot}(\rho)))/\pi i$. 

Similarly, the charges $a'_k$, $b'_k$ and $c'_k$ are given by
$$b_0' = a_1' - k_c(m)\ ,\ b_1'=-1+2a_0'+a_1'+2k_c(m)-k_c(l)/2$$
$$c_0' = a_1'-2k_c(m)\ ,\ c_1' = 1-a_0' - a_1' -k_c(m)+k_c(l)/2.$$
\subsubsection{The state sum formulas} For any $(\bu,\bv,\bw)\in G(T,b)$ and any charge $c$ on $T$ we have (here we denote by $\Rr_N(\Delta,\bw)$ the tensor $\Rr_{N,*_{b},c}(\bw_0,\bw_1))$ 
\begin{align*}\label{form1}
\Hh_{N}'(T,b,c)(\bu,\bv,\bw) & = \sum_{a,\ldots,f=0}^{N-1} \Rr_N(\Delta^0,\bu)^{d,a}_{b,c} \Rr_N(\Delta^1,\bv)^{f,b}_{e,a} \Rr_N(\Delta^2,\bw)^{c,e}_{d,f} \\
& = (\bu_0^{-a_1'}\bu_1^{a_0'}\bv_0^{-b_1'}\bv_1^{b_0'}\bw_0^{-c_1'}\bw_1^{c_0'})^{\frac{N-1}{2}}\frac{[\bu_0][\bv_0][\bw_0]g(1)^3}{g(\bu_0)g(\bv_0)g(\bw_0)} \\
 & \times \sum_{b,c,d=0}^{N-1} \frac{\zeta^{(b+c-d)(c+d)}}{\omega(\bu_0/\zeta,\bu_1^{-1}\vert b-d)\omega(\bv_0/\zeta,\bv_1^{-1}\vert d)\omega(\bw_0/\zeta,\bw_1^{-1}\vert d-c)}
%     & =_{2N} (\bu_0^{2}\bu_1)^{\frac{N-1}{2}}\frac{[\bv_0][\bw_0]g(1)^2}{g(\bu_1/\zeta)g(\bv_0)g(\bw_0)} \sum_{c,d=0}^{N-1} \frac{\zeta^{c(c+d)}\omega(\bu_1^{-1},\bu_0\vert -c-d)}{\omega(\bv_0/\zeta,(\bv_1)^{-1}\vert d)\omega(\bw_0/\zeta,(\bw_1)^{-1}\vert d-c)}
\end{align*}
%where the double sum comes from the identity (see eg. the proof of Proposition 5.3 in \cite{GT})
%$$\frac{[\bu_0]g(1)}{g(\bu_0)}\sum_{b=0}^{N-1} \frac{\zeta^{(b-d)(c+d)}}{\omega(\bu_0/\zeta ,\bu_1^{-1} \vert b-d)} =_{2N} \frac{(\bu_0^2\bu_1)}{g(\bu_1/\zeta)}^{\frac{N-1}{2}} \omega(\bu_1^{-1},\bu_0\vert -c-d). $$

The QH state sum $\Hh_{N}'(T,b,c)(\bu,\bv,\bw)$ computes the invariant $\Hh_{N}(V, \rho,\kappa, k_c,h,h_c)$ and defines a (non-constant) rational function on $G_N(T,b)$. As for the figure-eight knot complement, and using the same notations, at a point $(\bu,\bv,\bw)$ such that $\kappa(m^pl^q)= k_c(m^pl^q)=0$, $\Hh_{N}(V, \rho,\kappa, k_c,h,h_c)$ coincides with the $N$-th quantum hyperbolic invariant of the triple $(M',L,\rho')$ obtained by hyperbolic $(p,q)$-Dehn filling of $M$.

\subsubsection{Some numerical results}\label{num2} Keeping the same notations as in Section \ref{num1}, we present below some values of $\vert \Hh_{15}(V, \rho(A,B),\kappa,k_c, h,h_c) \vert$ up to 4 digits, and some values of $G_{121}$. For each of the shape parameter $u_0(A,B):=u_{hyp}+A+iB$ that we consider, with the exception of $u_0(0,0)$, there are two distinct points $(u_0(A,B),v_0(A,B),w_0(A,B))_{\pm}\in G(T,b)$ such that both $u_0(A,B)$, $v_0(A,B)$ and $w_0(A,B)$ have negative imaginary parts. We take take again $\kappa(m)=k_c(m)=0$, but now $k_c(l)=0$ and hence $\kappa(l)=k_f(l) =L([w;f])(l)$, where $w:=w_{\pm}(A,B)$ is determined by $(u_0(A,B),v_0(A,B),w_0(A,B))_{\pm}$, and $f$ is a (compatible) flattening such that, at the complete hyperbolic structure, $k^0_f:=L([w(0,0);f])$ takes the value indicated on the left of each row. 

$$\begin{tabular}{||c||c|c|c|c|c||} \hline &
  \rule{0cm}{0.4cm}{$(0,0)$}& $(-0.5,0.5)_+$ & $(0.5,0.2)_+$ & $(0.7,0.3)_+$ & $(0.7,0.3)_-$\\ \hline

\rule{0cm}{0.8cm} \raisebox{0.2cm}{$-6\pi i$} & \raisebox{0.2cm}{$428.2809...$} &
\raisebox{0.2cm}{$437.3709...$} & \raisebox{0.2cm}{$474.4262...$} & \raisebox{0.2cm}{$447.6316...$} & \raisebox{0.2cm}{$488.3855....$}\\ 
\hline

\rule{0cm}{0.8cm} \raisebox{0.2cm}{$-4\pi i$} & \raisebox{0.2cm}{$55.2674...$} &
\raisebox{0.2cm}{$70.5945...$} & \raisebox{0.2cm}{$65.7471...$} &
\raisebox{0.2cm}{$90.4874...$} & \raisebox{0.2cm}{$68.6357....$}\\ 
\hline

\rule{0cm}{0.8cm} \raisebox{0.2cm}{$-2\pi i$} & \raisebox{0.2cm}{$470.6170...$} &
\raisebox{0.2cm}{$487.7969...$} & \raisebox{0.2cm}{$518.2639...$} &
\raisebox{0.2cm}{$450.2483...$} & \raisebox{0.2cm}{$532.0331...$}\\ 
\hline

\rule{0cm}{0.8cm} \raisebox{0.2cm}{$0$} & \raisebox{0.2cm}{$71.9995...$} &
\raisebox{0.2cm}{$74.9737...$} & \raisebox{0.2cm}{$76.7206...$} &
\raisebox{0.2cm}{$96.6462...$} & \raisebox{0.2cm}{$77.9037....$}\\ 
\hline

\rule{0cm}{0.8cm} \raisebox{0.2cm}{$2\pi i$} & \raisebox{0.2cm}{$470.6170...$} &
\raisebox{0.2cm}{$494.8640...$} & \raisebox{0.2cm}{$515.6044...$} &
\raisebox{0.2cm}{$416.9562...$} & \raisebox{0.2cm}{$527.6381...$}\\ 
\hline

\rule{0cm}{0.8cm} \raisebox{0.2cm}{$4\pi i$} & \raisebox{0.2cm}{$55.2674...$} &
\raisebox{0.2cm}{$45.3442...$} & \raisebox{0.2cm}{$62.4615...$} &
\raisebox{0.2cm}{$109.4872...$} & \raisebox{0.2cm}{$65.4050...$}\\ 
\hline

\rule{0cm}{0.8cm} \raisebox{0.2cm}{$6\pi i$} & \raisebox{0.2cm}{$428.2809...$} &
\raisebox{0.2cm}{$456.9464...$} & \raisebox{0.2cm}{$466.5497...$} &
\raisebox{0.2cm}{$356.3560...$} & \raisebox{0.2cm}{$476.1700...$}\\ 
\hline

\end{tabular}$$
\medskip

$$\begin{tabular}{||c||c|c|c|c|c|c|c||} \hline &
  \rule{0cm}{0.4cm}{$-6\pi i$}& $-4\pi i$ & $-2\pi i$ & $0$ & $2\pi i$ %& $4\pi i$ %& $6\pi i$ 
  \\ \hline

\rule{0cm}{0.8cm} \raisebox{0.2cm}{$G_{121}$} & \raisebox{0.2cm}{$2.85352...$} &
\raisebox{0.2cm}{$0.20745...$} & \raisebox{0.2cm}{$2.83431...$} &
\raisebox{0.2cm}{$0.15541...$} & \raisebox{0.2cm}{$2.83431...$}%& \raisebox{0.2cm}{$0.20745...$}%& \raisebox{0.2cm}{$2.85352...$}
\\ \hline
\end{tabular}$$
\medskip

Further numerical computations show that $G_N(\pm 2\pi i)$ converges as $N>\!>1$ to $Vol(M) \approx 2.8281220883...$.

\subsection{The sister of the figure-eight knot complement} ($M=m003$ in Snappea's census).\label{sister} This cusped manifold is obtained by $(5, 1)$ Dehn surgery on a component of the Whitehead link (see \cite{Weeks}, \cite{MF}). The figure eight knot complement and its sister $M$ have the smallest volume of any orientable, cusped hyperbolic $3$-manifold. Hyperbolic Dehn surgery on $M$ induced by $(5, 1)$ and $(5, 2)$ Dehn surgeries on the Whitehead link yields the Fomenko-Matveev-Weeks manifold, which has the smallest volume of any closed orientable hyperbolic 3-manifold (approximately 0.9427...) \cite{GMM}. 

The Epstein-Penner decomposition of $M$ is the ideal triangulation $T$ with two tetrahedra considered in Example \ref{GT1} and \ref{moreF8S}. We have already described a weak branching $\tilde b$ on $T$ where both tetrahedra have positive branching orientations, and computed the charges, the flattenings at positive points of the gluing variety $G(T,\tilde b)$, and the QH state sums carried by $(T,\tilde b)$. 

\subsubsection{The QH gluing variety} The edge relations of $(T,\tilde b)$ are $u_0u_1^2v_0v_1^2 = 1$ and $u_0u_2^2v_0v_2^2 = 1$. By the tetrahedral relations they are equivalent, so the gluing variety $G(T,\tilde b)$ is the irreducible plane curve with coordinates $(u_0,v_0)$ and defining equation $u_0u_1^2v_0v_1^2 =1$.  The point $(e^{\frac{\pi i}{3}}, e^{\frac{\pi i}{3}})$ realizes the complete hyperbolic structure (thus $ \Delta^0$ and $\Delta^1$ are regular ideal tetrahedra, as in the case of the figure-eight knot complement). There is a basis $(l,m)$ of $\pi_1(\partial \bar M)$ where the dilation factors of the holonomy $\rho:= \rho(u_0,v_0)$ are given by 
\begin{equation}\label{eqhol}
{\rm hol}_l(\rho) = u_0u_1v_0^{-1}v_1^{-1} = \frac{u_0(1-v_0)}{(1-u_0)v_0}\ ,\ {\rm hol}_m(\rho) = u_1^{-2}v_2^{2} = \frac{(1-u_0)^2(v_0-1)^2}{v_0^2}.
\end{equation}
The tetrahedral and edges relations between q-shape parameters are
$$\bu_0\bu_1\bu_2 =  e^{-\frac{\pi i}{N}}\ ,\ \bv_0\bv_1\bv_2 =  e^{-\frac{\pi i}{N}}$$
$$\bu_0\bu_1^2\bv_0\bv_1^2 = e^{-\frac{2\pi i}{N}}\ ,\  \bu_0\bu_2^2\bv_0\bv_2^2 = e^{-\frac{2\pi i}{N}}.$$
Together with the identity $(\bu_{i+1})^N = 1/(1-(\bu_i)^N)$ and the similar one for the $(\bv_{i})^N$'s, these relations define $G_N(T,\tilde b)$. Using the tetrahedral relations it is easy to see that the two edge relations are equivalent. The reduced weights satisfy the compatibility relations
$$\bu_0\bu_1\bv_0^{-1}\bv_1^{-1}  = e^{\frac{k(l)}{N}+\pi i h(l)}\ ,\ \bu_{1}^{-2}\bv_2^{2} = e^{\frac{k(m)}{N}+\pi i h(m)}.$$
Here we note that $H_1(M;\mz) \cong \mz \oplus \mz/5\mz$, and that the curve $m$ or $l$ generates $H_1(M;\mz/2\mz) \cong \mz/2\mz$. So the bulk weight $h$ is determined by $\kappa$, and similarly $h_c$ is determined by $k_c$. Moreover $h(m)=0$ mod$(2)$, so $h$ is non zero if and only if $h(l) = a_0+a_1+b_0+b_1=a_2+b_2=1$ mod$(2)$.

%At the level of edge colors we have the tetrahedral and edge relations 
%$$\sum_{k=0,1,2} \log(u_k)+\pi i a_k = -\pi i\ ,\  \sum_{k=0,1,2} \log(v_k)+\pi i b_k= -\pi i$$
%$$\log(u_0)+2\log(u_1)+\log(v_0)+2\log(v_1)+ \pi i (a_0+2a_1+b_0+2b_1) = -2\pi i$$
%$$\log(u_0)+2\log(u_2)+\log(v_0)+2\log(v_2)+ \pi i (a_0+2a_2+b_0+2b_2) = -2\pi i,$$
%and the compatibility relations with the weight $\kappa$,  
%$$\log(u_0)+\log(u_1)-\log(v_0)-\log(v_1)+\pi i (a_0+a_1-b_0-b_1) = \kappa(l)$$
%$$-2\log(u_1)+2\log(v_2)+\pi i (-2a_1+2b_2) = \kappa(m).$$
At the complete hyperbolic structure we see that $\kappa(l)\in \pi i \mz$ and $\kappa(m)\in 2\pi i \mz$, and the tetrahedral relations, the edge relations, and the compatibility relations with $\kappa$ imply the following expressions of the edge colors $a_k$ and $b_k$ in terms of $a_1$ and $\kappa$:
\begin{equation}\label{standardf3bis} 
\begin{array}{c}  a_0 = \frac{\kappa(l)}{\pi i} - \frac{\kappa(m)}{2\pi i} -2a_1-2\ ,\ a_2 = a_1 -\frac{\kappa(l)}{\pi i} + \frac{\kappa(m)}{2\pi i}
\\ b_0 = \frac{\kappa(l)}{\pi i} - \frac{3}{2}\frac{\kappa(m)}{\pi i}-2a_1-2\ ,\ b_1 = -\frac{\kappa(l)}{\pi i} + \frac{\kappa(m)}{\pi i}+a_1\ ,\ b_2 = \frac{\kappa(m)}{2\pi i} + a_1.\end{array}
\end{equation}
For points $(\bu,\bv)\in G_N(T,\tilde b)$ lying over a small neighborhood of $(e^{\frac{\pi i}{3}}, e^{\frac{\pi i}{3}}) \in G(T,\tilde b)$, the same expressions are valid if one replaces $\kappa(\cdot)/\pi i$ with $(\kappa(\cdot )-\log({\rm hol}_{\cdot}(\rho)))/\pi i$. 
\subsubsection{The state sum formulas} For any $(\bu,\bv) \in G_N(T,\tilde b)$ and any charge $c$ on $T$ we have (again we denote by $\Rr_N(\Delta,\bw)$ the tensor $\Rr_{N,*_{b},c}(\bw_0,\bw_1)$)
\begin{align*}
\Hh_{N}'(T,\tilde b,c)(\bu,\bv) = & \sum_{i,j,k,l,I,J=0}^{N-1} 
 \Rr_N(\Delta^0,\bu)^{i,j}_{k,l}\Rr_N(\Delta^1,\bv)_{j,i}^{I,J}
(\Qq^{2})^l_I \Qq^k_J \\
& =  (\bu_0^{-a_1'}\bu_1^{a_0'}\bv_0^{-b_1'}\bv_1^{b_0'})^{\frac{N-1}{2}}
\frac{g(\bu_0)g(\bv_0)}{N\phi_Ng(1)^2}\\
& \quad \quad  \times \sum_{i,j,k,I=0}^{N-1} \zeta^{(i+k)(j+I)-ik+\frac{j^2+I^2}{2}} \omega(\bu_0,\bu_1^{-1} \vert i-k)\omega(\bv_0,\bv_1^{-1}\vert I-j)
%& =_{2N} \frac{(\bu_0)^{N-1}g(\bv_0)}{N\phi_Ng(1)g((\bu_1\zeta)^{-1})}\sum_{j,k,I=0}^{N-1} \frac{\zeta^{2k(j+I)+Ij-k^2}}{\omega((\bu_1\zeta)^{-1},\bu_0\vert j+I-k)\omega((\bv_0\zeta)^{-1},\bv_2\vert j-I)}
\end{align*}
%where the triple sum comes from the identities $\omega(\bv_0,\bv_1^{-1}\vert I-j) = \zeta^{-\frac{(I-j)^2}{2}}\omega((\bv_0\zeta)^{-1},\bv_2\vert j-I)^{-1}$ and
%$$\sum_{i=0}^{N-1} \zeta^{i(j+I-k)}\omega(\bu_0,\bu_1^{-1}\vert i-k) =_{2N} \frac{(\bu_0)^{N-1}g(1)}{g(\bu_0)g((\bu_1\zeta)^{-1})} \ \frac{\zeta^{k(j+I-k)} }{\omega((\bu_1\zeta)^{-1},\bu_0\vert j+I-k)}.$$
The QH state sum $\Hh_{N}'(T,\tilde b,c)(\bu,\bv)$ computes the invariant $\Hh_{N}(V, \rho,\kappa, k_c,h,h_c)$ and defines a (non-constant) rational function on $G_N(T,\tilde b)$. As for the previous examples, at a point $(\bu,\bv)$ such that $\kappa(m^pl^q)= k_c(m^pl^q)=0$, $\Hh_{N}(V, \rho,\kappa, k_c,h,h_c)$ coincides with the $N$-th quantum hyperbolic invariant of the triple $(M',L,\rho')$ obtained by hyperbolic $(p,q)$-Dehn filling of $M$.  
\subsubsection{Some numerical results}\label{num3} We keep the same notations as in Section \ref{num1} and we take again $\kappa(m)=k_c(m)=0$, but now both $u_0(A,B)$ and $v_0(A,B)$ have positive imaginary part, and the value of $\kappa(l)$ is given by $k_c(l)$ (on the left of each row), $w(A,B)$, and a flattening $f$ such that $L([w(0,0);f])=0$. The following table presents values of $\vert \Hh_{15}(V, \rho(A,B),\kappa,k_c, h,h_c) \vert$ up to 4 digits.

$$\begin{tabular}{||c||c|c|c|c|c|c||} \hline &
  \rule{0cm}{0.4cm}{$(0,0)$}& $(1,0)$ & $(-1,0)$ & $(0,1)$ & $(1,1)$ %& $(2,7)$ 
  \\ \hline 

\rule{0cm}{0.8cm} \raisebox{0.2cm}{$2$} & \raisebox{0.2cm}{$4.7755...$} & \raisebox{0.2cm}{$5.4995...$} & \raisebox{0.2cm}{$4.5346...$} & \raisebox{0.2cm}{$3.1943...$}  &
\raisebox{0.2cm}{$2.6295...$} %&  \raisebox{0.2cm}{$1.9024...$}
\\ 
\hline

\rule{0cm}{0.8cm} \raisebox{0.2cm}{$1$} & \raisebox{0.2cm}{$173.2621...$} & \raisebox{0.2cm}{$182.5736...$} & \raisebox{0.2cm}{$173.3850...$} & \raisebox{0.2cm}{$139.8353...$}  &
\raisebox{0.2cm}{$126.3516...$} %& \raisebox{0.2cm}{$101.7810...$} 
\\ \hline

\rule{0cm}{0.8cm} \raisebox{0.2cm}{$0$} &
\raisebox{0.2cm}{$0.2500...$} & \raisebox{0.2cm}{$0.6231...$} & \raisebox{0.2cm}{$0.3624...$} & \raisebox{0.2cm}{$0.2367...$}  & \raisebox{0.2cm}{$0.1498...$} %&  \raisebox{0.2cm}{$0.3701...$} 
\\ \hline

\rule{0cm}{0.8cm} \raisebox{0.2cm}{$-1$} &
\raisebox{0.2cm}{$173.2621...$} & \raisebox{0.2cm}{$114.3211...$} &
\raisebox{0.2cm}{$173.3850...$} & \raisebox{0.2cm}{$214.7305...$}  & \raisebox{0.2cm}{$236.1129...$} %&\raisebox{0.2cm}{$302.5352...$} 
\\ \hline

\rule{0cm}{0.8cm} \raisebox{0.2cm}{$-2$} &
\raisebox{0.2cm}{$4.7755...$} & \raisebox{0.2cm}{$6.8535...$} &
\raisebox{0.2cm}{$4.5346...$} & \raisebox{0.2cm}{$7.3011...$}  & \raisebox{0.2cm}{$9.3949...$} %&\raisebox{0.2cm}{$14.6569...$} 
\\ \hline
\end{tabular}$$

\end{document}